\documentclass[article]{amsart}
\RequirePackage{amsthm,amsmath,amsfonts,amssymb}

\usepackage[dvipdfmx]{graphicx}

\usepackage{bm}
\usepackage{mathrsfs}

\newtheorem{theorem}{Theorem}[section]
\newtheorem{proposition}{Proposition}[section]%

\newtheorem{lemma}{Lemma}[section]

\theoremstyle{remark}
	\newtheorem{remark}{Remark}[section]
\newtheorem{definition}{Definition}[section]
	\newtheorem{example}{Example}[section]

\newcommand\DN{\newcommand}

\DN\lref[1]{Lemma~\ref{#1}}
 \DN\tref[1]{Theorem~\ref{#1}}
 \DN\pref[1]{Proposition~\ref{#1}}
 \DN\sref[1]{Section~\ref{#1}}
 \DN\ssref[1]{Subsection~\ref{#1}}
 \DN\dref[1]{Definition~\ref{#1}}
 \DN\rref[1]{Remark~\ref{#1}} 
 \DN\eref[1]{Example~\ref{#1}}

\numberwithin{equation}{section}
\newcounter{Const} 
\numberwithin{Const}{section}
\DN\Ct{\refstepcounter{Const}c_{\theConst}}
\DN\cref[1]{c_{\ref{#1}}}
\DN\bs{\bigskip}

		\DN\ARc{\mathbb{A}_{\rR }^c}
		\DN\XB{(\mathbf{X},\mathbf{B})}
		\DN\SSSsde{\mathbf{S}_{\mathrm{sde}}}
		\DN\Bt{\mathcal{B}_t }
		\DN\zti{ 0 \le t < \infty }\DN\zzti{ 0 < t < \infty }
		\DN\WSN{W ^{\mathbb{N}}} 			
		\DN\FtB{$\{ \mathcal{F}_t \}$-Brownian motion }
		\DN\sigmaXms{\sigma \xm }
\DN\bbbXms{\mathit{b} \xm }
\DN\bbbXmshat{\mathit{b} \xmhat }
\DN\xmhat{_{\hat{\mathsf{X}}}^m}

\DN\SSsdemt{\mathsf{S}_{\mathrm{sde}}^m(t,\mathsf{X})}
\DN\SSsdemtw{\mathsf{S}_{\mathrm{sde}}^m(t,\ww )}
\DN\SSSsdemt{\mathbf{S}_{\mathrm{sde}}^m(t,\mathsf{X})}
\DN\SSSsdemtg{\mathbf{S}_{\mathrm{sde}}^m(t,\mathsf{X})}
\DN\SSSsdemtw{\mathbf{S}_{\mathrm{sde}}^m(t,\ww )}

	\DN\SSsde{\mathsf{S}_{\mathrm{sde}}}
	\DN\PPPm{\Pt ^m }
	\DN\Pt{\widetilde{P}} 
	\DN\lBlhatm{\mathbf{B}^m,\mathbf{X}^{m*}}
	\DN\WWdm{\WRdzm \ts \WRN }
	\DN\WRdzm{ W _{\mathbf{0}}^m}
	\DN\WRdm{ W ^m} 
\DN\Ehat{\mathcal{C}}
\DN\Ehatm{\Ehat^{m}}
\DN\Ehatmt{\Ehatm _t}
	\DN\WRNz{ W _{\mathbf{0}} (\RdN )}
	\DN\WRN{ W ^{\mathbb{N}}} 
	\DN\Btm{\mathcal{B}_t^m }
	\DN\uPs{under $ \Ps $}
\DN\Fm{F^m}
\DN\Fms{F_{\mathbf{s}}^m} 
\DN\Fmss{F_{\mathbf{s}}^m}

	\DN\OFPF{(\Omega ,\mathcal{F}, P ,\{ \mathcal{F}_t \} )}
	\DN\OFpsF{(\Omega ,\mathcal{F}, \Ps , \{ \mathcal{F}_t \} )}

\DN\sIn{$\mathbf{SIN}$}
\DN\iFc{$\mathbf{IFC}$}
\DN\iFcs{\As{$\mathbf{IFC}$}$ _{\mathbf{s}}$}
\DN\nbj{$\mathbf{NBJ}$}
\DN\muAC{$\mathbf{AC}$}
\DN\muTT{$\mathbf{TT}$}

	\DN\SSSsdeone{\mathsf{S}_{\mathrm{sde}}^{[1]}}
	\DN\xm{_{\mathsf{X}}^m}

\DN\labx{\lab (\xx )}
\DN\labi{\lab ^i}
\DN\labii{\lab ^{i+1}}
\DN\labix{\labi (\xx )}
\DN\labiix{\labii (\xx )}
\DN\Pls{{P}_{\lab (\xx )}^{\infty }}
\DN\Els{{E}_{\lab (\xx )}^{\infty }}
\DN\Os{\mathsf{upr}}
\DN\La{\mathsf{lwr}}
\DN\FBL{F(\XRp , \bRp ,\LRp )}
\DN\FBLL{F(\Xp , \mathbf{b}^{\perp } , 0 )}
\DN\nc{\nabla _j \psi ( \mathbf{X}_u^{\Rp , m } )}
\DN\BB{\mathsf{B}}
\DN\XRl{\mathbf{X}^{R,\La }}
\DN\XRs{\mathbf{X}^{\Rp }}

\DN\pL{\Big(} \DN\pR{\Big)_{i=1}^{\infty}}

\DN\CT{C([0,T];\Rtwo )^{\mathbb{N}}}
\DN\Cinfty{C([0,\infty);\Rtwo )^{\mathbb{N}}}
\DN\CiRdN{C([0,\infty);\Rd )^{\mathbb{N}}}
\DN\CRtwo{C([0,\infty);\Rtwo )}
\DN\CRtwoN{C([0,\infty);\Rtwo )^{\mathbb{N}}}

\DN\CTSS{C([0,T];\sSS )}
\DN\CiSS{C([0,\infty );\sSS )}

 \DN\rspRN{\rsp ,R \in \mathbb{N} }
\DN\Xsm{\mathbf{X}^{\perp , m}}
\DN\Xsn{\mathbf{X}^{\perp , n}}
\DN\si{\perp , i }
\DN\XRmbar{\overline{\mathbf{X}}^{\Rp , m}}
\DN\Xsmbar{\overline{\mathbf{X}}^{\perp ,m}}
\DN\Xmbar{\overline{\mathbf{X}}^{m}} 
\DN\XRidt{\mathsf{X}_t^{\RpiD }}
\DN\XRidu{\mathsf{X}_u^{\RpiD }}
\DN\XRid{\mathsf{X}^{\RpiD }}
\DN\XRjt{\XRj _t}
\DN\XRit{\XRi _t}
\DN\XRj{X^{\Rp ,j}}
	\DN\XRi{X^{\Rp ,i}}
\DN\Rp{\rR , \perp }
\DN\iD{i \diamond }
\DN\RpiD{\Rp , i \diamond }
\DN\RpkD{\Rp , k \diamond }
\DN\piD{\perp , i \diamond }
\DN\pkD{\perp , k \diamond }

\DN\sigR{\sigma _a^{i}}
\DN\sigmaA{\taui (\Xp )}
\DN\sigRpi{\taui (\X ^{\Rp })}

\DN\llL{\lim_{\rsp }}

\DN\ZQR{\mathcal{Z}^{\qQ , \rR }}	\DN\ZQ{\mathcal{Z}^{\qQ } }

\DN\vp{\varphi }
\DN\XXXsit{(\Xsit ,\Xsidt )}
\DN\XXXsiu{(\Xsiu ,\Xsidu )}
\DN\Xsidt{\mathsf{X}_t^{\piD }}
\DN\Xsidu{\mathsf{X}_u^{\piD }}
\DN\Xsit{\Xsi _t} 
\DN\Xsiu{\Xsi _u} 
\DN\XXXRit{(\XRit ,\XRidt )}
\DN\XXXRpiz{(\XRiz ,\XRidz )}
\DN\XXXRonet{(\XRt , \YYRt )}
\DN\XXXRoneu{(\XRu , \YYRu )}
\DN\XXXRone{( X ^{\Rp } , \YYR )}
\DN\XRt{X_t^{\rR ,\perp }} 
\DN\XRu{X_u^{\rR ,\perp }} 

\DN\YYRt{\YY _t^{\rR ,\perp }}
\DN\YYRu{\YY _u^{\rR ,\perp }}

 \DN\YYR{\YY ^{\rR ,\perp }}
\DN\XRz{X_0^{\rR ,\perp }} \DN\YYRz{\YY _0^{\rR ,\perp }}

\DN\XXXRiu{(\XRiu ,\XRidu )}
\DN\XXXRi{(\XRi ,\XRid )}
\DN\Sr{\sS _r}
\DN\Xsi{X^{\perp ,i}}
 \DN\ER{\E _R}
\DN\muRs{\mu ^{\Rp }}
\DN\XRiz{\XRi _0}
\DN\XRiu{\XRi _u}
\DN\XRju{\XRj _u}
\DN\XXXRiz{(\XRiz ,\XRidz )}
\DN\XXXRonez{(\XRz ,\YYRz )}
 \DN\XRidz{\mathsf{X}_0^{\RpiD }}

\DN\Qbbrsp{1_{\SQ }(x) ( \bbb (x,\yy ) - \bbb _{\rsp }(x,\yy )) }
\DN\bbrsp{\bbb - \bbb _{\rsp }}
\DN\bbrs{\bbb }
 \DN\Rsss{ \rR , \perp }
\DN\EE{E}
\DN\SQX{1_{\SQ }(X_t^{\Rp , i })}
\DN\SQXu{1_{\SQ }(X_u^{\Rp , i })}
\DN\SQXone{1_{\SQ }(\XRt )}
\DN\SQXoneu{1_{\SQ }(\XRu )}
\DN\SQXonez{1_{\SQ }(\XRz )}
\DN\SQx{1_{\SQ }(X_t^{\Rp , i })}

\DN\Rd{\mathbb{R} ^d}\DN\RdN{(\Rd )^{\mathbb{N} }}
\DN\Rdd{\mathbb{R} ^{\dd }}
\DN\Rtwo{\mathbb{R} ^2}
\DN\RtwoN{(\Rtwo )^{\mathbb{N}}}
\DN\RRtwo{\Rtwo }

\DN\ASN{\As{AC}, \As{SIN}, and \As{NBJ}}
\DN\IASN{\As{IFC}, \As{AC}, \As{SIN}, and \As{NBJ}}

\DN\limi[1]{\lim_{#1\to\infty}} 	\DN\limz[1]{\lim_{#1\to0}}
\DN\limsupi[1]{\limsup_{#1\to\infty}} 	\DN\limsupz[1]{\limsup_{#1\to0}}
\DN\liminfi[1]{\liminf_{#1\to\infty}} 	\DN\liminfz[1]{\liminf_{#1\to0}}
\DN\sumii[1]{\sum_{#1=1}^{\infty}}\DN\sumi[1]{\sum_{#1=0}^{\infty}}

\DN\map[3]{#1\!:\!#2\!\to\!#3} \DN\ot{ \otimes } \DN\ts{ \times }
\DN\PD[2]{\frac{\partial#1}{\partial#2}}
\DN\half{\frac{1}{2}}

\DN\oplusR{\oplus _{\rR }}

\DN\As[1]{$ ($\textbf{#1}$)$} \DN\Ass[1]{$ \{ $\textbf{#1}$\}$}
\DN\uL[1]{\underline{#1}}\DN\oL[1]{\overline{#1}}
\DN\PF{\begin{proof}} \DN\PFEND{\end{proof}} 
\DN\ssp{\smallskip} 
 \DN\ms{\medskip}

\DN\IRT{I_{\rR , T }}
\DN\IRRT{I_{\rR + 1 , T }}
 \DN\IQT{I_{\qQ , T }}\DN\IrT{I_{r , T }}
\DN\RT{\rR , T } \DN\QT{\qQ , T } \DN\QR{\qQ , \rR }
\DN\Rpk{\rR , \perp , k }
\DN\Rpi{\rR , \perp , i }
\DN\Rpj{\rR , \perp , j }
\DN\Rks{\rR , k , \sss }

\DN\LambdaR{\Lambda _{\rR } }
\DN\LambdaRm{\LambdaR ^m }

\DN\nuRR{\nu _{\rR ' ,\yy }}
\DN\nui[1]{ \mu (#1 \vert \Iinfty )}
\DN\nuR[1]{ \mu (#1 \vert \IR )}

\DN\muRRc{\mu _{\rR ,\rR }^c}
\DN\mutmR{\mut _{m , \rR}}
\DN\XAs{For $ \mutm $-a.e.\,$ \xx \in \SSm $}
\DN\Xas{for $ \mutm $-a.e.\,$ \xx \in \SSm $}
\DN\xSSmB{$ \xx \in \SSmB $}
\DN\xSSm{$ \xx \in \SSm $}
\DN\xasone{all $ \xx \in \SSone $}

\DN\Yas{$ \muxx $-a.s.\,$ \yy \in \SSmD $}
\DN\yas{$ \mutmD $-a.s.\,$ \yy $}

\DN\NNNR{\mathbf{N}_{\rR }^m }
\DN\nQ{\nn _{\qQ }}
\DN\nnn{\mathbf{n}}
\DN\nn{\mathit{n}}
\DN\0{\half \sumQ \Big\lvert \DsftQ \psilaP - \deltaPQ \Big\rvert ^2 }
\DN\1{\Big\lvert }
\DN\2{ \zeta \times \muz }
 \DN\3{ f_{ \rrr ,\sss }^m }
	\DN\4{\varpi _{\rR ,\yy }}
\DN\varpiR{\varpi _{\rR }}
\DN\varpiRR{\varpi _{\rR '}}
	\DN\TTRyn{\mathsf{T}_{\rR ,\yy }^{m , \nnn }}
	\DN\TTRynn{\mathsf{T}_{\rR ,\yy }^{m , \nnn '}}
	 \DN\UURy{\mathsf{U}_{\rR ,\yy }^m }
\DN\IR{\mathscr{I}_{\rR ,\yy } } \DN\IRR{\mathscr{I}_{\rR ',\yy } } \DN\Iinfty{\mathscr{I}_{\infty ,\yy } }
\DN\5{ 1 \le p \le d}
\DN\7{\varsigma _{\rR } } 
\DN\8{\half \sum_{ m =1}^{\infty} 1_{\UUpm } (\sss )} 

\DN\ssR{\sss _{\rR }}
\DN\ssRR{\sss _{\rR '}}

\DN\Upsone{\Upsilon _{\rR ,1}^{m}}
\DN\Upstwo{\Upsilon _{\rR ,2 } ^{ m , n }}
\DN\Upsthree{\Upsilon _{\rR ,3 } ^{ m , n }}

\DN\ve{\varepsilon }

\DN\kSQ{\sum_{k=1}^{\sss (\SQover ) } }
\DN\iSQ{\sum_{i=1}^{\sss (\SQover ) } }

\DN\BR{\mathcal{B}^{\rR } }
\DN\KR{\mathcal{K}^{\rR } }
\DN\LR{\mathcal{L}^{\rR } }
\DN\BRt{\mathcal{B}_t^{\rR } }
\DN\BQt{\mathcal{B}_t^{\qQ } }
\DN\KRt{\mathcal{K}_t^{\rR } }
\DN\LRt{\mathcal{L}_t^{\rR } }

\DN\maxT{\max_{0\le t \le T} } \DN\supT{\sup_{0\le t \le T} }

\DN\rsp{\mathsf{n}}
\DN\rspR{\rsp , \rR } 

\DN\intSS{\int_{\sSS }} 

\DN\dzf{\dz } 
\DN\dbf{\db }
\DN\Ps{P_{\mathbf{s}}}
\DN\PRs{P_{\rR ,\sss }^{\mu }}
\DN\PPRs{\mathbf{P}_{\rR , \lab (\sss ) }}
\DN\PPs{\mathbf{P}_{\lab (\sss ) }}

\DN\TT{\mathsf{T}}

\DN\aaa{\mathrm{a}} \DN\bbb{\mathsf{b}}
\DN\sss{\mathsf{s}}
\DN\sssiD{\mathsf{s}_i^{\diamond }}
\DN\ssskD{\mathsf{s}_k^{\diamond }}
\DN\mm{\sss } 

\DN\xx{\mathsf{x}}\DN\yy{\mathsf{y}} \DN\zz{\mathsf{z}}
\DN\xxx{\xx '}

\DN\dd{d}
\DN\dom{\mathscr{D}} 
\DN\domtrn{\dom ^{\mathrm{trn}}}

\DN\Xp{\X ^{\perp }}
\DN\XXp{\XX ^{\perp }}
\DN\X{\mathbf{X} }
\DN\XRp{ \X ^{\rR ,\perp } }
\DN\LRp{ \mathbf{L} ^{\Rp }} 
\DN\XR{\X ^{\rR }}
\DN\XX{\mathsf{X}} \DN\YY{\mathsf{Y}}
\DN\qQ{Q} 	\DN\rR{R} 	\DN\sS{S}
\DN\la{\lambda}

\DN\m{ { \sss ( \SRover ) } } 

\DN\XXRp{\mathsf{X}^{\Rp }}
\DN\XXRpD{\mathsf{X}^{\Rp , \diamond }}
\DN\XXpD{\mathsf{X}^{\perp , \diamond }}
\DN\XXpDk{\mathsf{X}^{\pkD }} 
\DN\XXRpDk{ \mathsf{X}^{ \RpkD } }
\DN\XXRpDi{ \mathsf{X}^{ \RpiD } }
\DN\XXpDi{\mathsf{X}^{ \piD } }

\DN\XRpi{X^{ \Rp , i}} 
\DN\XRpk{X^{ \Rp , k}} 
\DN\XRpD{\XRpk , \XXRpDk } 
\DN\XRpDi{\XRpi , \XXRpDi } 
\DN\XRpDu{X_u^{ R , \perp ,k} , \mathsf{X}_u^{ \RpkD } } 
\DN\XRpDui{X_u^{ R , \perp ,i} , \mathsf{X}_u^{ \RpiD } } 

\DN\XpDu{X_u^{ \perp ,k} , \mathsf{X}_u^{ \perp ,k\diamond} } 
\DN\XpDui{X_u^{ \perp , i } , \mathsf{X}_u^{ \perp , i \diamond} } 
\DN\XpDuk{X_u^{ \perp ,k } , \mathsf{X}_u^{ \perp , k \diamond} } 
\DN\XRpDt{X_t^{ R , \perp ,k} , \mathsf{X}_t^{ \RpkD } } 
\DN\XRpDti{X_t^{ R , \perp ,i} , \mathsf{X}_t^{ \RpiD } } 

\DN\XpDt{X_t^{ \perp ,k} , \mathsf{X}_t^{ \perp ,k\diamond} } 
\DN\XpDz{X_0^{ \perp ,k} , \mathsf{X}_0^{ \perp ,k\diamond} } 

\DN\XRDu{X_u^{ R ,i} , \mathsf{X}_u^{ R , \iD } } 
\DN\XDu{X_u^{ i} , \mathsf{X}_u^{ \iD } } 
\DN\XRDt{X_t^{ R , i} , \mathsf{X}_t^{ R , \iD } } 
\DN\XDt{X_t^{ i} , \mathsf{X}_t^{ \iD } } 
\DN\XDz{X_0^{ i} , \mathsf{X}_0^{ \iD } } 

\DN\nR{\mathbf{n}^{\rR }}
\DN\BBrspRk{\mathsf{b}_{\rsp }^{\Rpk }}

\DN\brspRp{\mathbf{b}_{\rsp }^{\rR , \perp }}
\DN\bRp{\mathbf{b}^{\rR , \perp }}
\DN\brsp{\mathbf{b}_{\rsp }^{\perp }}
\DN\bp{\mathbf{b}^{\perp }}

\DN\XRm{\mathbf{X}^{\rR , \perp , \infty }}
\DN\LRm{\mathbf{L}^{\rR , \perp }}
\DN\XXR{( \XRpDt )_{k\in\mathbb{N}}}

\DN\dib{\mathscr{D} _{\circ \mathrm{b}}}
\DN\dibone{C_{0}^{\infty} (\Rd )\ot \dib }
\DN\HRm{\mathcal{H}_{\rR }}
\DN\HxRm{\mathcal{H}_{ x , \rR }^m}

\DN\pirc{\pi _r^c}\DN\pisc{\pi _s^c}\DN\pitc{\pi _t^c}
\DN\pir{\pi _r}\DN\pis{\pi _s}

\DN\sumPQ{\sum_{ p , q = 1}^{\dd }} \DN\sumP{\sum \Pdd } \DN\sumQ{\sum \Qdd }
\DN\Pdd{_{ p =1 }^{\dd }}	\DN\Qdd{_{ q =1 }^{\dd }}

\DN\eP{\mathbf{e}_{ p }} 
\DN\eQ{\mathbf{e}_{ q }} 

	\DN\xP{x _{ p }} \DN\xQ{x _{ q }}
	\DN\deltaPQ{\delta _{ p , q }}

\DN\xiL{\xi _{L }}
\DN\xiLP{\xiL (\xP )}
\DN\xiLQ{\xiL (\xQ )}
\DN\xiLh{\xiL (\xP + h \eQ )}

	\DN\dYU{\oL{\dom }^{\YY }}
	\DN\DDDR{\DDD _{\rR }}
	\DN\DDDRR{\DDD _{\rR +1}}

\DN\TTTp{\mathbb{T}^{\perp }} \DN\UUUp{\mathbb{U}^{\perp }} 
\DN\VVVp{\mathbb{T}^{\gamma }} 
\DN\Dp{\DDD ^{\perp }}

\DN\psila{\psi _{\la }}
\DN\psilaa{\psi _{\la '}}
\DN\psilaP{\psi _{\la , p }}
\DN\psilaPP{\psi _{\la , p '}}
\DN\psilaPt{\tilde{\psi } _{\la , p }}
\DN\psilaaP{\psi _{\la ' , p }}
\DN\psilaaPt{\tilde{\psi }_{\la ' , p }}
\DN\psilaaQ{\psi _{\la ' , q }}
\DN\psilaQ{\psi _{\la , q }}

	\DN\psie{\epsilon \psi _{\epsilon ^2}}
	\DN\psieW{\epsilon \widetilde{\psi }_{\epsilon ^2}}
	\DN\psiP{\psi _{ p }}
	\DN\psiPP{\psi _{ p '}}

\DN\M{\mathcal{M}}
\DN\Mla{\M _{\la }}
\DN\Mlaa{\M _{\la '}}
\DN\Mlap{\M _{\la , p }}
\DN\Mlapp{\M _{\la , p '}}
\DN\Mlaq{\M _{\la , q }}

\DN\fQ{f _{\qQ }}
\DN\FQ{F _{\qQ }}

\DN\XY{X\YY }\DN\XYR{\XY , \rR }

\DN\capa{\mathrm{Cap}}

\DN\Nz{\{ 0 \} \cup \mathbb{N} } \DN\Ni{\mathbb{N} \cup \{ \infty \} } \DN\Nzi{\mathbb{N} \{ 0, \infty \} } 

\DN\ChiLp{\chi _{ L , p }}
\DN\ChiLpt{\tilde{\chi }_{ L , p }}

 \DN\chiL{\chi _{ L }}
\DN\chiLp{\chi _{ L , p }}
\DN\chiLpE{\chi _{ L , p ,\epsilon }}

	\DN\elaw{\stackrel{\mathrm{law}}{=}}	
	
	\DN\NR{\mathrm{N}_{\rR }}
	\DN\NRa{\mathrm{N}_{\rR , a}}
	 \DN\NTr{\mathrm{N}_{\TR }}
	 \DN\MTr{\mathrm{M}_{\TR }}
\DN\mucGi{\muz ( \cdot \vert \Gi ) }
\DN\muAGi{\muz ( A \vert \Gi ) }
\DN\muAGR{\muz ( A \vert \GR ) }

\DN\Ry{\rR , \yy } 
\DN\iy{\infty , \yy } 
		\DN\is{\infty , \sss } 

	\DN\muz{\mu _{0}}
	\DN\mux{\mu _{x}}
	\DN\muya{\mu _{\vartheta_{ \aaa } (\yy ) }}
\DN\muxxSSmD{\mu (\SSmD \Vert \xx )}
\DN\muxx{\mu (\cdot \Vert \xx )} 
\DN\muxxx{\mu (\cdot \Vert \xxx )} 
\DN\muyy{\mu (\cdot \Vert \yy )}
\DN\muy{\mu (\cdot \Vert \yy )}

\DN\muystar{\mu _{\yy }^* }

\DN\piRy{\piR (\yy )}\DN\piRcy{\piRc (\yy )}
\DN\circpiRc{ \circ (\piRc )^{-1}}

\DN\muRm{\mu (\cdot \cap \{ \sss (\SR ) \ge m \} ) \circpiRc }
\DN\muRmA{\mu ( \ARc \cap \{ \sss (\SR ) \ge m \} ) }
		\DN\muRc{\mu \circ (\piRc )^{-1}}

\DN\muR{\mu \circ \piR ^{-1}}
\DN\mutR{\mut \circ \piR ^{-1}}
\DN\mut{\widetilde{\mu }}
\DN\mutm{\mut _m}

\DN\mutmD{\mu _m ^{\diamond }}

\DN\mutone{\mut _1}
\DN\mutoneD{\mu _1^{\diamond }}

	\DN\muC{\sigma } 
	\DN\muiyC{\muC _{ \yy }} 
	\DN\muiyCC{\muC _{ \vartheta _{\aaa }(\yy ) }} 
	\DN\muiyCh{\muC _{\vartheta _{ h \eQ } (\yy )}}
	\DN\muyC{\muC _{\yy }}
\DN\muRyy{\mu _{\rR , \yy }}
\DN\muRss{\mu _{\rR , \sss }}
\DN\muRyym{\muRyy ^{ m }}
\DN\muRRyy{\mu _{\rR ', \yy }}
\DN\muRRRyy{\mu _{\overline{\rR }, \yy }}

	\DN\GR{\mathcal{G}_{\rR } }
	\DN\Gi{\mathcal{G}_{\infty } }
	\DN\MMR{\mathrm{M} _{\rR }}
	\DN\MRa{\mathrm{M} _{\rR , a }}
	\DN\piR{\pi _{\rR }} 

\DN\piRc{\piR ^{c}}
\DN\piRonec{\piR ^{[1],c}}
\DN\piTR{\pi _{\TR }}
\DN\piTQ{\pi _{\TQ }}
\DN\piRcyy{\piRc (\yy )}

\DN\Sm{\sS ^{ m }}
\DN\Rdm{(\Rd )^m}

\DN\SSmB{\SSm }
\DN\SSm{\sSS _{ m }} 

\DN\SSn{\sSS _{ n }}
\DN\SSone{\sSS _{1}}

\DN\SSz{\sSS _{ 0 }} \DN\SSi{\sSS _{ \infty }}
	\DN\SSmDW{\widetilde{\sSS}_{ m } ^{\diamond}}
	\DN\SSmD{\sSS_{ m } ^{\diamond}}
	\DN\SSnD{\sSS_{ n } ^{\diamond}}
	\DN\SSoneD{\sSS_{ 1 } ^{\diamond}}

\DN\SSzB{\SSz } 
\DN\SSoneB{\SSone } 
\DN\SSB{\sSS } 

	\DN\SSmDDk{ \{ \SSmD \}_{m=0}^{ k }}
	\DN\SSDDone{\{ \SSzD , \SSoneD \}}

\DN\SSyy{\sSS (\yy )}
\DN\SSfyy{\sSS ^{\mathrm{f}}(\yy )}
\DN\SSxx{\sSS (\xx )}

\DN\SSzD{\sSS_{ 0 } ^{\diamond}}
\DN\SSzDyy{\sSS_{ 0 } ^{\diamond } (\yy )}
\DN\SSzDfyy{\sSS_{ 0 } ^{\diamond } (\yy )}
\DN\SSzDxx{\sSS_{ 0 } ^{\diamond } (\xx )}
\DN\SSzDfxx{\sSS_{ 0 } ^{\diamond } (\xx )}
\DN\SSzDP{\sSS _{0}^{\diamond } ( }
\DN\SSoneDW{\widetilde{\sSS }_{1} ^{\diamond}}

\DN\ySSzD{\yy \prec \sss ,\, \sss \in \SSzD}

\DN\diaone{1 \diamond }
\DN\diai{i \diamond }

\DN\trn{\mathrm{trn}}
\DN\Dsft{D^{\trn} }
\DN\DsftP{D^{\trn}_{ p } }
\DN\DsftQ{D^{\trn}_{ q } }
\DN\dsft{\dom^{\trn}}

\DN\DXY{\mathbb{D}^{\XY }}
\DN\DY{\mathbb{D}^{\YY }}
\DN\nablaR{\nabla } %

\DN\CY{C^{\infty} (\Rd )\ot\dYL } 
\DN\dP{\mathscr{D}^{\perp } }

\DN\dzY{\dz ^{\YY }} 
\DN\dbY{\db ^{\YY }}
\DN\dY{\dom ^{\YY }}
\DN\dQbP{\dom _{ \rrrr }^{\perp }}
\DN\dbYp{\db ^{\YY ,\perp }}
\DN\dbP{\db ^{\perp }}

\DN\dbYRL{\uL{\mathscr{D}}_{\YY , \rR \bullet } }

\DN\dbXYR{\mathscr{D}_{\XYR \bullet } }

\DN\dbXYRL{\uL{\mathscr{D}}_{\XYR \bullet } }
\DN\dRXYL{\uL{\mathscr{D}}_{\XYR } }
\DN\dImu{\uL{\dom }^{\mu }}

\DN\dbXY{\db ^{\XY }}

\DN\dzmu{\dz ^{\mu }}

\DN\dbmug{\db ^{\mug }}

\DN\dzmuR{\dom _{\rR ,\circ } ^{\mu }}
\DN\dbmuR{\dbR ^{\mu }}
\DN\dbmu{\db ^{\mu }}

\DN\dRmu{\uL{\dom }_{\rR }^{\mu }}
\DN\dRmug{\uL{\dom }_{\rR }^{\mug }}
\DN\dRRmu{\uL{\dom }_{\rR +1}^{\mu }}

\DN\doneO{\oL{\dom }^{[1]}}
\DN\dL{\uL{\dom }}
\DN\dR{\dom _{\rR}}
\DN\dz{\dom _{\circ }}
\DN\dzR{\dom _{ \rR , \circ }} 
\DN\dbR{\dom _{ \rR , \bullet }} 
\DN\dRoneL{\dL _{\rR }^{[1]}}

\DN\dzone{\dz ^{[1]}}

\DN\db{\dom _{\bullet }}
\DN\dbb{\dom _{\bullet \mathrm{b}}}
\DN\dbone{\dom _{\bullet }^{[1]}}

\DN\dsm{\mathscr{D}_{\infty } } 

\DN\dRpL{\underline{\dom }_{\rr }^{\perp }} \DN\dpL{\underline{\dom }^{\perp }}

\DN\di{C_0^{\infty} (\Rd )}

 \DN\dYL{\uL{\dom }^{\YY }}

\DN\dtY{\dt ^{\YY }} 
\DN\dYtwo{\dtY ^{2}}
\DN\dt{\tilde{\mathscr{D}}}
 \DN\dzYt{{\dt }_{\YY \circ }}
 \DN\dYLt{{\dt }^{\YY }}

\DN\dXYL{\uL{\dom }^{\XY }}
\DN\dXYU{\oL{\dom }^{\XY }}

\DN\dboneL{{\dom }_{\bullet }^{[1]}}
\DN\dbXYL{{\dom }_{\bullet }^{\XY }}

\DN\doneL{\uL{\dom }^{[1]}}

\DN\dzYZ{\dz (\sSS \ts \mathbb{R}^2 )}
\DN\dbYZ{\db (\sSS \ts \mathbb{R}^2 )}

\DN\EXY{\mathscr{E}^{\XY }}
\DN\EXYtwo{\mathscr{E}^{\XY , 2}}
\DN\EXYstar{\mathscr{E}_{\XY \star}}

\DN\EXYRone{\mathscr{E}_{\XYR }} 
\DN\EXYR{\mathscr{E}^{\XY }}

\DN\xRm{_{ \x , \rR }^{m}}
\DN\xRonem{_{ \x , \rR }^{[1], m}}
\DN\Rsonem{_{ \rR , \sss }^{[1], m}}
\DN\Rsone{_{ \rR , \sss }^{[1]}}
\DN\ExRm{\mathscr{E}\xRm }
\DN\ExRonem{\mathscr{E}\xRonem }
\DN\ERsonem{\mathscr{E}\Rsonem }
\DN\xRsm{_{ \x , \rR , \sss }^{m}}

\DN\dYp{\mathscr{D}^{\YY , \perp } }

\DN\EY{\mathscr{E}^{\YY }}
\DN\EYL{\underline{\mathscr{E}}^{\YY }}
\DN\EYR{\mathscr{E}_{\rR }^{\YY }}
\DN\EYRL{\uL{\mathscr{E}}_{\rR }^{\YY }}

\DN\EYRR{\mathscr{E}_{\rR +1}^{\YY }}
\DN\EYRtwo{\mathscr{E}_{\rR }^{\YY }^2 }
\DN\EYone{\mathscr{E}^{\YY ,1}}
\DN\EYtwo{\mathscr{E}^{\YY ,2}}
\DN\EP{\mathscr{E}^{\perp } }
\DN\EYi{\mathscr{E}^{\YY , i }}

\DN\EYU{\underline{\mathscr{E}}^{\YY }}

\DN\ED{( \Emu , \overline{\dom }^{\mu })}
\DN\EDg{( \Emug , \overline{\dom }^{\mug })}
\DN\EDgP{( \EP , \overline{\dom }^{\perp })}

\DN\EYt{\tilde{\mathscr{E}}^{\YY }}
\DN\EYtone{\tilde{\mathscr{E}}^{\YY ,1}}
\DN\EYttwo{\tilde{\mathscr{E}}^{\YY ,2}}
\DN\EYtrgamma{\EYtr ^{2,\gamma }}

\DN\EYLt{ {\tilde{\mathscr{E}}}^{\YY }}
\DN\EYLtone{ {\tilde{\mathscr{E}}}^{\YY ,1}}
\DN\EYLttwo{ {\tilde{\mathscr{E}}}^{\YY ,2}}

\DN\Ek{\mathscr{E}^{[k]}}
\DN\Eone{\mathscr{E}^{[1]}}
\DN\EkL{\uL{\mathscr{E}}^{[k]}}
\DN\EoneL{{\mathscr{E}}^{[1]}}

\DN\ulab{\mathfrak{u} }
\DN\upath{\ulab _{\mathrm{path}}}
\DN\lab{\mathfrak{l}}
\DN\lpath{\mathfrak{l}_{\mathrm{path}}}

\DN\uu{u}
\DN\uN{\uu ^{N}}
\DN\vN{v^{N}}
\DN\vNsinfty{\int_{s\le \lvert x-y\rvert }\vN (x,y)dy }
\DN\vNs{\int_{\lvert x-y\rvert < s}\vN (x,y)dy }
\DN\vvv{v}
\DN\vsinfty{\int_{s \le \lvert x-y\rvert }\vvv (x,y)dy }
\DN\vs{\int_{\lvert x-y\rvert < s}\vvv (x,y)dy }

\DN\gin{\mathrm{Gin}}
\DN\mug{\mu _{\gin }}
\DN\mugRs{\mug ^{\rR \sss }}
\DN\mugRsone{\mug ^{\rR \sss ,[1]}}

\DN\LMz{L^{2}(\muz )}

\DN\mugz{\mu _{\gin , 0 }}
\DN\muga{\mu _{\gin ,\mathsf{a}}}
\DN\mugone{\mu _{\gin }^{[1]}}

\DN\mugxx{\mug (\cdot \Vert \xx )} 
\DN\mugyy{\mug (\cdot \Vert \yy )}
\DN\mugx{\mug (\cdot \Vert \xx )} 
\DN\mugy{\mug (\cdot \Vert \yy )}

 \DN\muone{\mu ^{[1]}}
 \DN\rg{\rho _{\gin }}

 \DN\Kg{\mathsf{K}_{\gin }}

	\DN\taui{\tau ^{i}}
	\DN\Ke{\mathsf{K}_{\epsilon }}

\DN\sSS{\mathsf{S}}

\DN\SSR{\sSS _{\rR }}
\DN\SSs{\sSS _{\mathrm{s}}}
\DN\SSsi{\sSS _{\mathrm{s,i}}}
\DN\WSsi{W (\SSsi )}
\DN\WSsiNE{W_{\mathrm{NE}} (\SSsi )}
\DN\ww{\mathsf{w}}

\DN\KKK{\mathsf{K}^{*}}
\DN\KKKRe {\KKK _{\rR , \epsilon }}
\DN\KKKR {\KKK _{\rR }} 
\DN\KKKI{\KKK _{\infty }} 
\DN\KKKIone{\mathsf{K} _{\infty }^{*,1}} 
\DN\KKKIe{\KKK _{\infty , \epsilon }} 

\DN\KKKIxx{\KKKI (\xx )}

\DN\SN{(\Rd )^{\mathbb{N}}}
\DN\SkS{\Sk \ts \sSS }

\DN\KKR{\mathsf{K}_{\epsilon , \rR }}
\DN\Sk{(\Rtwo )^{k}}

\DN\SRm{\SR ^m }
	\DN\SSRm{\SSR ^m }\DN\SSRn{\SSR ^n } 
\DN\SSRcz{ \sSS _{\rR ,c}^{0}}
\DN\SSRmm{\SSR ^{m-1} }
 \DN\SSRz{\SSR ^0 }
\DN\Se{\sS _{\epsilon }}

\DN\SR{\sS _{\rR }}
\DN\SRc{\SR ^c}
\DN\SQ{\sS _{Q}}
\DN\SRover{\overline{\sS }_{\rR }}
\DN\SQover{\overline{\sS }_{ Q }}
\DN\SRR{\sS _{\rR +1}}

\DN\Sone{\sS _{1}}
\DN\SRSS{ \SR \times \sSS }
\DN\RtwoSS{\mathbb{R} ^2 \ts \sSS }
\DN\RdSS{\Rd \ts \sSS }

\DN\rrrr{\rR }\DN\rrr{ \rR } \DN\rr{\rR }

\DN\TQ{ T _{ \qQ }}
\DN\TQQ{ T _{ \qQ +1}}

\DN\TR{ T _{\rR }} 
\DN\TRR{ T _{\rR +1}} 
 \DN\Tone{ T _1}

\DN\Up{\sS _{\rrr }} 	
\DN\Upp{\sS _{\rrr -1}}	
\DN\UUp{\sSS _{\rrr }} \DN\UUpm{\UUp ^m}

\DN\dlog{\mathsf{d}}
\DN\dgin{\dlog ^{\mug }}
\DN\dginone{\dlog ^{\mug ^1}}

\DN\dmu{\dlog ^{\mu }}

\DN\Lmz{L^{2}(\muz )}
\DN\Lmg{L^{2}(\mug )}
\DN\Lmgz{L^{2}(\mugz )}

\DN\Lone{L^1(\mu )}
\DN\Lm{L^{2}(\mu )}
\DN\Lmugone{L^{2}(\mugone )}
\DN\Lmugonep{L^{p}(\mugone )}
\DN\Lmug{L^{2}(\mug )}
\DN\Lmone{L^{2}(\muone )}

\DN\Llocone{L_{\mathrm{loc}}^{1}}
\DN\Llocmugone{L_{\mathrm{loc}}^{2} (\mugone )}
\DN\Llocmugonep{L_{\mathrm{loc}}^{p} (\mugone )}
\DN\LlocmugRsone {L_{\mathrm{loc}}^1 (\mugRsone )}

\DN\LlocmuRsone{L_{\mathrm{loc}}^1 ( \muRss ^{[1]})}
\DN\Llocmuone{L_{\mathrm{loc}}^1 ( \muone )}

\DN\E{\mathscr{E}}
\DN\Ep{\E ^{\perp }}
\DN\EpL{\underline{\E}^{\perp }}

\DN\Emug{\E ^{\mug }}
\DN\Emu{\E ^{\mu }}
\DN\ERmu{\Emu _{\rR }}
\DN\ERmug{\Emug _{\rR }}
\DN\ERRmu{\Emu _{\rR +1}}
\DN\Emuz{\E ^{\muz }}

\DN\DDD{\mathbb{D}}

\begin{document}
\title{Ginibre interacting Brownian motion in infinite dimensions is sub-diffusive }

\maketitle
\markboth{Hirofumi Osada}{Ginibre interacting Brownian motion in infinite dimensions is sub-diffusive}

\begin{center}
\textsf{Hirofumi Osada}\\
Chubu University \\
\texttt{osada@isc.chubu.ac.jp, osada@math.kyushu-u.ac.jp}
\end{center}

\begin{abstract}We prove that the tagged particles of infinitely many Brownian particles in $ \Rtwo $ interacting via a logarithmic (two-dimensional Coulomb) potential with inverse temperature $ \beta = 2 $ are sub-diffusive. The associated unlabeled diffusion is reversible with respect to the Ginibre random point field, and the dynamics are thus referred to as the Ginibre interacting Brownian motion. 
If the interacting Brownian particles have interaction potential $ \Psi $ of Ruelle class and the total system starts in a translation invariant equilibrium state, then the tagged particles are always diffusive 
if the dimension $ \dd $ of the space $ \mathbb{R}^{\dd } $ is greater than or equal to two. 
That is, the tagged particles are always non-degenerate under diffusive scaling. 
Our result is, therefore, contrary to known results. 
The Ginibre random point field has various levels of geometric rigidity. Our results reveal that the geometric property of infinite particle systems affects the dynamical property of the associated stochastic dynamics. 
\end{abstract}

\bs 


\medskip

	\tableofcontents

\section{Introduction}\label{s:1} 
We consider a system $ \X =( X ^i )_{i\in\mathbb{N} }$ of infinitely many Brownian particles moving in $ \Rd $ and interacting through a translation invariant, two-body potential 
$ \Psi (x) $. 
$ \X $ is described by the infinite-dimensional stochastic differential equation (ISDE) 
\begin{align}&\notag 
 X _t^i - X _0^i = B_t^i - \frac{\beta }{2} \int_0^t \sum_{j\not=i}^{\infty} \nabla \Psi ( X _u^i- X _u^j) du 
,\end{align}
where $ B^i $ ($ i\in\mathbb{N} $) denotes independent $ \dd $-dimensional Brownian motions and $ \beta \ge 0 $ is the inverse temperature, which is taken as a constant. 
The solution $ \X $ provides a description of the interacting Brownian motion \cite{lang.1,lang.2,fritz,tane.1}. 

The unlabeled process $ \XX =\{ \XX _t \}_{t\in[0,\infty)} $ 
associated with $ \X $ is given by 
\begin{align} &\notag 
\XX _t = \sum_{i\in\mathbb{N} } \delta_{ X _t^i}
,\end{align}
where $ \delta_a $ is the delta measure at $ a \in \Rd $ and 
$ \XX $ is a configuration-valued process by definition. 
 
We suppose that the unlabeled process $ \XX $ is reversible with respect to a translation invariant equilibrium state $ \mu ^{\Psi , \beta }$.
 In many cases, we expect the existence of such an equilibrium state. For example, if $ \Psi $ is a Ruelle-class potential (i.e., it is super-stable and regular in the sense of Ruelle), then the associated translation invariant canonical Gibbs measures exist. 
Here, super-stability is a condition that prevents infinitely many particles agglomerating in a bounded domain, and regularity means the integrability of interactions at infinity and therefore provides the Dobrushin--Lanford--Ruelle equation \cite{ruelle.2}. 

We investigate the tagged particles 
$ X ^i = \{ X _t^i \}_{t\in[0,\infty)} $ in the system. 
Although the total unlabeled system $ \XX $ is 
a $ \mu ^{\Psi , \beta }$-reversible Markov process, 
each tagged particle $ X ^i$ is a non-Markov process 
because the total system affects it in a complicated fashion. 
Applying the Kipnis--Varadhan theory, it can nevertheless be seen that 
the motion of each tagged particle always reverts to Brownian motion under diffusive scaling \cite{gp,De,o.inv1,o.inv2,tane.93}. 
That is, 
\begin{align} &\notag 
\limz{\epsilon } \epsilon X_{t/\epsilon ^2}^i = \sigma B_t 
.\end{align}
The constant may depend on the initial configuration, and we therefore introduce $ \alpha $, defined as the expectation of $ (1/2) \sigma ^2 $ with respect to the reduced Palm measure of the reversible measure $ \mu ^{\Psi , \beta }$. 
The constant matrix $ \alpha $ is called the self-diffusion matrix. 

Once such convergence of motion is established under this fairly general situation, it is natural and important to inquire about the positivity of the self-diffusion matrix $ \alpha $. 

Historically, there was a conjecture that $ \alpha = 0 $ for hard-core potentials and sufficiently large activities in multi-dimensional grand canonical Gibbs measures; cf. \cite{abl}. 
This conjecture seemed to be plausible because the presence of a hard core should suppress the motion of tagged particles. However the contrary was proved in \cite{o.p}. 
Indeed, $ \alpha $ is always positive definite if 
$ \dd \ge 2 $ and $ \Psi $ is a Ruelle-class potential corresponding to a hard core. 

In the setting of discrete spaces, the counterpart is a tagged particle problem of exclusion processes in $ \mathbb{Z}^{\dd } $. As for the simple exclusion processes, Kipnis--Varadhan \cite{KV} 
proved that $ \alpha $ is always positive definite except for the nearest neighborhood jump in one space dimension. Spohn \cite{spohn.ls} proved that $ \alpha $ is always positive definite for general exclusion processes with Ruelle-class potentials when $ \dd \ge 2 $. 

We note that the set of the Gibbs measures with Ruelle-class potentials is the standard class of random point fields in both continuous and discrete spaces. 
We hence consider, on good grounds, that it is reasonable to believe that self-diffusion matrices are always positive definite for $ \dd \ge 2 $. 
Nevertheless, we present the antithesis in the present paper. 


The Ginibre interacting Brownian motion $ \X = ( X ^i)_{i\in\mathbb{N}}$ 
is a system of infinite-many Brownian particles moving in $ \Rtwo $ and 
interacting via the two-dimensional Coulomb potential 
$ \Psi (x) = - \log \lvert x\rvert $ with inverse temperature $ \beta = 2 $. 
The stochastic dynamics $ \X = \{ \X _t \} $ 
are then described by the ISDE 
\begin{align}& \label{:10a} 
 X _t^i - X _0^i = B_t^i + 
\int_0^t 
\limi{\rR }
\sum_{\lvert X _u^i- X _u^j \rvert <\rR ,\, j\not=i} 
\frac{ X _u^i- X _u^j }{\lvert X _u^i- X _u^j \rvert ^2} du 
.\end{align}
The associated unlabeled process $ \XX $ is reversible with respect to the Ginibre random point field $ \mug $. By definition, $ \mug $ is a random point field on $ \Rtwo $ for which the $ n $-point correlation function $ \rg ^n $ with respect to the Lebesgue measure is given by 
\begin{align} &\notag 
\rg ^n (x_1,\ldots,x_n)= 
\det [\Kg (x_i,x_j)]_{1\le i,j \le n} \quad \text{ for each } n \in \mathbb{N} 
,\end{align} 
where $ \map{\Kg }{\mathbb{R} ^2 \ts \mathbb{R} ^2 }{\mathbb{C}}$ is the exponential kernel defined by 
\begin{align} &\notag 
\Kg (x,y) = \pi ^{-1} 
e^{-\frac{\lvert x\rvert ^{2}}{2}-\frac{\lvert y\rvert ^{2}}{2}}\cdot 
e^{x \bar{y}}
.\end{align}
Here, we identify $ \mathbb{R} ^2 $ as $ \mathbb{C}$ by the correspondence 
$ \mathbb{R} ^2 \ni x=(x_1,x_2)\mapsto x_1 + \sqrt{-1} x_2 \in \mathbb{C}$, and 
$ \bar{y}=y_1-\sqrt{-1} y_2 $ gives the complex conjugate of $ y $ under this identification. 

It is known that $ \mug $ is translation and rotation invariant. 
Furthermore, $ \mug $ is tail trivial \cite{o-o.tt,ly.18}. 
Let $ \sSS $ be the configuration space over $ \Rtwo $ defined by \eqref{:12s}. 
Let $ \sSS _*$ be a Borel subset of $ \sSS \backslash \{ \mathsf{0} \} $, where $ \mathsf{0}$ is the zero measure. 
A measurable map $ \map{\lab }{\sSS _*}{\RtwoN \cup\{ \sum_{m=1} ^{\infty} (\Rtwo )^m \} } $ 
called a label if $ \lab (\sss )=(\labi (\sss ))_{i}$ satisfies 
$ \sum_i \delta_{\labi (\sss )}= \sss $. 
Typically, we take a label $ \lab $ such that $ \lvert \labi (\sss )\rvert \le \lvert \labii (\sss )\rvert $ 
for all $ i $. 

For $ \mug \circ \lab ^{-1}$-a.s.\! $ \mathbf{s} = (s_i)_{i\in\mathbb{N} }$, 
\eqref{:10a} has a solution $ \X = ( X ^i)_{i\in\mathbb{N} }$, 
of which the unlabeled process $ \XX $ is $ \mug $-reversible \cite{o.isde}. 
The ISDE has a unique strong solution starting at $ \mug \circ \lab ^{-1}$-a.s.\,$ \mathbf{s} $ 
under a reasonable constraint. 
We refer to Section 7 in \cite{o-t.tail} and references therein for the existence and uniqueness of solutions of \eqref{:10a}. 

\begin{theorem}	\label{l:11}
Let $ \Ps $ be the distribution of the solution $ \X = ( X ^i)_{i\in\mathbb{N} }$ of \eqref{:10a} staring at $ \mathbf{s} \in (\Rtwo )^{\mathbb{N} }$. 
Then for each $ i \in \mathbb{N}$, under $ \Ps $ in $ \mug \circ \lab ^{-1}$-probability, 
\begin{align} &\notag 
\limz{\epsilon} \epsilon X _{\cdot /\epsilon ^2}^i = 0 \text{ weakly in $ C([0,\infty );\mathbb{R} ^2 ) $}
.\end{align}
\end{theorem}

\begin{remark}\label{r:13} The claim in \tref{l:11} means that 
for any bounded continuous function $ F $ on $ C([0,\infty ) ; \Rtwo ) $ and $ \kappa > 0 $, it holds that 
 for each $ i \in \mathbb{N} $, 
\begin{align}&\notag & 
\limz{\epsilon } \mug \circ \lab ^{-1} 
\Big(\Big\{ \mathbf{s} \in (\Rtwo )^{\mathbb{N} } ; 
\Big\lvert \int 
 F( \epsilon X _{ \cdot /\epsilon ^2}^i )) 
d\Ps - F (0) 
\Big\lvert > \kappa \Big\} 
 \Big) = 0 
.\end{align}
Here $ 0 = \{ 0_t \} $ of $ F (0) $ denotes the constant path with value $ 0 $. 
\end{remark}

Recently, it has become clear that the Ginibre random point field has 
various geometric rigidities, specifically 
a small variance property according to Shirai \cite{shirai}, 
the number rigidity according to Ghosh and Peres \cite{ghosh-peres}, and 
the dichotomy of reduced Palm measures \cite{o-s.abs}. 
These properties are different from those of the Poisson random point field and the Gibbs measure with a Ruelle-class potential, which have been extensively studied as the standard class of random point fields appearing in statistical physics. 

These geometric properties affect the dynamical properties. 
Indeed, from these rigidities, 
our theorem demonstrates that geometric rigidities yield dynamical rigidity in the sense of the sub-diffusivity of each tagged particle of the natural infinite-particle system given by \eqref{:10a}. 
%
%
%
%
%

The self-diffusion matrix $ \alpha [\mug ] $ is given by the solution to the Poisson equation of 
the quotient Dirichlet form on the configuration space (cf. \cite{o.inv1,o.inv2}). 
To prove $ \alpha [\mug ] = O $, we use the geometric rigidities of the Ginibre random point field we shall introduce in the following. 

The dichotomy of the reduced Palm measures of $ \mug $ was proved in \cite{o-s.abs}, and is 
the critical geometric rigidity that we use to prove \tref{l:11}. 
Let 
$ \mug (\cdot \Vert \xx ) $ be the reduced Palm measure of 
$ \mug $ conditioned at $ \xx = \sum_{i=1}^m \delta_{x_i} $ (see \eqref{:13x}). 
\begin{lemma}[{\cite[Theorem 1.1]{o-s.abs}}] \label{l:12}
Assume that $ \xx (\Rtwo )= m $ and $ \yy (\Rtwo ) = n $ for $ m , n \in \Nz $, where 
we take $\mugx = \mug $ if $ m = 0 $. The following then holds. 
\\\thetag{1} 
If $ m \not= n $, then $\mugx $ and $\mugy $ are singular relative to each other. 
\\\thetag{2} 
If $ m = n $, then $\mugx $ and $\mugy $ are mutually absolutely continuous. 
\end{lemma}

We find that $ \mugx $ is continuous in $ \mathbf{x}=(x_1,\ldots,x_m) \in (\Rtwo )^m$, where 
$ \xx = \sum_{i=1}^{ m } \delta_{x_i}$. 
That is, $ (\mathbb{R}^2)^{ m } \ni \mathbf{x} \mapsto \int _{\sSS } f d\mugx $ is continuous for any $ f \in C_b(\sSS )$. 
This follows from the explicit formula of the Radon--Nikodym density $ {d \mugxx }/{d \mugyy} $ in \lref{l:1X}. 
\begin{lemma}[{\cite[Theorems 1.2]{o-s.abs}}] \label{l:1X} 
Let $ \xx (\Rtwo ) =\yy (\Rtwo ) = m \in \mathbb{N} $. Then 
\begin{align} &\notag %
\frac{d \mugxx }{d \mugyy} = 
\frac{1}{\mathcal{Z}_{\xx , \yy } } \limi{\rR } \prod_{\lvert s_j\rvert < \rR }
\frac{\lvert \mathbf{x} - s_j\rvert ^2}{\lvert \mathbf{y} - s_j\rvert ^2}
,\end{align}
where $ \lvert \mathbf{x} - s_j\rvert = \prod_{i=1}^m \lvert x_i- s_j\rvert $. 
The normalization constant $ \mathcal{Z}_{\xx , \yy }$ is given by 
\begin{align} \notag 
\mathcal{Z}_{\xx , \yy } = &
\frac{\det [\Kg (x_i,x_j)]_{i,j=1}^m}{\det [\Kg (y_i,y_j)]_{i,j=1}^m}
\frac{\lvert \Delta (\mathbf{y})\rvert ^2}{\lvert \Delta (\mathbf{x})\rvert ^2}
.\end{align}
Here, 
$ \mathcal{Z}_{\xx , \yy } $ is the unique smooth function on $ (\Rtwo )^m \ts (\Rtwo )^m $ 
defined by continuity when the denominator has vanished. 
Furthermore, $ \Delta $ denotes the difference product for $ m \ge 2$ and 
$ \Delta (\mathbf{x}) = 1 $ for $ m = 1 $. 
\end{lemma}

Intuitively, the dichotomy in \lref{l:12} indicates the following phenomena. Suppose that we remove a finite unknown number $ m $ of particles $ \{ s_{i_1},\ldots,s_{i_m} \}$ from a sample point $ \sss = \sum_i \delta_{s_i}$ of $ \mug $. We then deduce the number $ m $ from information of $ \sss_{\diamond } := \sum_{i \in \mathbb{N} \backslash \{i_1,\ldots,i_m\}} \delta_{s_i}$. Such a structure is the same as periodic random point fields. Although a sample point $ \sss $ of $ \mug $ has enough randomness as seen from the simulation in Fig 1, we can infer the number of the removed particles exactly for $\mug $-a.s.\,$\sss $. 

For $ \dd = 1$, we proved that the non-collision of particles always implies sub-diffusivity \cite{o.nc}. (See also \cite{harris,spohn.2}.) 
Using the variational formula of the self-diffusion constant, 
the proof in \cite{o.nc} relies on the construction of a sequence of functions that reduces this constant to zero. This crucially uses the total order structure of non-collision particle systems in $ \mathbb{R} $, which is specific in one-dimension. 

A key point of the proof of \tref{l:11} is to construct such a sequence of functions without using the total order structure. Indeed, we shall use the above-mentioned geometric rigidity of the Ginibre random point field to accomplish this procedure. 

We shall present general theorems concerning the sub-diffusivity of the interacting Brownian motions for $ \dd \ge 2 $ and prove \tref{l:11} as a specific example of the general theorems 
(\tref{l:13} and \tref{l:14}).

Let $ \SR =\{ x \in \Rd ; \lvert x\rvert < \rR \} $. 
Let $\sSS $ be the configuration space over $\Rd $. 
\begin{align}\label{:12s}&
\sSS = \{ \sss = \sum _i \delta _{s_i} ; 
\sss ( \SR ) < \infty \text{ for all } \rR \in \mathbb{N} \} 
.\end{align}
We endow $\sSS $ with the vague topology, under which $\sSS $ is a Polish space. 
Let $ \mathcal{B}( \sSS ) $ be the Borel $ \sigma$-field of $ \sSS $. 
A probability measure $ \mu $ on $ (\sSS , \mathcal{B}( \sSS ) )$ is called a random point field and also a point process. 

Let $ \{ \vartheta_x \}_{x\in\Rd }$ be the translation operator on $ \sSS $ such that for 
$\sss =\sum_i\delta_{s_i}$, 
\begin{align}\label{:12t}&
 \vartheta_x (\sss ) = \sum_i \delta_{s_i - x}
.\end{align}
Then, $ \map{\vartheta_x }{\sSS }{\sSS }$ is a homeomorphism for each $ x \in \Rd $, and 
$ \sss \mapsto \vartheta_x (\sss )$ is a continuous function of $ x \in \Rd $ for each $ \sss \in \sSS $. 
Furthermore, $ (x , \sss ) \mapsto \vartheta_x (\sss )$ is continuous. 
A random point field $ \mu $ on $ \Rd $ is called translation invariant if 
\begin{align}\label{:12u}&
\mu = \mu \circ \vartheta_x ^{-1} \quad \text{ for all } x \in \Rd 
.\end{align}
We assume the following. 
\smallskip 

\noindent \As{A1} $ \mu $ is translation invariant and $ \mu (\{ \sss (\Rd ) = \infty \} ) = 1 $. 

\smallskip 

The translation invariance implies $ \mu (\{ \sss (\Rd ) = \infty \} ) = 1 $ if $ \mu $ is not a zero measure. Thus, the second assumption in \As{A1} yields no restriction in practice. 

A symmetric and locally integrable function $ \map{\rho ^n }{( \Rd ) ^n}{[0,\infty ) } $ 
is called the $ n $-point correlation function of $ \mu $ 
with respect to the Lebesgue measure if 
\begin{align} & \label{:12v} 
\int_{A_1^{k_1}\ts \cdots \ts A_m^{k_m}} \rho ^n (x_1,\ldots,x_n) dx_1\cdots dx_n 
 = \int _{\sSS } \prod _{i = 1}^{m} 
\frac{\sss (A_i) ! }
{(\sss (A_i) - k_i )!} \mu (d\sss )
 \end{align}
for any sequence of disjoint bounded measurable sets 
$ A_1,\ldots,A_m \in \mathcal{B}(\sS ) $ and a sequence of natural numbers 
$ k_1,\ldots,k_m $ satisfying $ k_1+\cdots + k_m = n $. 
When $ \sss (A_i) - k_i < 0$, according to our interpretation, 
${\sss (A_i) ! }/{(\sss (A_i) - k_i )!} = 0$ by convention. 
We make an assumption. 

\smallskip 

\noindent \As{A2} $ \mu $ has a locally bounded $ k $-point correlation function for each $ k \in \mathbb{N} $. 
\smallskip 

 We set the projections $ \map{ \piR , \piRc }{\mathsf{S} }{\mathsf{S} }$ such that 
 \begin{align} &\label{:12w}
\piR (\sss ) = \sss (\cdot \cap \SR ) , \quad \piRc (\sss ) = \sss (\cdot \cap \SR ^c) 
.\end{align}
For two measures $ \nu _1$ and $ \nu _2 $ on a measurable space 
$ (\Omega , \mathcal{B})$, we write $ \nu _1 \le \nu _2 $ 
if $ \nu _1(A)\le \nu _2(A)$ for all $ A\in\mathcal{B}$. 
Let $ \map{\Psi }{\Rd }{\mathbb{R} \cup \{ \infty \} }$ be a measurable function satisfying 
$ \Psi (x) = \Psi (- x )$. We take $ \Psi $ as an interaction potential of $ \mu $. 
Let $ \SSRm = \{ \sss \in \sSS ; \sss (\SR )= m \} $. 
We set $ \LambdaRm = \Lambda (\cdot \cap \SSRm )$, 
where $ \Lambda $ is the Poisson random point field whose intensity is the Lebesgue measure. 

For $ \mathsf{x} = \sum_i \delta_{x_i} \in \sSS $, let $ \HRm $ be the Hamiltonian on $ \SR $ such that 
\begin{align} & \label{:12z}
\HRm (\mathsf{x}) =  \sum_{ x_j, x_k \in \SR ,\, j < k } \Psi (x_j-x_k)
.\end{align}
\begin{definition}[{\cite{o.rm}}]\label{d:11} 
We say a random point field $ \mu $ on $ \Rd $ is a $ \Psi $-quasi-Gibbs measure with inverse temperature $ \beta \ge 0 $ if there exists a sequence of measures $ \{\mu_{\rR , k }^{m}\} $ on $ \sSS $ such that, for each 
$ \rR , m \in \mathbb{N}$, 
\begin{align}&\notag 
\mu_{\rR , k }^{m} \le \mu_{\rR , k +1}^{m} \text{ for all }k , \quad 
\limi{k} \mu_{\rR , k }^{m} = \mu (\cdot \cap \SSRm ) \text{ weakly}
\end{align}
and, for all $ \rR , k , m \in \mathbb{N}$ and $ \mu $-a.s.\! $ \sss $, the regular conditional measures 
\begin{align} & \label{:12x}
 \mu _{ \Rks }^{m} = 
\mu _{ \rR , k }^{m} (\,\piR (\xx ) \in \cdot \lvert \,\piRc (\mathsf{x}) = \piRc (\mathsf{\sss }) )
\end{align}
satisfy 
\begin{align} & \label{:12y}&
\cref{;10Q}^{-1} e^{- \beta \HRm (\mathsf{x}) } \LambdaRm (d\mathsf{x}) \le 
 \mu _{ \Rks }^{m} (d\mathsf{x}) \le 
\cref{;10Q} e^{- \beta \HRm (\mathsf{x}) } \LambdaRm (d\mathsf{x}) 
.\end{align}
 Here, $ \Ct \label{;10Q} $ 
is a positive constant depending on $ \beta $, $\rR $, $ k $, $ m $ and $ \piRc (\sss ) $. 
%
\end{definition}

\begin{remark}\label{r:11} 
The definition of quasi-Gibbs measure in \dref{d:11} is a particular case of that of \cite{o.rm}. 
Because we consider the tagged particle problem in the present paper, we adopt a more restrictive definition of quasi-Gibbs measures as above. 
\end{remark}


\noindent 
\As{A3} $ \mu $ is a $ \Psi $-quasi-Gibbs measure with inverse temperature $ \beta \ge 0 $ such that 
$\Psi ( x ) < \infty $ for $ x \ne 0 $ and that there exist an upper semi-continuous function $ \hat{\Psi} $ locally bounded from below and a constant $\Ct > 0 \label{;4z}$ satisfying 
$ \cref{;4z}^{-1} \hat{\Psi}(x) \le \Psi(x) \le \cref{;4z}\hat{\Psi}(x) $. 

\smallskip

Let $ \SSsi $ be the set consisting of infinite, single configurations such that 
\begin{align} &\notag 
 \SSsi = \{ \sss \in \sSS ; \sss (\{ x \} ) \in \{ 0,1 \} \text { for all }x \in \Rd 
 \, ,\, \sss (\Rd ) = \infty \} 
.\end{align}
We endow $ \SSsi $ with the vague topology. 
Let $ \WSsi $ be the set consisting of $ \SSsi $-valued continuous paths on $ [0,\infty)$. 
We write $ \ww = \{ \ww _t \} \in \WSsi $ as 
\begin{align}\label{:12a}&
 \ww _t = \sum_{i=1}^{\infty} \delta_{w^i(t)}, \quad w^i 
 \in C(I_i ; \Rd )
.\end{align}
Here,
$ I_i $ is an interval of the form $ [0,b_i)$ or $ (a_i,b_i)$, where $ 0 \le a_i < b_i \le \infty $. 
For $ \ww \in \WSsi $, the set $ \{ ( w ^i , I_i ) \}_{i \in \mathbb{N}} $ 
is uniquely determined (except labeling). 
For reader's convenience, we give a proof of this representation of $ \ww $ in \lref{l:Z1}. 

We write $ \ww = \{ ( w ^i , I_i ) \}_{i \in \mathbb{N}} $ if the representation of $ \ww $ 
is $ \{ ( w ^i , I_i ) \}_{i \in \mathbb{N}} $. 
Let 
\begin{align} \label{:12b}&
 \WSsiNE = \{\ww = \{ ( w ^i , I_i ) \}_{i \in \mathbb{N}} \in \WSsi ; I_i = [0,\infty ) \text{ for all $ i \in \mathbb{N}$}\} 
.\end{align}
For $ \ww \in \WSsiNE $, $ w^i \in C([0,\infty);\mathbb{R}^d)$ holds for all $ i \in \mathbb{N}$. 
Thus, $ \WSsiNE $ is the set consisting of non-exploding and non-entering paths. 
For a label $ \lab $ on $ \SSsi $, we have a unique map $ \map{\lpath }{ \WSsiNE }{ C([0,\infty);(\mathbb{R}^d)^{\mathbb{N}}) }$ such that 
\begin{align}\label{:12c}&
\lab (\ww _0) = \mathbf{w}_0, \quad \ww = \{ ( \sum_{i\in \mathbb{N}} \delta_{w^i(t)}) \} _{t\in [0,\infty) }
 \longmapsto \lpath (\ww ) = \mathbf{w} = ( w ^i)_{i \in \mathbb{N}}
.\end{align}
We give a proof of the construction of $ \lpath $ in \lref{l:Z2}. 

In \ssref{s:21}, we introduce the Dirichlet forms 
$ (\Emu , \dImu ) $ and $ \ED  $ on $ \Lm $. 
We call $ (\Emu , \dImu ) $ and $ \ED  $ the lower and upper Dirichlet forms, respectively. 
Such Dirichlet forms exist under \As{A2} and \As{A3}. 
In \ssref{s:4}, we introduce the perpendicular Dirichlet form $ (\EP , \dP ) $, which is a new Dirichlet form provided in the present paper and plays vital role in the proof of the main theorems. 

From \As{A2} and \As{A3}, we have a $ \mu $-reversible diffusion $ (P_{\sss }, \mathsf{X}_t)$ associated with $ \ED  $ on $ \Lm $ \cite{o.rm}. 
%
We set $ P _{\mu } (\cdot )= \int P_{\sss } (\cdot ) \mu (d\sss ) $. 
From \As{A1}--\As{A3} and $ \dd \ge 2 $, the $ \mu $-reversible diffusion $ (P_{\sss }, \mathsf{X}_t)$ 
has the non-explosion and non-collision properties (see Lemma 10.2 in \cite{o-t.tail}): 
\begin{align}\label{:12e}&
P _{\mu } ( \mathsf{X} \in \WSsiNE ) = 1
.\end{align}
From \eqref{:12c} and \eqref{:12e}, we have a continuous labeled process $ \X = \lpath (\mathsf{X})$. 
By construction, $ \X _0 = \lab (\mathsf{X}_0)$, $ \X = (X^i)_{i\in \mathbb{N} } $, and 
$ \mathsf{X} = \{ \mathsf{X}_t \}_{t\in [0,\infty)} $ is 
such that $ \mathsf{X}_t = \sum_{i\in \mathbb{N} } \delta_{X_t^i}$.

To prove sub-diffusivity, we introduce the new Dirichlet form $ (\EP , \dP )$ in \lref{l:41}. 
From \lref{l:41}\thetag{3} and \lref{l:21}, we have 
\begin{align}\label{:12g}&
(\EP , \dP ) \le ( \Emu , \dImu ) \le \ED  
.\end{align}
Here, for non-negative, symmetric bilinear forms $ (\mathscr{E}^i , \mathscr{D}^i ) $, $ i=1,2$, 
we write 
$ (\mathscr{E}^1, \mathscr{D}^1 ) \le (\mathscr{E}^2, \mathscr{D}^2 )$ if 
$ \mathscr{D}^1 \supset \mathscr{D}^2 $ and 
$ \mathscr{E}^1 (f,f) \le \mathscr{E}^2 (f,f) $ for all $ f \in \mathscr{D}^2 $. 

Taking \eqref{:12g} into account, we make an assumption. 

\smallskip 
\noindent 
\As{A4} $ (\EP , \dP ) = \ED  $. 
\smallskip 

We shall prove \As{A4} for the Ginibre random point field in \pref{l:K8}. 

For $ \xx , \sss \in \sSS $, we write $ \xx \prec \sss $ if 
$ \xx (\{ x \} ) \le \sss (\{ x \} )$ for all $ x \in \Rd $. 
For $ \sss , \xx \in \sSS $ such that $ \xx \prec \sss $, the difference $ \sss - \xx $ belongs to $ \sSS $. 
We set 
\begin{align}\label{:13r}&
\mathsf{A}- \yy = \{ \sss - \yy ; \yy \prec \sss , \sss \in \mathsf{A} \}
.\end{align}
By definition $ \mathsf{A} - \yy = \emptyset $ if no $ \sss \in \mathsf{A}$ satisfies 
$ \yy \prec\sss $. 

Let $ \SSm = \{ \sss \in \sSS ; \sss (\Rd ) = m \} $ for $ m \in \Nz $. 
By definition, $ \SSz $ consists of the zero measure. 
Let $ \ulab (\mathbf{x}) = \sum_i\delta_{x_i}$ for $ \mathbf{x} = (x_i)_i$. 
Then $ \SSm = \ulab ( (\mathbb{R}^d)^m) $. 
Let $ \check{\mu }^m $ be the $ m$th factorial moment measure of $ \mu $ such that 
\begin{align} \label{:13s}&
\check{\mu }^m (A_1\ts \cdots \ts A_m ) = \int_{A_1 \ts \cdots \ts A_m } \rho^m(x_1,\ldots ,x_m) 
dx_1 \cdots dx_m
\end{align}
for $\{ A_i \} $ such that $ A_i\cap A_j = \emptyset $ for $i\ne j $. 
Let $ \mutm = \check{\mu }^m \circ \ulab ^{-1 }$ be the measure supported on $ \SSm $. 
For $ \mutm $-a.e.\,$ \xx \in \SSm $, 
the reduced Palm measure $ \muxx $ exists and satisfies, by definition, 
for any $\mathsf{A} \in \mathcal{B}(\SSm )$, and $ \mathsf{B} \in \mathcal{B}(\sSS ) $, 
\begin{align}\label{:13x}& 	
\int_{\mathsf{A}} \mu ( \mathsf{B} \Vert \xx ) \mutm (d\xx ) 
 = 
\int_{\mathsf{A}} \mu ( \mathsf{B}+\xx \vert \xx \prec \sss ) \mutm (d \xx ) 
.\end{align}
%
For each $ \mathsf{B} \in \mathcal{B}(\sSS ) $, $ \mu (\mathsf{B} \Vert \xx )$ becomes a $ \mathcal{B}(\SSm )$-measurable function in $ \xx $. 

We write $ \nu_1 \ll \nu_2$ if $ \nu _1 $ is absolutely continuous with respect to $ \nu _2$, 
where $ \nu _1$ and $ \nu _2 $ are measures. 
We write $ \nu _1 \approx \nu _2 $ if $ \nu _1 \ll \nu _2 $ and $ \nu _2 \ll \nu _1 $. 

Taking the dichotomy of the Ginibre random point field in \lref{l:12} into account, we introduce the following concept. 
\begin{definition} \label{d:D1} 
We call $ \mu $ $ k $-decomposable with $ \SSmDDk $ if for $ 0 \le m \le k $ 
\begin{align} & \label{:13a}
 \SSmD \cap \SSnD = \emptyset \quad \text{ for }n \ne m , 0 \le n \le k 
,\\\label{:13b}&
\SSzD \subset \SSmD + \SSm 
,\\ \label{:13c}&
\SSmD \in \overline{\mathcal{B} (\sSS )}^{\muxx } \text{ and } \mu (\SSmD \Vert \xx ) = 1 
 \text{ for all } \xx \in \SSmB 
.\end{align}
We call $ \mu $ irreducibly $ k $-decomposable with $ \SSmDDk $ if, in addition, for $ 1 \le m \le k $ 
\begin{align}\label{:13f} & 
\muxx \approx \muxxx \text{ for all } \xx , \xxx \in \SSmB 
.\end{align}
\end{definition}

\begin{remark}\label{r:14}
\thetag{1} Let $ \mathsf{0}$ be the zero measure. 
Then $ \SSz = \{ \mathsf{0} \} $ and $ \mu (\cdot \Vert \mathsf{0}) = \mu $. 
Hence, we have $ \mut _0 = \mu (\cdot \Vert \mathsf{0}) = \mu $. 
Thus, \eqref{:13c} implies $ \mu ( \SSzD ) = 1 $. 
\\\thetag{2} From \eqref{:13a} and \eqref{:13c}, $ \mu $ and $ \muxx $ are singular relative to each other for $ \xx \in \SSm $ and $ \mu (\SSmD ) = 0 $ for $ m \ge 1 $. 
\end{remark}

\begin{example}[Ginibre random point field] \label{d:D3}
From \lref{l:12} and \lref{l:1X},  we easily see that the Ginibre random point field is irreducibly $ k$-decomposable for all $ k \in \mathbb{N}$. 
\end{example}

In \lref{l:X4}, we construct a reduced Palm measure $ \muyy $ for $ \yy $ such that $ \yy (\Rd ) = \infty $ if $ \mu $ is irreducibly decomposable in the sense of \dref{d:D1}. Such Palm measures are called dual reduced Palm measures. 
The concept of dual reduced Palm measures plays an important role of the proof of main theorems. 


\ssp 
\noindent 
\As{A5} \noindent $ \mu $ is irreducibly one-decomposable with $ \SSDDone $. 

\begin{theorem}	\label{l:13}
Assume $ \dd \ge 2 $. Assume \As{A1}--\As{A5}. 
Let $ \X = \lpath ( \mathsf{X}) = ( X ^i)_{i\in\mathbb{N} }$ be the labeled process defined after \eqref{:12e} with $ \X _0 = \lab (\sss)$. 
Then, for each $ i \in \mathbb{N}$, 
\begin{align} &\notag 
\limz{\epsilon} \epsilon X _{\cdot /\epsilon ^2}^i = 0 \text{ weakly in $ C([0,\infty );\Rd ) $ under $ P_{\sss } $ in $ \mu $-probability}
.\end{align}
\end{theorem}
\begin{remark}\label{r:15}\thetag{1}
In Theorems \ref{l:13}--\ref{l:15}, we assume $ d \ge 2$. This assumption is used only for the non-collision property of tagged particles. 
\\\thetag{2} 
 In \lref{l:X6}, we deduce \As{A5} from \As{A6} introduced below. 
 Thus, we obtain \As{A5} for the Ginibre random point field from \lref{l:12} and \lref{l:X6}. 
\end{remark}

We write 
$ \mux = \mu ( \cdot \Vert \delta_x ) $. 
From \As{A1} and \As{A2}, $ \mux $ exists for all $ x \in \Rd $. 
We can and do take 
$ \mux \circ \vartheta _x ^{-1} = \muz $ for all $ x \in \mathbb{R}^d $. 
We make an assumption. 

\noindent 
\As{A6} 
\thetag{1} $ \mu $ and $ \muz $ are singular relative to each other. \\
\thetag{2} $ \muz \approx \mux $ for all $ x \in \Rd $.

\begin{theorem}	\label{l:14}
Assume $ \dd \ge 2 $. Assume \As{A1}--\As{A4} and \As{A6}. 
We then obtain the same result as in \tref{l:13}. 
\end{theorem}

We prepare a set of notations for functions on $ \sSS $ following \cite{o.dfa}. 

Let $ \SSRm = \{ \sss \in \mathsf{S} ; \sss (\SR ) = m \} $ for $ \rR \in\mathbb{N} $ and 
$ m \in \{ 0 \}\cup \mathbb{N} $. 
Then 
\begin{align}& \label{:14z}
\sSS = \sum_{m=0}^{\infty} \SSRm 
.\end{align}
Let $ \SRm = \SR \ts \cdots \ts \SR $ be the $ m$-product of $ \SR $. We call $\mathbf{x}_{\rR }^m (\sss ) \in \SRm $ 
 an $\SRm $-coordinate of $\sss \in \SSRm $ if $ \piR (\sss )=\sum_{i=1}^m \delta_{x_{\rR }^{i}(\sss )}$, 
where $ \mathbf{x}_{\rR }^m (\sss )=(x_{\rR }^{i}(\sss ))_{i=1}^m $. 

For $f: \mathsf{S} \to \mathbb{R} $ and $\rR , m \in\mathbb{N} $, let 
$ f_{ \rR ,\sss }^m (\mathbf{x}) $ be the function satisfying 
\begin{align}\label{:14a}&
f_{ \rR , \cdot }^m (*) : \mathsf{S} \times \SRm \to \mathbb{R} 
\text{ such that } (\sss , \mathbf{x}) \mapsto f_{ \rR ,\sss }^m (\mathbf{x}) 
,\\
\label{:14b}&
\text{$f_{ \rR ,\sss }^m (\mathbf{x}) $ is permutation invariant in $ \mathbf{x}$ on $\SRm $ for each $\sss \in \SSRm $}
,\\\label{:14c}&
\text{$f_{ \rR ,\sss (1)}^m (\mathbf{x}) = f_{ \rR ,\sss (2)}^m (\mathbf{x}) $ 
if $ \piRc (\sss (1))=\piRc (\sss (2))$ for $\sss (1), \sss (2)\in \SSRm $}
,\\\label{:14d}&
\text{$f_{ \rR ,\sss }^m(\mathbf{x}_{\rR }^m (\sss ))=f(\sss )$ for $\sss \in \SSRm $}
,\\\label{:14e}&
\text{$f_{ \rR ,\sss }^m (\mathbf{x}) =0$ for $\sss \notin \SSRm $}
.\end{align}

Note that $ f_{ \rR ,\sss }^m$ is unique and 
$ f(\sss )=\sum_{m=0}^\infty f_{ \rR ,\sss }^m(\mathbf{x}_{\rR }^m (\sss ))$ for each $ \rR \in \mathbb{N} $ and $ \sss \in \sSS $. 
Here by convention, $ \mathbf{x}_{\rR }^0 (\sss ) = \emptyset $ for $ \sss \in \SSRz $. 
We see 
$ f_{ \rR ,\sss }^0 (\emptyset )= f ( \mathsf{0}) $ because $ \SSRz $ consists of the zero configuration 
$ \mathsf{0}$. 
%
The function $ f_{ \rR ,\sss }^0 $ is thus constant on $ \SSRz $. 
Although the $\SRm $-coordinate $ \mathbf{x}_{\rR }^m (\sss )$ of $\sss $ is not unique, $ f_{ \rR ,\sss }^m $ is well defined by \eqref{:14b}. 
For a bounded set $ A $, we set 
$ \mathbf{x}_{A }^m (\sss ) $ and 
$ f_{A ,\sss }^m (\mathbf{x} )$ similarly as above 
by replacing $ \SR $ by $ A $. 

A function $ \map{f}{\sSS }{\mathbb{R} }$ is called smooth if 
$f_{ \rR ,\sss }^m (\mathbf{x})$ is smooth in $ \mathbf{x} $ on $\SRm $ for all 
$\rR , m \in \mathbb{N} $, $\sss \in \mathsf{S} $, and local if $ f $ is 
$ \sigma[\pi _{\rR }]$-measurable for some $ \rR \in\mathbb{N} $. 
Let 
\begin{align} & \notag 
\db =\{  f ; 
\text{$ f $ is $ \mathcal{B}(\sSS ) $-measurable and smooth}\},  \quad 
\dz =\{f \in \db ; \text{$ f $ is local}\} 
,\\&  \label{:14g} 
\dbb = \{ f \in \db ; \text{$ f $ is bounded} \} ,\quad 
 \dib = \{ f \in \dz ; \text{$ f $ is bounded} \} 
.\end{align}
Let $ \muone $ be the one-Campbell measure of $ \mu $ such that 
\begin{align} \label{:14h}&
\muone (dx d\sss ) = \rho ^1 (x) \mux (d\sss )dx 
,\end{align}
where $ \rho ^1 $ is the one-point correlation function of $ \mu $ with respect to 
the Lebesgue measure and $ \mux = \mu ( \cdot \Vert \delta_x ) $. 
 $ \rho ^1$ exists and is constant by \As{A1} and \As{A2}. 


We now recall the concept of the logarithmic derivative of $ \mu $ from \cite{o.isde}. 
\begin{definition}[\cite{o.isde}] \label{d:12}
The logarithmic derivative $ \dmu $ of $\mu $ is an $ \Rd $-valued function such that 
$ \dmu \in \Llocone (\muone ) ^{\dd } $ and that, for all $\varphi \in \dibone $, 
\begin{align} & \label{:13i}&
\int _{\Rd \times \mathsf{S} } \dmu (x,\sss )\varphi (x,\sss ) 
\muone (dx d\sss ) = - \int _{\Rd \times \mathsf{S} } 
 \nabla_x \varphi (x,\sss ) \muone (dx d\sss ) 
.\end{align}
Here 
$ \Llocone (\muone ) = \bigcap_{\rR =1}^{\infty} L^1 (\muone _{\rR } )$ and 
$ \muone _{\rR } = \muone (\,\cdot\, \cap \{\SRSS \} )$. 
\end{definition}%

Once $ \dmu $ is calculated, we obtain the ISDE describing the labeled process 
$ \X = ( X ^i)_{i\in\mathbb{N} }$. 
Let $ \XX _t^{\diai } = \sum_{j \not= i } \delta_{X_t^j}$. 
We consider the ISDE 
\begin{align}\label{:14j}&
 X _t^i- X _0 ^i = B_t^i + \frac{1}{2} \int_0^t \dmu ( X _u^i, \XX _u^{\diai } ) du 
\  ( i \in \mathbb{N}) 
,\quad 
\X _0= \lab (\sss )
.\end{align}
Then, under \As{A2} and \As{A3}, \eqref{:14j}  has a weak solution for $ \mu $-a.s.\,$ \sss $ such that the associated unlabeled process $ \XX $ is a $ \mu $-reversible diffusion associated with the Dirichlet form $ \ED  $ 
on $ \Lm $ \cite{o.isde}. 
Under mild constraints, a weak solution of \eqref{:14j} is unique in law for $ \mu $-a.s.\,$ \sss $ (see \lref{l:74}, \cite{o-t.tail,k-o-t.ifc}). 


The logarithmic derivative 
$ \mathsf{d}^{\mug } $ of $ \mug $ is given by 
\begin{align} \label{:14l}&
\mathsf{d}^{\mug } (x , \sss ) = \limi{ \rR } \sum_{\lvert x -s_i \rvert < \rR } 
\frac{ 2 (x -s_i ) }{\lvert x -s_i \rvert ^2 }
\quad \text{ in } \Llocmugonep , 1\le p < 2 
.\end{align}
Hence, taking $ \mu = \mug $ in \eqref{:14j}, we obtain the ISDE \eqref{:10a} (see \cite{o.isde,o-t.tail}). 

We explain the idea of the proof of the main theorems (Theorems \ref{l:11}--\ref{l:15}). 
Let $ \muz $ be the reduced Palm measure conditioned at the origin. 
It is known that the self-diffusion matrix 
$ \alpha = (\alpha _{ p , q })_{ p , q =1}^d $ satisfies the variational formula 
\begin{align} &\label{:14L}
\alpha _{ p , p } = \inf 
 \Big\{ \int_{\sSS } 
\half \sumQ \Big\lvert \DsftQ f - \deltaPQ \Big\rvert ^2 + \DDD [ f ,f ] \, d\muz ; 
 f \in \dbY 
\Big\}
.\end{align}
Here, $ \Dsft = (\DsftP )_{ p = 1}^d $ is the generator of the translation $ \vartheta_x $ on $ \sSS $ defined by 
\begin{align}\label{:14m}&
\DsftP f (\sss ) = \limz{\epsilon} \frac{1}{\epsilon} 
\{ f (\vartheta_{\epsilon \eP } (\sss )) - f (\sss )\} 
,\end{align}
where $ \eP $ is the unit vector in the $ p $-direction, 
$ \deltaPQ $ is the Kronecker delta, 
 $ \DDD $ is the carr\'{e} du champ defined by \eqref{:21r}, 
and $ \dbY $ is the subset of $ \db $ given by \eqref{:71y}. 
In \cite{o.inv1,o.inv2}, $ \alpha _{ p , p }$ was given by replacing $ \dbY $ by $ \dzY $, 
where $ \dzY $ is given after \eqref{:71z}. Because we shall prove that the closures of 
$ \dbY $ and $ \dzY $ coincide in \lref{l:71}, $ \alpha _{ p , p }$ in \eqref{:14L} 
equals that given in \cite{o.inv1,o.inv2}. 
We set 
\begin{align}\notag &
\Dsft [f,g] = \half ( \Dsft f , \Dsft g )_{\Rd } ,\quad \EYtwo (f,g) = \int_{\sSS } \DDD [f,g] d\muz 
,\\  \label{:14M}  &
\EY (f,g) = \int_{\sSS } \Dsft [f,g] + \DDD [f,g] d\muz 
.\end{align}
In \lref{l:83}, we shall derive $ \alpha _{ p , p } = 0 $ from the following assumption. 
 
 \noindent 
\As{A7} 
For $ \5 $, 
we find an $ \EY $-Cauchy sequence $ \{ \ChiLp \} $ in $ \dbY $ such that 
\begin{align}\label{:14n}&
\limi{L} \Big\{ \int_{\sSS } 
\half \sumQ 
\Big\lvert \DsftQ \ChiLp - \deltaPQ \Big\rvert ^2 
\, d\muz + 
\EYtwo ( \ChiLp ,\ChiLp ) \Big\}
= 0 
.\end{align}

\begin{theorem}	\label{l:15}
\thetag{1} Let $ \dd \ge 2 $. Assume \As{A1}--\As{A5}. Then, \As{A7} holds. 
\\\thetag{2} 
Let $ \dd \ge 2 $. Assume \As{A1}--\As{A4} and \As{A7}. 
We then obtain the same result as in \tref{l:13}. 
\end{theorem}

To verify \As{A7}, we shall construct a sequence of functions $ \chiLp $ such that 
\begin{align}\label{:14o}&
\limi{L } \DsftQ \chiLp (\sss ) = \deltaPQ 
,\\\label{:14p}&
\limi{L } \DDD [ \chiLp , \chiLp ] (\sss ) = 0 
.\end{align}
At first glance, it is difficult to construct such a sequence of functions satisfying these two conditions. This is because the second condition suggests that the limit function is a constant, while the first condition states that it is not. To resolve this issue, we focus on 
the tail $ \sigma $-field $ \mathrm{Tail}(\sSS ) = \bigcap_{\rR \in \mathbb{N} } \sigma [\piRc ]$. 

We note that, from \eqref{:21p}--\eqref{:21r}, all tail measurable functions $ f $ satisfy 
\begin{align} &\notag 
 \DDD [f,f] = 0 
.\end{align}
We also remark that a tail measurable function $ f \in \db $ is not necessarily continuous under the vague topology. Indeed, it happens that, in general, 
\begin{align} & \label{:14r}
\limi{\rR } f (\piR (\sss )) \ne f (\sss )
\end{align} 
even if $f_{ \rR ,\sss }^m (\mathbf{x}) $ in \eqref{:14a}--\eqref{:14e} for $ f $ is constant for each $ \rR , m \in \mathbb{N} $ and $ \sss \in \sSS $. 

Let $ \SSi = \{ \sss \in \sSS ; \sss (\Rd ) = \infty \} $. 
If $ f $ is tail measurable, then $ f $ is constant on $ \sSS \backslash \SSi $. 
In contrast, $ f $ is not necessarily constant on $ \SSi $, as we see in \eqref{:14r}, even if $ f $ is tail measurable. 
We note $ \mu (\SSi ) = 1$ by \As{A1}. 
Thus, it may possible to construct a sequence of tail measurable functions $ f $ satisfying both \eqref{:14o} and \eqref{:14p}. 

The Ginibre random point field is tail trivial \cite{ly.18,o-o.tt}. 
Hence, a tail measurable function becomes a constant for $ \mug $-a.s.\,and thus does not satisfy \eqref{:14o}. Hence, we shall introduce the $ \sigma $-field $ \Gi $ in \eqref{:J1w}. 
This $ \sigma $-field is larger than $ \mathrm{Tail}(\sSS )$ and we can construct a sequence of $ \Gi $-measurable functions satisfying both \eqref{:14o} and \eqref{:14p}. 

From \As{A5}, we shall construct the function $ \chiLp ( \sss ) $ in \eqref{:Q1r}. 
Using \As{A1}, we deduce that $ \{\chiLp \} $ satisfies \eqref{:14o}. 
Furthermore, the function $ \chiLp ( \sss ) $ is $ \Gi $-measurable, and thus satisfies \eqref{:14p}.

The remainder of the paper is organized as follows. 
In \sref{s:2}, we introduce the three types of Dirichlet forms for the unlabeled dynamics. 
In \sref{s:III}, we present various diffusion processes related to the tagged particle problem and the associated Dirichlet forms: namely, one-labeled processes (\ssref{s:5}), tagged particle processes (\ssref{s:6}), and environment processes (\ssref{s:7}). 
In Subsection \ref{s:8}, we present a sufficient condition such that the limit self-diffusion matrix vanishes. 
In Subsection \ref{s:9}, we recall the Kipnis--Varadhan theory and prove an invariance principle of the additive functional of reversible diffusion processes. 
In Subsection \ref{s:X}, we introduce the concept of the dual reduced Palm measures conditioned at infinitely many particles for irreducibly decomposable random point fields. 
This concept is one of the main tools of our analysis. 
In Subsection \ref{s:J}, we introduce the mean-rigid $ \sigma $-field $ \Gi $, which yields the mean-rigid conditioning of random point fields. 
This $ \sigma $-field is also a key point of the proof of the main theorems. 
In Subsection \ref{s:Q}, we complete the proof of Theorems \ref{l:13}--\ref{l:15}. 
In Subsections \ref{s:K}--\ref{s:65}, we complete the proof of \tref{l:11}. 
\sref{s:VII} consists of appendices: 
In \ssref{s:Z}, we prove \eqref{:12a}. 
In \ssref{s:72}, we prepare the concept of strongly local, quasi-regular Dirichlet forms. 
In \ssref{s:73}, we quote the IFC condition and the result on the uniqueness of weak solutions of ISDEs. 
\section{Dirichlet forms of unlabeled dynamics} \label{s:2} 

In \sref{s:2}, we introduce three types of the Dirichlet forms describing the unlabeled dynamics: the perpendicular, lower, and upper Dirichlet forms. 
From \As{A4}, we deduce these three Dirichlet forms are the same in \lref{l:42}. 
This result yields the identity of the corresponding three Dirichlet forms of the environment process in \lref{l:71}. 

We say that a non-negative symmetric bilinear form $ (\E , \mathscr{D}_0 ) $ is closable on $ \Lm $ if $ \limi{n}\E (f_n , f_n ) = 0 $ for any $ \E $-Cauchy sequence 
$ f_n\in \mathscr{D}_0$ such that $ \limi{n} \| f_n \| _{\Lm } = 0 $. 
If $ (\E , \mathscr{D}_0) $ is closable on $ \Lm $, then there exists 
a closed extension of $ (\E , \mathscr{D}_0) $. 
The smallest closed extension $ (\E , \mathscr{D} ) $ of $ (\E , \mathscr{D}_0) $ is called the closure of $ (\E , \mathscr{D}_0) $. 
We refer to \cite{fot.2} for detail. 

For non-negative symmetric bilinear forms $ (\E _1, \mathscr{D}_1 ) $ and $ (\E _2, \mathscr{D}_2 ) $ on $ \Lm $, we say $ (\E _2, \mathscr{D}_2 ) $ is an extension of $ (\E _1, \mathscr{D}_1 ) $ if 
\begin{align}\label{:20a}&
 \mathscr{D}_1 \subset \mathscr{D}_2 , \quad \E _1 (f,f) = \E _2 (f,f) \quad \text{ for all }f \in \mathscr{D}_1 
.\end{align}
Suppose that $ (\E _2, \mathscr{D}_2 ) $ is an extension of $ (\E _1, \mathscr{D}_1 ) $. Then, $ (\E _2, \mathscr{D}_2 ) \le (\E _1, \mathscr{D}_1 ) $. The following simple fact will be used repeatedly in the present paper. 
\begin{lemma} \label{l:20}
Let $ (\E _1, \mathscr{D}_1 ) $ and $ (\E _2, \mathscr{D}_2 ) $ be non-negative symmetric bilinear forms on $ \Lm $. 
Let $ (\E _2, \mathscr{D}_2 ) $ be an extension of $ (\E _1, \mathscr{D}_1 ) $. Let $ (\E _2, \mathscr{D}_2 ) $ be closable on $ \Lm $. Then, $ (\E _1, \mathscr{D}_1 ) $ is closable on $ \Lm $. 
\end{lemma}
\begin{proof}
If $ \{ f_n \} $ is a Cauchy sequence of $ (\E _1, \mathscr{D}_1 ) $, then $ \{ f_n \} $ is a Cauchy sequence of $ (\E _2, \mathscr{D}_2 ) $ by \eqref{:20a}. 
Then $ \limi{n} \E _2 (f_n , f_n ) = 0 $ because of the closability of $ (\E _2, \mathscr{D}_2 ) $ on $ \Lm $. 
Hence, $ \limi{n} \E _1 (f_n , f_n ) = 0 $ from \eqref{:20a}. This completes the proof. 
\end{proof}
We say a closed non-negative symmetric bilinear form $ (\E , \mathscr{D} ) $ on $ \Lm $ is 
a symmetric Dirichlet form \cite{fot.2} if any $ u \in \dom $ satisfies 
\begin{align*}& 
v := \min \{ 1 , \max\{ 0 , u \} \} \in \dom \text{ and }
\E (v,v) \le \E (u,u)
.\end{align*}
For a symmetric Dirichlet form, there exists an associated symmetric Markovian $ L^2$-semi-group. 
In addition, if the Dirichlet form is strongly local and quasi-regular (see \ssref{s:72}) and the state space is homeomorphic to a complete separable metric space, then the associated symmetric diffusion process exists \cite{c-f}. In general, a Dirichlet form is not necessarily symmetric. In the present paper, a Dirichlet form means a symmetric Dirichlet form.

\subsection{Dirichlet forms associated with the unlabeled diffusions}\label{s:21}
In \ssref{s:21}, we prepare results for the Dirichlet forms associated with the unlabeled diffusions from \cite{k-o-t.udf,o.dfa,o.rm}. 
Let $ \SSRm = \{ \sss \in \sSS ; \sss (\SR )= m \} $ as before. 
Let $ \db $ be as in \eqref{:14g}. For $f,g \in \db $ and $ \sss =\sum_i \delta_{s_i}$, we set 
\begin{align} \label{:21p} &
\mathbb{D}_{\rR }^{m} [f,g] (\sss ) = 
1_{\SSRm } (\sss ) 
 \frac{1}{2} \sum_{ s_i\in \SR } 
 \big(
 \PD{}{s_{i}} f_{ \rR ,\sss }^m , \PD{}{s_{i}} g_{ \rR ,\sss }^m \big)_{\Rd } 
(\mathbf{x}_{\rR }^m (\sss )) 
. \end{align}
Here $ \PD{}{s_{i}} = (\PD{}{s_{i,1}},\ldots,\PD{}{s_{i,d}})$, 
$ f_{ \rR ,\sss }^m $ is as in \eqref{:14a}--\eqref{:14e} for $ f \in \db $, and 
$\mathbf{x}_{\rR }^m (\sss )$ is an $\SRm $-coordinate of $\sss $ introduced after \eqref{:14z}. 
Note that $\DDDR ^{m}[f,g] (\sss )$ is independent of the choice of 
the $\SRm $-coordinate $\mathbf{x}_{\rR }^m (\sss )$ and is well defined. Let 
\begin{align} & \label{:21q}
\DDDR =\sum_{m=1}^\infty \DDDR ^m 
.\end{align}
Then $ \DDDR [ f , f ] (\sss ) $ is non-decreasing in $ \rR $ for all $ f \in \db $ and $ \sss \in \sSS $. 
Hence, we set 
\begin{align}\label{:21r}&
\DDD [f , f ] ( \sss )= \limi{ \rR } \DDDR [ f , f ] (\sss ) \le \infty 
.\end{align}
We set the carr\'{e} du champs $ \DDD [ f , g ]$ by polarization. 
We set 
\begin{align}\label{:21d}&
\Emu (f,g) = \int_{\sSS } \DDD [f,g] d\mu , \quad 
\ERmu (f,g) = \int_{\sSS } \DDDR [f,g] d\mu 
,\\ & \notag
 \dbmu = \{ f \in \dbf ; \Emu (f,f) < \infty ,\, f \in \Lm \} 
,\\& \label{:21b} 
\dbmuR = \{ f \in \dbf ; \ERmu (f,f) < \infty ,\, f \in \Lm \} 
.\end{align}
Using the method in \cite{o.dfa,o.rm,k-o-t.udf}, we deduce from \As{A3} that
 $ (\ERmu ,\dbmuR )$ is closable on $ \Lm $. 
Hence, we denote by $ (\ERmu ,\dRmu ) $ its closure. 
Clearly, the sequence of the closed forms $ (\ERmu ,\dRmu ) $ is increasing in the sense that 
\begin{align} & \notag 
\ERmu (f,f) \le \ERRmu (f,f) \quad \text{ for all } f \in \dRRmu ,\quad 
\dRmu \supset \dRRmu 
.\end{align}
Let $ (\Emu , \dImu ) $ be the increasing limit of $ (\ERmu ,\dRmu ) $, $ \rR \in \mathbb{N}$, such that 
\begin{align}\notag &
\Emu (f,f) = \limi{\rR } \ERmu (f,f), 
\\ &\label{:21f}
 \dImu = \{f \in \bigcap_{\rR =1}^{\infty} \dRmu ; \limi{\rR } \ERmu (f,f) < \infty \} 
.\end{align}
Then $ (\Emu , \dImu ) $ is a closed form on $ \Lm $. 
We easily see that $ (\Emu , \dbmu )$ is closable on $ \Lm $ and its closure coincides with $ (\Emu , \dImu ) $. 

Let $ \dz $ be as in \eqref{:14g}. Clearly, $ \dz \subset \db $. We set 
\begin{align} \label{:21h}&
\dzmuR = \{ f \in \dzf ; \ERmu (f,f) < \infty , f \in \Lm ,\ 
 \text{$ f $ is $ \sigma [\piR ]$-measurable}
\} 
.\end{align}
Note that $ (\ERmu ,\dbmuR )$ is an extension of $ (\ERmu ,\dzmuR )$ and that $ (\ERmu ,\dbmuR )$ is closable on $ \Lm $. 
Hence, $ (\ERmu ,\dzmuR )$ is closable on $ \Lm $ by \lref{l:20}. 
Then we denote the closure of $ (\ERmu ,\dzmuR )$ on $ \Lm $ as $ ( \ERmu , \overline{\dom }_{\rR }^{\mu })$. 
By construction, $ (\ERmu , \dImu_{\rR } ) $ is an extension of $ ( \ERmu , \overline{\dom }_{\rR }^{\mu }) $. 
In particular, we have 
\begin{align} & \label{:21j} 
 ( \ERmu , \underline{\dom }_{\rR }^{\mu }) \le ( \ERmu , \overline{\dom }_{\rR }^{\mu }) 
.\end{align}
If $ f \in \dzmuR $, then $ f $ is $ \sigma [\piR ]$-measurable. Then 
$ \DDDR [f,f] = \DDDRR [f,f]$. Hence
\begin{align}\label{:21k}
 \ERmu (f,f) = \ERRmu (f,f) = \Emu (f,f) 
\text{ for all } f \in \overline{\dom }_{\rR }^{\mu } ,\quad 
 \overline{\dom }_{\rR }^{\mu } \subset \overline{\dom }_{\rR +1 }^{\mu }
.\end{align}
From \eqref{:21k}, $ ( \ERmu , \overline{\dom }_{\rR }^{\mu })$ is decreasing in $ \rR $. 
Let $ ( \Emu , \cup_{\rR \in \mathbb{N} }\overline{\dom }_{\rR }^{\mu }) $ be the decreasing limit. 
Note that $ (\Emu , \dImu ) $ is an extension of $ ( \Emu , \cup_{\rR \in \mathbb{N} }\overline{\dom }_{\rR }^{\mu }) $ and $ (\Emu , \dImu ) $ is a closed form on $ \Lm $. 
Hence, the decreasing limit $ ( \Emu , \cup_{\rR \in \mathbb{N} }\overline{\dom }_{\rR }^{\mu }) $ is closable on $ \Lm $ by \lref{l:20}. 

We denote the closure of $ ( \Emu , \cup_{\rR \in \mathbb{N} }\overline{\dom }_{\rR }^{\mu })$ on $ \Lm $ by $ \ED  $: 
\begin{align}\label{:21l}&
\overline{\dom }^{\mu } := \overline{\cup_{\rR \in \mathbb{N}}\overline{\dom }_{\rR }^{\mu }}^{\mu} 
.\end{align} 
\begin{lemma} \label{l:21} Assume \As{A1}--\As{A3}. 
Let $ \dz ^{\mu } = \{ f \in \dz ; \Emu ( f , f ) < \infty , f \in \Lm \} $. 
Then $ (\Emu , \dz ^{\mu } )$ is closable on $ \Lm $ and its closure 
$ ( \Emu , \overline{\dz ^{\mu }} ) $ satisfies 
$\overline{\dz ^{\mu }} = \overline{\dom }^{\mu } $ and 
\begin{align}\label{:21a}
 (\Emu , \dImu ) \le \ED  
\end{align}
\end{lemma}
\PF 
From \eqref{:21f}--\eqref{:21l}, it is clear that \eqref{:21a} holds. 
If $ f $ is $ \sigma [\piR ]$-measurable, then $ \DDD [f,f] = \DDDR [f,f]$. 
Hence from \eqref{:21h} and \eqref{:21l}, 
\begin{align}\label{:21n}&
 \dz ^{\mu } = \bigcup_{\rR \in \mathbb{N}}\dzmuR 
\subset \bigcup_{\rR \in \mathbb{N} }\overline{\dom }_{\rR }^{\mu } 
\subset \overline{\dom }^{\mu } 
.\end{align}
Then $ (\Emu , \dz ^{\mu } )$ is closable on $ \Lm $ by \lref{l:20} 
because $ ( \Emu , \overline{\dom }^{\mu } )$ is an extension of $ (\Emu , \dz ^{\mu } )$ from \eqref{:21n} and $ ( \Emu , \overline{\dom }^{\mu } )$ is closed on $ \Lm $. 
Thus, from this and \eqref{:21n}, 
\begin{align}\label{:21o}&
 \overline{\dz ^{\mu }} \subset \overline{\dom }^{\mu }
.\end{align}
Because $ ( \ERmu , \overline{\dom }_{\rR }^{\mu } )$ is the closure of $ (\ERmu ,\dzmuR )$, 
$ \dzmu \supset \dzmuR $, and \eqref{:21k} holds, we have 
$ \overline{\dz ^{\mu }} \supset \overline{\dom }_{\rR }^{\mu } $ for all $ \rR $. 
Hence, $ \overline{\dz ^{\mu }} \supset \cup_{\rR \in \mathbb{N}}\overline{\dom }_{\rR }^{\mu } $. 
From this and \eqref{:21l}, we deduce 
\begin{align}\label{:21O}&
 \overline{\dz ^{\mu }} \supset 
\overline{\cup_{\rR \in \mathbb{N}}\overline{\dom }_{\rR }^{\mu }}^{\mu} 
= \overline{\dom }^{\mu }
.\end{align}
We thus obtain 
 $ \overline{\dz ^{\mu }} = \overline{\dom }^{\mu } $ from \eqref{:21o} and \eqref{:21O}. 
\PFEND

It is known that $ \ED  $ on $ \Lm $ is a quasi-regular Dirichlet form and the associated $ \mu $-reversible diffusion $ (P_{\sss }, \mathsf{X}_t)$ satisfies \eqref{:12e} (see \cite{o.rm}, Lemma 10.2 in \cite{o-t.tail}, Lemma 2.5 in \cite{k-o-t.udf}). 
Hence, we have the labeled process $ \X = \lpath (\XX )$. 

\subsection{Perpendicular carr\'{e} du champs $ \Dp $} \label{s:3}
In \ssref{s:3}, we introduce the concept of the carr\'{e} du champ perpendicular to the generator of 
the translation operator $ \Dsft = ( \DsftP )\Pdd $ defined by \eqref{:14m}. 

For a set $ A \subset \Rd $ and $ \sss =\sum_i \delta_{s_i} \in \sSS $, we set 
\begin{align} &\notag 
\PD{}{\Gamma ( A )}= \sum_{s_i\in A } \PD{}{s_i} 
.\end{align}
%
Note that the orthogonal projection of $\partial / \partial s_i $ 
onto the subspace perpendicular to $\partial / \partial \Gamma ( A ) $ is then given by 
\begin{align} &\notag 
\PD{}{s_i} - \frac{1}{\sss ( A ) } \PD{}{\Gamma ( A )}
.\end{align}
We note that $ \sss ( A )$ becomes the number of particles in $ A $ for $ \sss = \sum_i \delta_{s_i}$. 

The orthogonality above is with respect to the inner product such that 
\begin{align} &\notag 
\Big( \PD{}{s_{i,p}} ,\, \PD{}{s_{j,q}} \Big) = \delta_{i,j}\delta_{p,q}
,\end{align}
where $ s_i = (s_{i,p})_{p=1}^d ,\, s_j = (s_{j,q})_{q=1}^d \in \Rd $. 
For each $ s_i \in A $, we have 
\begin{align}& \label{:31u} 
\Big( \PD{}{s_i} - \frac{1}{\sss ( A ) } \PD{}{\Gamma ( A )} ,\, 
 \frac{1}{\sss ( A ) } \PD{}{\Gamma ( A )}
\Big) = 0 
.\end{align}
Let $ A \subset B $. Then 
\begin{align}\label{:31v}&
\Big(
 \frac{1}{\sqrt{\mm ( A )}} 
 \PD{}{\Gamma ( A )} - \frac{\sqrt{\mm ( A )}}{\mm ( B )} \PD{}{\Gamma ( B ) }
, 
\frac{1}{\sqrt{\mm ( B )}} 
 \PD{}{\Gamma ( B )}
\Big) 
 = 0 
.\end{align}

Let $ T _1 = \sS _{1} $ and $ \TR = \SR \backslash \sS _{{\rrr -1}} $ for $ \rrr \ge 2 $. 
Note that $ \{ \TR \}_{\rrr \in \mathbb{N} } $ is a partition of $ \Rd $. 
For each $ \rrr \in \mathbb{N} $, let $ \{ \TT _{\rrr }^{ m } \}_{ m \in \{ 0 \} \cup \mathbb{N} } $ be the partition of $ \sSS $ such that $ \TT _{\rrr }^{ m } = \{ \sss \in \sSS ; \sss (\TR )= m \} $. 
We set the carr\'{e} du champs such that 
\begin{align} \label{:31y}
\TTTp _{\rrr } [f,f] (\sss ) &= 
\half 
\sum_{ m =1}^{\infty} 1_{\TT _{\rrr }^{ m } } (\sss )
\sum_{s_i\in \TR } 
\Big\lvert \Big(\PD{}{s_i} - 
\frac{1}{\mm (\TR ) }\PD{}{\Gamma (\TR )} \Big) f_{ \TR ,\sss }^m 
 \Big\rvert ^2 
,\\ \label{:31z}
\VVVp _{\rrr } [f,f ] (\sss ) &
= \half 
\sum_{ m =1}^{\infty} 1_{\TT _{\rrr }^{ m } } (\sss ) 
 \frac{ 1 }{{\mm ( \TR )}} \Big\lvert \PD{}{\Gamma ( \TR ) }
 f_{ \TR ,\sss }^m 
 \Big\rvert ^2 
\quad \text{ for $ f \in \db $}
.\end{align}
Here for $ f $, the functions $ f_{ \TR ,\sss }^m$ are given by \eqref{:14a}--\eqref{:14e}. 
\begin{lemma} \label{l:31} 
For each $ f \in \db $ and $ \rR \in \mathbb{N} $, 
\begin{align}\label{:31a} &
\DDDR [f,f] = \sum_{ Q = 1 }^{\rR } \TTTp _{ Q } [f,f] + 
 \sum_{ Q = 1 }^{\rR } \VVVp _{ Q } [f,f] 
.\end{align}
\end{lemma}
\PF 
For $ \sss = \sum_i \delta_{s_i}\in \TT _{ Q }^m $, we see from \eqref{:31u} 
\begin{align}\notag 
\sum_{s_i\in \TQ } 
\Big\lvert \PD{}{s_i}& f_{ \TQ ,\sss }^m 
 \Big\rvert ^2 = 
\sum_{s_i\in \TQ } 
\Big\lvert \Big(\PD{}{s_i} - 
\frac{1}{\mm (\TQ ) }\PD{}{\Gamma (\TQ )} \Big) f_{ \TQ ,\sss }^m 
+ \Big(
\frac{1}{\mm (\TQ ) }\PD{}{\Gamma (\TQ )} \Big) f_{ \TQ ,\sss }^m 
 \Big\rvert ^2 
\\ \label{:31f} =&
\sum_{s_i\in \TQ } 
\Big\lvert \Big(\PD{}{s_i} - 
\frac{1}{\mm (\TQ ) }\PD{}{\Gamma (\TQ )} \Big) f_{ \TQ ,\sss }^m 
 \Big\rvert ^2 
+ 
\frac{1}{\mm (\TQ ) }
 \Big\lvert 
\PD{}{\Gamma (\TQ )} f_{ \TQ ,\sss }^m 
 \Big\rvert ^2 
.\end{align}
Hence from \eqref{:31y}, \eqref{:31z}, and \eqref{:31f}, we deduce 
\begin{align}\notag &
\half \sum_{ m =1}^{\infty} 1_{\TT _{ Q }^{ m } } (\sss ) 
\sum_{s_i\in \TQ } \Big\lvert \PD{}{s_i} f_{ \TQ ,\sss }^m \Big\rvert ^2 
\\ \notag 
 = &
\half \sum_{ m =1}^{\infty} 1_{\TT _{ Q }^{ m } } (\sss ) 
\Big\{ 
\sum_{s_i\in \TQ } 
\Big\lvert \Big(\PD{}{s_i} - 
\frac{1}{\mm (\TQ ) }\PD{}{\Gamma (\TQ )} \Big) f_{ \TQ ,\sss }^m 
 \Big\rvert ^2 
+ 
\frac{1}{\mm (\TQ ) }
 \Big\lvert 
\PD{}{\Gamma (\TQ )} f_{ \TQ ,\sss }^m 
 \Big\rvert ^2 
\Big\}
\\ \notag =&
\TTTp _{ Q } [f,f] (\sss ) + \VVVp _{ Q } [f,f ] (\sss ) 
.\end{align}
Summing both sides over $ Q = 1,\ldots , \rR $, we obtain \eqref{:31a}. 
\PFEND

Let $ \UUpm = \{ \sss \in \sSS ; \sss (\Up )= m \}$. 
Note that $ \Upp \cup \TR = \Up $ and $ \Upp \cap \TR = \emptyset $ for $ \rrr \in \mathbb{N} $. 
We set, for $ \rrr \in \mathbb{N} $ such that $ 2 \le \rR $, 
\begin{align} \notag &
\UUUp _{\rrr } [f,f ] (\sss ) 
= \half 
\sum_{ m =1}^{\infty} 1_{\UUpm } (\sss )
\Big\lvert 
 \Big ( 
 \frac{1}{\sqrt{\mm ( \Upp )}} \PD{}{\Gamma ( \Upp )} -
 \frac{\sqrt{\mm ( \Upp )}}{\mm ( \Up )} \PD{}{\Gamma ( \Up ) }
\Big ) 
 f_{ \rrr ,\sss }^m 
 \Big\rvert ^2
\\ \label{:32v} & \quad \quad \quad + 
\half 
\sum_{ m =1}^{\infty} 1_{\UUpm } (\sss )
\Big\lvert 
 \Big ( 
 \frac{1}{\sqrt{\mm ( \TR )}} \PD{}{\Gamma ( \TR )} -
 \frac{\sqrt{\mm ( \TR )}}{\mm ( \Up )}
 \PD{}{\Gamma ( \Up ) }
\Big ) 
 f_{ \rrr ,\sss }^m 
 \Big\rvert ^2 
,\\ \label{:32x}&
\mathbb{U}_{\rrrr }^{\gamma }[f,f] (\sss ) = 
 \half \sum_{m=1}^{\infty} 1_{\UUpm } (\sss )
 \frac{ 1 }{\mm ( \Up )} \Big\lvert \PD{}{\Gamma ( \Up ) } f_{ \rrrr , \sss }^m (\sss )
 \Big\rvert ^2
\quad \text{ for $ \rR \in \mathbb{N} $}
.\end{align}

\begin{lemma} \label{l:2!}
Let $ \VVVp _{\rr }$ be as in \eqref{:31z}. For $ \rr \ge 2$, we have 
\begin{align}\label{:32c}&
 \mathbb{U}_{\rr -1}^{\gamma }[f,f] (\sss ) + \VVVp _{\rr } [f,f ] (\sss ) = 
\UUUp _{\rr } [f,f ] (\sss ) + \mathbb{U}_{\rr }^{\gamma }[f,f] (\sss ) 
.\end{align}
\end{lemma}
\begin{proof}
From \eqref{:32v} and \eqref{:32x}, we have 
\begin{align} \notag &
\UUUp _{\rrr } [f,f ] (\sss ) + \mathbb{U}_{\rrrr }^{\gamma }[f,f] (\sss ) 
\\ 
\notag 
& = \half 
\sum_{ m =1}^{\infty} 1_{\UUpm } (\sss ) 
\Big\lvert 
 \frac{1}{\sqrt{\mm ( \Upp )}} \PD{\3}{\Gamma ( \Upp )} -
 \frac{\sqrt{\mm ( \Upp )}}{\mm ( \Up )} \PD{\3}{\Gamma ( \Up ) }
 \Big\rvert ^2 %
\\ \notag & + 
\half 
\sum_{ m =1}^{\infty} 1_{\UUpm } (\sss ) 
\Big\lvert 
 \frac{1}{\sqrt{\mm ( \TR )}} \PD{\3}{\Gamma ( \TR )} -
 \frac{\sqrt{\mm ( \TR )}}{\mm ( \Up )}
 \PD{\3}{\Gamma ( \Up ) } 
 \Big\rvert ^2 
\\&+ \label{:!1a} 
 \half \sum_{m=1}^{\infty} 1_{\UUpm } (\sss ) 
 \frac{1}{\mm ( \Up )} 
\Big\lvert \PD{\3 }{\Gamma ( \Up ) } \1 ^2 
.\end{align}

Using $ \mm ( \Upp )+ \mm ( \TR ) = \mm (\SR )$ and 
$ \PD{}{\Gamma ( \Upp ) } + \PD{}{\Gamma ( \TR ) }= \PD{}{\Gamma ( \Up ) } $, 
we see 
\begin{align} \notag &
\Big\lvert 
 \frac{1}{\sqrt{\mm ( \Upp )}} \PD{\3}{\Gamma ( \Upp )} -
 \frac{\sqrt{\mm ( \Upp )}}{\mm ( \Up )} \PD{\3}{\Gamma ( \Up ) }
 \Big\rvert ^2 %
\\ \notag &\quad \quad 
+ 
\Big\lvert 
 \frac{1}{\sqrt{\mm ( \TR )}} \PD{\3}{\Gamma ( \TR )} -
 \frac{\sqrt{\mm ( \TR )}}{\mm ( \Up )}
 \PD{\3}{\Gamma ( \Up ) } 
 \Big\rvert ^2 
\\ \notag =&
\frac{1}{\mm ( \Upp )} \1\PD{\3 }{\Gamma ( \Upp )} \Big\rvert ^2 
+ 
 \frac{{\mm ( \Upp )}}{\mm ( \Up )^2} 
\1 \PD{\3 }{\Gamma ( \Up ) }\1 ^2 
 - \frac{2}{\mm ( \Up )} \Big( 
 \PD{\3 }{\Gamma ( \Upp )} 
, 
 \PD{\3 }{\Gamma ( \Up ) }\Big) _{\Rd }
\\ \notag & \quad \quad + 
 \frac{1}{{\mm ( \TR )}} \1 \PD{\3}{\Gamma ( \TR )} 
\1^2 
+
 \frac{\mm ( \TR )}{\mm ( \Up )^2}
\1 
 \PD{\3}{\Gamma ( \Up ) }
\1^2
- 
 \frac{2}{\mm ( \Up )} \Big( \PD{\3 }{\Gamma ( \TR )} , \PD{\3 }{\Gamma ( \Up ) }
\Big) _{\Rd }
\\ \notag =&
\frac{1}{\mm ( \Upp )} \1\PD{\3 }{\Gamma ( \Upp )} \Big\rvert ^2 
+ 
 \frac{1}{{\mm ( \TR )}} \1 \PD{\3}{\Gamma ( \TR )} 
\1^2 + 
 \frac{\mm ( \Upp )+ \mm ( \TR )}{\mm ( \Up )^2} \1 \PD{\3 }{\Gamma ( \Up ) }\1 ^2 
\\ \notag & \quad \quad 
- 
 \frac{2}{\mm ( \Up )} \Big( 
 \PD{\3 }{\Gamma ( \Upp )} + \PD{\3 }{\Gamma ( \TR )} , \PD{\3 }{\Gamma ( \Up ) }
\Big) _{\Rd }
\\ \label{:!1b} =&
\frac{1}{\mm ( \Upp )} \1\PD{\3 }{\Gamma ( \Upp )} \Big\rvert ^2 
+ 
 \frac{1}{{\mm ( \TR )}} \1 \PD{\3}{\Gamma ( \TR )} 
\1^2 
- 
 \frac{1}{\mm ( \Up )} 
\1 \PD{\3 }{\Gamma ( \Up ) } \1 ^2 
.\end{align}
Putting \eqref{:!1b} into \eqref{:!1a} and using \eqref{:31z} and \eqref{:32x}, we have 
\begin{align} &\notag 
\UUUp _{\rrr } [f,f ] (\sss ) + \mathbb{U}_{\rrrr }^{\gamma }[f,f] (\sss ) 
\\ \notag =&
\8 \frac{1}{\mm ( \Upp )} \1\PD{\3 }{\Gamma ( \Upp )} \Big\rvert ^2 
+
\8 
 \frac{1}{{\mm ( \TR )}} \1 \PD{\3}{\Gamma ( \TR )} 
\1^2 
\\ \notag 
=&
\8 \frac{1}{\mm ( \Upp )} \1 
\PD{ f_{ \rrr -1 ,\sss }^m }{\Gamma ( \Upp )} \Big\rvert ^2 
+ \VVVp _{\rr } [f,f ] (\sss ) 
\\ \notag 
=&
 \mathbb{U}_{\rr -1}^{\gamma }[f,f] (\sss ) + \VVVp _{\rr } [f,f ] (\sss ) 
.\end{align}
We have thus completed the proof of \eqref{:32c}. 
\end{proof}

Let $ \Dp _1 = \TTTp _{ 1 }$ and 
\begin{align} \label{:32w}
\Dp _{\rrrr } [f,f] (\sss ) &= 
\sum_{ Q = 1 }^{\rrrr } \TTTp _{ Q } [f,f] (\sss ) + 
\sum_{ Q = 2 }^{\rrrr } \UUUp _{ Q } [f,f] (\sss ) 
\quad \text{ for $ \rR \ge 2 $}
.\end{align}

\begin{lemma} \label{l:32}
For each $ f \in \db $ and $ \rR \in \mathbb{N} $, 
\begin{align}\label{:32a}&
\DDD _{\rrrr } [ f , f ] (\sss ) = \Dp _{\rrrr } [f,f] (\sss ) + \mathbb{U}_{\rrrr }^{\gamma }[f,f] (\sss ) 
.\end{align}
\end{lemma}

\PF 
Using \eqref{:32c}, we have 
\begin{align}\label{:32d}&
 \mathbb{U}_{1}^{\gamma }[f,f] (\sss ) +
\sum_{Q =2}^{\rR } \VVVp _{Q } [f,f ] (\sss ) = 
\sum_{Q =2}^{\rR } \UUUp _{Q } [f,f ] (\sss ) + \mathbb{U}_{\rR }^{\gamma }[f,f] (\sss ) 
.\end{align}
Note that $ \sS _1 = T _1 $ by definition. Hence from \eqref{:31z} and \eqref{:32x}, we have 
\begin{align}\label{:32i}&
 \VVVp _{1} [f,f ] (\sss ) = \mathbb{U}_{1}^{\gamma }[f,f] (\sss ) 
.\end{align}
Combining \eqref{:32d} and \eqref{:32i}, we deduce 
\begin{align}\label{:32j}&
\sum_{Q =1}^{\rR } \VVVp _{Q } [f,f ] (\sss ) =
\sum_{Q =2}^{\rR } \UUUp _{Q } [f,f ] (\sss ) + \mathbb{U}_{\rR }^{\gamma }[f,f] (\sss ) 
.\end{align}

Using \eqref{:31a}, \eqref{:32w}, and \eqref{:32j}, we obtain \eqref{:32a}. 
\PFEND

We now introduce the perpendicular carr\'{e} du champs $ \Dp $ such that 
\begin{align}\label{:33x}&
\Dp [f,f] (\sss ) = 
\sum_{ Q = 1 }^{\infty} \TTTp _{ Q } [f,f] (\sss ) + 
\sum_{ Q = 2 }^{\infty} \UUUp _{ Q } [f,f] (\sss ) 
.\end{align}
Then from \eqref{:32w}, $ \Dp _{\rrrr } [f,f] (\sss ) $ is increasing in $ \rrrr $ and satisfies
\begin{align}\label{:33y}&
\limi{\rrrr }\Dp _{\rrrr } [f,f] (\sss ) = \Dp [f,f] (\sss ) 
.\end{align}

\begin{lemma} \label{l:33}
 For each $ f \in \db $ and $ \sss \in \sSS $, 
\begin{align}\label{:33a}&
 \Dp [f,f] (\sss ) \le \DDD [ f , f ] (\sss ) 
.\end{align}
In particular, $ \Dp [f,f] (\sss ) < \infty $ holds for each $ f \in \dz $ and $ \sss \in \sSS $. 
\end{lemma}
\PF 
Note that both $ \Dp _{\rrrr } [f,f] (\sss ) $ and $ \DDD _{\rrrr } [ f , f ] (\sss ) $ are increasing in $ \rrrr $ and thus have the limits $ \Dp [f,f] (\sss ) $ and $ \DDD [ f , f ] (\sss ) $. Then, we obtain \eqref{:33a} from \eqref{:32a}. 

For each $ f \in \dz $, we there exists an $ \rR \in \mathbb{N} $ such that $ f $ is $ \sigma [\piR ]$-measurable. 
Because $ \DDD [f,f] (\sss ) = \DDDR [f,f] (\sss ) $ and $ \sss (\SR ) < \infty $, we have $ \DDD [f,f] (\sss ) = \DDDR [f,f] (\sss ) < \infty $. 
Hence from \eqref{:33a}, we obtain $ \Dp [f,f] (\sss ) < \infty $. 
\PFEND

\subsection{Identity of the perpendicular, lower, and upper Dirichlet forms} \label{s:4}
In \ssref{s:4}, we introduce the concept of the perpendicular Dirichlet form. 
We shall prove that the perpendicular, lower, and upper Dirichlet forms coincide under \As{A4}. 
We set for $ \sss = \sum_i \delta_{s_i}$ 
\begin{align} \notag 
& \Upsone = \Big\{ 
\upsilon = 
\PD{}{s_i} - \frac{1}{\mm (\TR ) }\PD{}{\Gamma (\TR )} 
 ; s_i \in \TR ,\, \sss \in \TT _{\rrr }^m 
\Big\} 
,\quad \rrr \in \mathbb{N} 
.\\ \notag &
\Upstwo = 
 \Big \{ \upsilon = 
 \frac{1}{\sqrt{\mm ( \Upp )}} 
 \PD{}{\Gamma ( \Upp )} - \frac{\sqrt{\mm ( \Upp )}}{\mm ( \Up )} \PD{}{\Gamma ( \Up ) }
 ; \sss \in \sSS _{\rrr - 1 }^m \cap \TT _{\rrr }^n 
 \Big\} 
,\\ \notag & 
\Upsthree = 
 \Big \{ \upsilon = 
 \frac{1}{\sqrt{\mm ( \TR )}} \PD{}{\Gamma ( \TR )} - \frac{\sqrt{\mm ( \TR )}}{\mm ( \Up )}
 \PD{}{\Gamma ( \Up ) } ; \sss \in \TT _{\rrr }^m \cap \UUp ^n 
 \Big\} 
,\quad \rrr \ge 2 
.\end{align}

If $ \upsilon \in \Upsone $, then $ \upsilon $ is a partial derivative on $ \TR $. 
To be precise, $ \upsilon \in \Upsone $ is a partial derivative on $ \TR ^{ m }$. 
We disregard $ m $ and simply call a partial derivative on $ \TR $ for convenience. 
Because $ \upsilon $ is a local operator, we can regard $ \upsilon $ as a partial derivative on any domain $ A $ including $ \TR $. 

For $\map{f}{\sSS }{\mathbb{R} }$ and $ \upsilon \in \Upsone $, 
we denote $ f \in \mathrm{Dom}(\upsilon ) $ if 
$ f_{\TR , \sss }^m (\mathbf{x})$ is in the domain of $ \upsilon $, 
where $ f_{\TR , \sss }^m $ is the representation of $ f $ defined by \eqref{:14a}--\eqref{:14e}. 

Let $ \ulab (\mathbf{x}) = \sum_i\delta_{x_i}$ for $ \mathbf{x} = (x_i)_i$ as before. 
We set $ \upsilon f (\sss ) = 0 $ for $ \sss \notin \TT _{\rr }^m $ and 
\begin{align} &\notag 
\upsilon f (\sss ) = \upsilon f_{\TR , \sss }^m(\mathbf{x}) 
\quad \text{ for $ \sss \in \TT _{\rr }^m $ such that }\pi_{\TR }(\sss ) = \ulab (\mathbf{x}) 
.\end{align}
For $\map{f}{\sSS }{\mathbb{R} }$ and $ \upsilon \in \Upsilon _{\rR , k}^{m,n}$, $ k=2,3 $, we define $ f \in \mathrm{Dom}(\upsilon )$ similarly. 
For $ \upsilon \in \Upstwo $, we set $ \upsilon f (\sss ) = 0 $ 
 for $ \sss \notin \sSS _{\rrr - 1 }^m \cap \TT _{\rrr }^n $ and 
\begin{align} &\notag 
\upsilon f (\sss ) = \upsilon f_{\SR , \sss }^{m+n}(\mathbf{x}) 
\quad \text{ for $ \sss \in \sSS _{\rrr - 1 }^m \cap \TT _{\rrr }^n $ 
such that }\pi_{\SR }(\sss ) = \ulab (\mathbf{x}) 
.\end{align}
For $ \upsilon \in \Upsthree $, we set 
$ \upsilon f (\sss ) = 0 $ 
 for $ \sss \notin \TT_{\rrr }^m \cap \SSR ^n $ and 
\begin{align} &\notag 
\upsilon f (\sss ) = \upsilon f_{\SR , \sss }^{n}(\mathbf{x}) 
\quad \text{ for $ \sss \in \TT_{\rrr }^m \cap \SSR ^n $ 
such that }\pi_{\SR }(\sss ) = \ulab (\mathbf{x}) 
.\end{align}

Let $ \DDD _{\upsilon} [f,g] = \half ( \upsilon f , \upsilon g )_{\Rd }$ be the carr\'{e} du champ 
generated by $ \upsilon $. 
Then, 
\begin{align} \label{:4r}&
\DDD ^{\perp } [ f , g ] = \sum_{\upsilon \in \Upsilon } \DDD _{\upsilon} [f,g] 
.\end{align}
Here $ \Dp $ is the perpendicular carr\'{e} du champ defined by \eqref{:33x} and $ \Upsilon $ is 
\begin{align} \label{:4p}&
\Upsilon = \{ \bigcup_{\rR , m \in \mathbb{N} } \Upsone \}\bigcup 
\{ \bigcup_{\rR \ge 2 ,\, m , n \in \mathbb{N} } \Upstwo \} \bigcup 
\{ \bigcup_{\rR \ge 2 ,\, m , n \in \mathbb{N} \atop m < n } \Upsthree \} 
,\end{align}
where the right-hand side is a disjoint union of $ \Upsone , \Upstwo $, and $\Upsthree $. 
Thus, $ \Upsilon $ is a collection of partial derivatives constituting $ \Dp $. 
We set 
\begin{align}\notag & 
 \EP _{\rrrr } (f,g)= \int_{\sSS } \Dp _{\rrrr } [f,g] d\mu 
, \quad 
\EP (f , f ) = \int_{\sSS } \DDD ^{\perp } [ f , f ] (\sss ) d \mu 
,\\ \label{:42q} &
 \dbP = \{ f \in \db ; \EP ( f , f ) < \infty , f \in \Lm 
\} 
.\end{align}
Here $ \Dp _{\rrrr }$ is as in \eqref{:32w}. In \eqref{:42q}, we suppress $ \mu $ from the notation. 
\begin{lemma} \label{l:41} 
Assume \As{A2} and \As{A3}. Then the following hold. 
\\\thetag{1} 
$ ( \EP _{\rrrr }, \dbP ) $ and $ ( \EP , \dbP ) $ are closable on $ \Lm $. 
\\\thetag{2} 
Let $ ( \EP _{\rrrr }, \dQbP ) $ and $ (\EP , \dP ) $ be the closures of $ ( \EP _{\rrrr }, \dbP ) $ 
and $ ( \EP , \dbP ) $ on $ \Lm $, respectively. Then, $ (\EP , \dP ) $ is the increasing limit of 
$ ( \EP _{\rrrr }, \dQbP ) $. 
\\\thetag{3} $ (\EP , \dP ) \le ( \Emu , \underline{\dom }^{\mu }) $. 
\end{lemma}
\PF
For $ \upsilon \in \Upsilon $, let 
$ \E ^{\upsilon } (f,g)= \int _{\sSS } \DDD _{\upsilon }[f , g] d\mu $. 
Then, $ ( \E ^{\upsilon } , \dbP )$ is closable on $ \Lm $ by \As{A3}. 
Thus from \eqref{:31y}, \eqref{:32v}, \eqref{:32w}, and \eqref{:33x}, 
$ ( \EP _{\rrrr }, \dbP ) $ and $ ( \EP , \dbP ) $ are countable sums of closable forms $ ( \E ^{\upsilon } , \dbP )$ on $ \Lm $. Hence, we obtain \thetag{1}. 

From \thetag{1} and \eqref{:33y}, 
$ ( \EP , \dbP ) $ is the increasing limit of closable forms $ ( \EP _{\rrrr }, \dbP ) $ on $ \Lm $. 
From this we obtain \thetag{2}. 

From \eqref{:21b}, \eqref{:33a}, and \eqref{:42q}, $ (\EP , \dbP ) \le (\Emu , \dbmu )$. This yields \thetag{3}. 
\PFEND 
\begin{lemma} \label{l:42} 
Assume \As{A2}--\As{A4}. Then 
$ (\EP , \dP ) = ( \Emu , \underline{\dom }^{\mu }) = \ED  $. 
\end{lemma}
\PF 
From \lref{l:41}\thetag{3} and \lref{l:21}, we obtain 
\begin{align}\label{:42f}&
(\EP , \dP ) \le ( \Emu , \underline{\dom }^{\mu }) \le \ED  
.\end{align}
By \As{A4}, we have $ (\EP , \dP ) = \ED  $. 
This and \eqref{:42f} complete the proof. 
\PFEND

Let $ T_t^{\mu }$ be the Markovian semi-group on $ \Lm $ associated with 
$ \ED  $ on $ \Lm $. 
Then there exists an unlabeled diffusion $ (P_{\sss }, \mathsf{X}_t)$ 
associated with the Dirichlet form $ \ED  $ on $ \Lm $ \cite{o.dfa,o.isde,o.rm,o-t.tail}. 
By construction, $ (P_{\sss }, \mathsf{X}_t)$ is $ \mu $-reversible and 
$ T_t^{\mu } f (\sss ) = E_{\sss }[ f (\XX _t )] $  for each $ f \in \Lm $. 
%
From \lref{l:42}, $ (P_{\sss }, \mathsf{X}_t)$ is also associated with $ (\EP , \dP ) $ and $ (\Emu , \dImu ) $. 

\section{One-labeled, tagged particle, and environment processes} \label{s:III}
In \sref{s:III}, we introduce the three stochastic processes related to the tagged particle problem: 
one-labeled, tagged particle, and environment processes. 

We can obtain these three stochastic processes from the labeled process $ \X =(X^i)_{i\in \mathbb{N}}$ by change of coordinate. These stochastic processes are diffusion processes and we identify the associated Dirichlet forms. 
We remark that the original labeled process $ \X $ does not have any associated Dirichlet form.

\subsection{One-labeled processes} \label{s:5} 

Let $ (P_{\sss }, \mathsf{X}_t)$ be the $ \mu $-reversible diffusion associated with $ \ED  $ on $ \Lm $ given in \ssref{s:21}. 
From \As{A1}--\As{A3}, $ (P_{\sss }, \mathsf{X}_t)$ satisfies \eqref{:12e}. 
Then using $ \lpath $ defined by \eqref{:12c}, we have the corresponding labeled process $ \X = \lpath (\XX )$. 
From the labeled process $ \X = ( X ^i)_{i\in\mathbb{N} }$, 
we construct the one-labeled processes $ ( X ^i , \XX ^{\diai })$, $ i \in \mathbb{N}$, where 
$ \XX _t^{\diai } = \sum_{j\ne i }^{\infty} \delta_{ X _t^j } $. 

The same construction is also possible for the unlabeled processes given by 
$ (\EP , \dP ) $ and $ ( \Emu , \underline{\dom }^{\mu }) $ on $ \Lm $. 
If we assume \As{A4} in addition, then these three processes are the same by \lref{l:42}. 
We shall present the Dirichlet form associated with the one-labeled processes $ ( X ^i , \XX ^{\diai })$, $ i \in \mathbb{N}$. The Dirichlet form is independent of $ i \in \mathbb{N}$ because of the symmetry of the particle system.

Let $ \nabla $ be the nabla in $\Rd $. We set for $ f , g \in C_0^{\infty}(\Rd ) $ 
\begin{align} & \notag 
 \nablaR [ f , g ] (x)= \frac{1}{2} ( \nabla f , \nabla g )_{\Rd } (x)
.\end{align}
We naturally regard $ \nablaR $ and $ \DDD $ as the carr\'{e} du champs on $ C_0^{\infty}(\Rd ) \ot \db $ in such a way that, 
for $ f = f_1 \ot f_2 $ and $ g = g_1 \ot g_2 $, 
\begin{align} & \notag 
 \nablaR [ f , g ] = \nablaR [f_1,g_1]f_2g_2 , \quad 
 \DDD [f,g] = f_1g_1 \DDD [f_2 ,g_2] 
.\end{align}
We regard $ \DDD _{ \rR }^{m} $, $ \DDD ^{\perp } $, and $ \DDDR ^{\perp } $ 
as the carr\'{e} du champs on $ C_0^{\infty}(\Rd ) \ot \db $ in the same fashion. 
Let $ \muone $ be the one-Campbell measure of $ \mu $ given by \eqref{:14h}. Let 
\begin{align} &\notag 
 \mathscr{E}_{\rR }^{[1]} (f,g)= \int_{\SR \ts \sSS } \{ \nablaR [ f , g ] + \DDD_{\rR } [f,g] \} d\muone 
\\\label{:50a}&
 \ER ^{[1],\perp } (f,g)= \int_{\SR \ts \sSS } \{ \nablaR [ f , g ] + \DDD_{\rR }^{\perp } [f,g] \} d\muone 
.\end{align}

\begin{lemma} \label{l:50} Assume \As{A1}--\As{A3}. Then 
$ ( \ER ^{[1],\perp } , \db ^{[1],\perp } )$, $ (\Eone _{\rR },\dbone ) $, and $ (\Eone _{\rR },\dzone ) $ 
are closable on $ L^2 (\muone ) $. 
\end{lemma}

\begin{proof}

Let $ \muRss $ be the regular conditional probability measure such that 
\begin{align} & \label{:51f}
\muRss = 
\mu (\,\piR (\xx ) \in \cdot \, \vert \,\piRc (\mathsf{x}) = 
\piRc (\sss ) ) 
.\end{align}
Let $ \muRss ^{[1]} $ be the one-Campbell measure of $ \muRss $. 
We set for $ m \in \mathbb{N} $
\begin{align}
& \label{:51i}
\muRss ^{[1], m } ( \cdot ) = \muRss ^{[1]} (\cdot \cap \SR \ts \SSRmm )
.\end{align}
Let $ \HRm (\mathsf{x}) $ be as in \eqref{:12z}. 
Recall that $ \mu $ is a $ \Psi $-quasi-Gibbs measure. 
Then using \eqref{:12y}, \eqref{:51f}, and \eqref{:51i}, we have 
for $ ( x , \xx ) \in \SR \ts \SSR ^{m-1} $ 
%
\begin{align} \label{:51j}&
\cref{;o1M}^{-1} e^{- \beta \HRm (\delta_x + \mathsf{x}) } 
 \LambdaR ^{[1], m } (d x d\mathsf{x})
 \le 
	\muRss ^{[1], m } ( d x d\mathsf{x}) 
 \le 
 \cref{;o1M} e^{- \beta \HRm (\delta_x + \mathsf{x}) } 
 \LambdaR ^{[1], m } (d x d\mathsf{x})
.\end{align}
Here, 
 $ \LambdaR ^{[1], 1 } (d x d\mathsf{x})= 1_{\SR }(x)dx \delta_{\mathsf{0}}(\mathsf{x})$ and 
$ \LambdaR ^{[1], m } (d x d\mathsf{x}) = 1_{\SR }(x)dx \LambdaR ^{m-1} (d\mathsf{x}) $ 
for $ m \ge 2 $. 
Furthermore, 
$ \Ct \label{;o1M} $ is a positive constant depending on $ \beta $, $ \rR $, $ m $, and $ \piRc (\sss ) $. 
We set 
\begin{align}\label{:51k} &
\ERsonem ( f , g ) = \int \{ \nablaR [ f , g ] + \DDD_{\rR } [f,g] \} d\mu \Rsonem 
.\end{align}
%
Using \As{A3}, \eqref{:51j}, and \eqref{:51k} and applying the method in \cite[Lemma 3.2]{o.dfa}, we see that 
$ (\ERsonem , \dbone ) $ is closable on $ L^2 (\mu \Rsonem )$. 
Let 
\begin{align}\label{:51m}&
 \mathscr{E}_{\rR }^{[1],m} ( f , g ) 
= \int \{ \nablaR [ f , g ] + \DDD_{\rR } [f,g] \} 
d\mu _{\rR }^{[1],m} 
.\end{align}
We write $ \piR (\sss ) = \xx $ and $ \piRc (\sss ) = \yy $. 
Thus $ \sss = (\xx , \yy )$. From \eqref{:51i}, we have 
\begin{align}\label{:51n}&
\muone (d x d\sss )= 
\int \sum_{m=1}^{\infty} \mu \Rsonem (d x d\mathsf{x}) \mu \circ (\piRc ) ^{-1} (d\yy )
.\end{align}
From \eqref{:51m} and \eqref{:51n}, we deduce that 
$ ( \mathscr{E}_{\rR }^{[1],m} ,\dbone ) $ on $ L^2 (\muone ) $ 
is the integral of closable bilinear forms $ (\ERsonem , \dbone ) $ on $ L^2 (\mu \Rsonem )$. 
Hence using the method in \lq\lq Proof of Theorem 4'' in \cite[p.130]{o.dfa}, 
we see that $ ( \mathscr{E}_{\rR }^{[1],m} ,\dbone ) $ is closable on $ L^2 (\muone ) $. 

From \eqref{:51m} and \eqref{:51n}, 
$  \mathscr{E}_{\rR }^{[1]} = \sum_{m=1}^{\infty} \mathscr{E}_{\rR }^{[1],m} $. 
Because $ ( \mathscr{E}_{\rR }^{[1]} ,\dbone ) $ is a countable sum of closable forms 
$ ( \mathscr{E}_{\rR }^{[1],m} ,\dbone ) $ on $ L^2 (\muone ) $, 
 $ ( \mathscr{E}_{\rR }^{[1]} ,\dbone ) $ is closable on $ L^2 (\muone ) $. 

The proof of the closability of $ ( \ER ^{[1],\perp } , \db ^{[1],\perp } )$ on $ \Lm $ is similar to that of $ (\Eone _{\rR } ,\dbone ) $ on $ \Lm $. Hence we omit it. 
\end{proof}

We set 
\begin{align}\label{:51x}&
\Eone (f,g) = \int_{\Rd \ts \sSS } \{ \nablaR [ f , g ] + \DDD [f,g] \} d\muone 
,\\ \label{:51y}&
\dbone = \{ f \in C_0^{\infty}(\Rd ) \ot \db ; 
\Eone (f,f) < \infty , \, f \in L^2 (\muone ) \} 
.\end{align}
We define $ \dzone $ by \eqref{:51y} through replacing $ \db $ with $ \dz $. 
We define $ \E ^{[1],\perp } $ 
by \eqref{:51x} through replacing $ \DDD $ with $ \DDD ^{\perp } $. 
We define $ \db ^{[1],\perp } $ by \eqref{:51y} 
through replacing $ \db $ and $ \Eone $ with $ \dbP $ and $ \E ^{[1],\perp } $, respectively. 
Here $ \dbP $ was given by \eqref{:42q}. 
We deduce from \eqref{:33a} and $ \db \supset \dz $ 
\begin{align}\label{:51z}&
 ( \E ^{[1],\perp } , \db ^{[1],\perp } ) \le (\Eone ,\dbone ) \le (\Eone ,\dzone ) 
.\end{align}

\begin{lemma} \label{l:5!}
Assume \As{A1}--\As{A3}. Then $ ( \E ^{[1],\perp } , \db ^{[1],\perp } )$, $ (\Eone ,\dbone ) $, and $ (\Eone ,\dzone ) $ 
are closable on $ L^2 (\muone ) $. 
\end{lemma}
\begin{proof}
Because $ (\Eone ,\dbone ) $ is the increasing limit of $ ( \mathscr{E}_{\rR }^{[1]} ,\dbone ) $ and 
 $ ( \mathscr{E}_{\rR }^{[1]} ,\dbone ) $ is closable on $ \Lmone $ by \lref{l:50}, 
$ (\Eone ,\dbone ) $ is closable on $ \Lmone $. 

Note that $ (\Eone ,\dbone ) $ is an extension of $ (\Eone _{\rR },\dzone ) $ and that 
$ (\Eone ,\dbone ) $ is closable on $ \Lmone $. 
Hence, $ (\Eone ,\dzone ) $ is closable on $ \Lmone $ by \lref{l:20}. 

We see that $ ( \E ^{[1],\perp } , \db ^{[1],\perp } ) $ is the increasing limit of $ ( \ER ^{[1],\perp } , \db ^{[1],\perp } )$ as $ \rR \to \infty $. 
From \lref{l:50}, $ ( \ER ^{[1],\perp } , \db ^{[1],\perp } )$ is closable on $ \Lmone $. 
Hence, $ ( \E ^{[1],\perp } , \db ^{[1],\perp } ) $ is closable on $ \Lmone $. 
\end{proof}

Let $ (\E ^{[1], \perp } , \dom ^{[1], \perp } ) $, $ (\EoneL , \doneL ) $, and $ (\Eone , \doneO ) $ 
be the closures of $ ( \E ^{[1],\perp } , \db ^{[1],\perp } )$, $ (\Eone ,\dbone ) $, and $ (\Eone ,\dzone ) $ on $ L^2 (\muone )$, respectively. 
\begin{lemma} \label{l:51}
Assume \As{A1}--\As{A4}. Then 
\begin{align}& \label{:51a}
 (\E ^{[1], \perp } , \dom ^{[1], \perp } ) = (\Eone ,\doneL ) = (\Eone ,\doneO ) 
.\end{align}
The diffusion $ ( X ^i , \XX ^{\diai })$ is associated with the closed forms in \eqref{:51a} on $ \Lmone $. 
\end{lemma}
\PF 
We see $ \overline{\dbP }= \dP $ by definition and 
$ \overline{\dom }^{\mu } = \overline{\dz ^{\mu } }$ by \lref{l:21}. 
From \lref{l:42}, 
$ (\EP , \dP ) = ( \Emu , \underline{\dom }^{\mu }) = \ED  $. 
Combining these, we have $ \overline{\dbP } = \overline{\dz ^{\mu } } $. 
Hence, 
\begin{align}\label{:51q}&
C_0^{\infty}(\Rd ) \ot \dbP \subset C_0^{\infty}(\Rd ) \ot \overline{\dbP } = 
C_0^{\infty}(\Rd ) \ot \overline{\dz ^{\mu } } 
.\end{align}
It is easy to see that 
\begin{align}\label{:51r}&
C_0^{\infty}(\Rd ) \ot \overline{\dz ^{\mu } } \subset \overline{C_0^{\infty}(\Rd ) \ot \dz ^{\mu } } \equiv \doneO 
.\end{align}
From \eqref{:51q} and \eqref{:51r}, we obtain 
\begin{align}\label{:51R}&
 (\E ^{[1],\perp } ,\dom ^{[1], \perp } ) \ge (\Eone ,\doneO ) 
.\end{align}
From \eqref{:51z}, we see 
\begin{align}\label{:51s}&
 ( \E ^{[1],\perp } , \dom ^{[1], \perp } ) \le (\Eone ,\doneL ) \le (\Eone ,\doneO ) 
.\end{align}
Hence combining \eqref{:51R} and \eqref{:51s}, we obtain \eqref{:51a}. 

Using \cite[Theorem 2.4]{o.tp}, we see that 
$ ( X ^i , \XX ^{\diai })$ is the diffusion associated with the Dirichlet form 
$ (\Eone , \doneO ) $ on $ \Lmone $. Hence, \eqref{:51a} yields the second claim. 
\PFEND 

\subsection{Tagged particle processes} \label{s:6}
The tagged particle problem of interacting Brownian motions is to prove the diffusive scaling limit of each particle in the system \cite{gp,De,o.inv2,o.p}. The standard device for this problem is to introduce the environment process seen from the tagged particle \cite{gp}.
We use this device and define it by changing the coordinates as follows. 
For the labeled process $ \X = (X^i)_{i\in\mathbb{N}}$ given after \eqref{:12e}, we set 
\begin{align}\label{:61u}& 
 X = X ^{1} ,\quad Y^{i} = X ^{i+1} - X ^{1} \quad (i\in\mathbb{N} )
.\end{align}
Here, $ X $ denotes the tagged particle and $ \mathbf{Y}=(Y^{i})_{i\in\mathbb{N}}$ 
is the labeled environment process seen from the tagged particle. 
The unlabeled environment process $ \YY = \{ \YY _t \}_t $ is given by 
\begin{align}\label{:61v}&
\YY _t = \sum_{i=1}^{\infty} \delta_{Y_t^i }
.\end{align}

For $ R,T \in \mathbb{N} $ and $ \mathbf{w}=(w^i)_{ i \in \mathbb{N} } \in \CiRdN $, we set 
\begin{align}& \label{:61w}
\IRT (\mathbf{w}) = 
\sup\{ i \in \mathbb{N}; 
\min_{t\in[0,T]} \lvert w^i(t)\rvert \le \rR \} 
.\end{align}
For $ \mathbf{w}=(w^i)_{ i \in \mathbb{N} } \in \CiRdN $, 
let $ \mathsf{w} $ be such that $ \mathsf{w}_t = \sum_{i=1}^{\infty} \delta_{w^i(t)}$. 
If 
$ \IRT (\mathbf{w}) < \infty $ for all $ R,T \in \mathbb{N}$, then 
$ \mathsf{w} $ is an $ \sSS $-valued continuous path. 
 See \lref{l:Z3} for proof. Conversely, if 
$ \IRT (\mathbf{w}) = \infty $ for some $ R$ or $T \in \mathbb{N}$, then 
$ \mathsf{w} $ is {\em not} necessarily an $ \sSS $-valued continuous path. 
See \cite[Remark 3.10]{o-t.tail}. 

\begin{example}\label{r:51}
We present the ISDE of $ ( X , \mathbf{Y} )$ for the Ginibre interacting Brownian motion. 
Using \eqref{:10a} and \eqref{:61u}, we see that $ X $ and $ \mathbf{Y}=(Y^i)_{i\in\mathbb{N}}$ satisfy 
\begin{align} \notag 
dX_{t} &= dB_t^1 - \lim_{ \rR \to\infty }
( \sum_{\lvert Y^{j}_{t}\rvert < \rR ,\, j \in \mathbb{N} } \frac{Y^{j}_{t}}{\lvert Y^{j}_{t}\rvert ^{2}} ) dt 
,\\ \notag 
dY^{i}_{t} &= \sqrt{2} \, d\tilde{B}_t^i + \frac{Y^{i}_{t}}{\lvert Y^{i}_{t}\rvert ^{2}} dt 
 + 
 \lim_{ \rR \to\infty } (\sum_{\lvert Y^{j}_{t}\rvert < \rR ,\, \atop j \in \mathbb{N} }
 \frac{Y^{j}_{t}}{\lvert Y^{j}_{t}\rvert ^{2}} ) dt 
 + 
\lim_{ \rR \to\infty } (\sum_{\lvert Y^{i}_{t} - Y^{j}_{t}\rvert < \rR ,\, \atop i \not= j ,\, j \in \mathbb{N} }
\frac{Y^{i}_{t} - Y^{j}_{t}}{\lvert Y^{i}_{t} - Y^{j}_{t}\rvert ^{2}} ) dt 
.\end{align}
Here, $ \tilde{B}_{t}^i = ({1}/{\sqrt{2}}) (B^{i+1}_{t} - B^{1}_{t}) $. 
 $ \{ \tilde{B}^i \}_{i\in\mathbb{N}} $ are not independent and each $ \tilde{B}^i $ is 
equivalent in law to the standard Brownian motion in $ \Rtwo $. 
The second equation above is self-contained as an equation of $ \mathbf{Y}$. 
We see that the unlabeled process $ \YY $ of $ \mathbf{Y} $ is a diffusion process with invariant probability measure $ \mugz $ in \sref{s:7}. 
This property is critical in the Kipnis--Varadhan theory to the tagged particle problem \cite{De,KV,o.inv2}. 
\end{example}

Although $ ( X , \mathbf{Y} )$ is a diffusion with state space 
$\Rd \ts (\Rd ) ^{\mathbb{N} } $, there exists no associated Dirichlet space. 
Indeed, suitable invariant measures are lacking for $ ( X , \mathbf{Y} )$. 
In contrast, $ ( X , \YY )$ is a diffusion with an invariant measure. 
As a result, it has the associated Dirichlet space. This fact is important in the analysis of the Dirichlet form version of the Kipnis--Varadhan theory in \cite{o.inv2}. We shall specify the Dirichlet form associated with $ ( X , \YY )$.

Let $ \Dsft $ be the generator of the translation operator given by \eqref{:14m}. Let 
\begin{align}\label{:61x}&
 (\nabla - \Dsft ) [f,g] = \half ( (\nabla - \Dsft )f , (\nabla - \Dsft ) g )_{\Rd } 
.\end{align}
We introduce the bilinear form $ (\EXY , \dbXY )$ such that 
\begin{align} \label{:61y}&
\EXY ( f , g ) = \int_{ \Rd \ts \sSS } (\nabla - \Dsft ) [f,g] + \DDD [ f , g ]
 dx \ts \muz 
,\\ \notag &
\dbXY = \{ f \in C_0^{\infty} (\Rd ) \otimes ( \dsft \cap \db ); 
\EXY (f,f) < \infty , f \in L^2 ( dx \ts \muz ) \} 
.\end{align}
Here, $ \dsft $ is the domain of $ \Dsft $: i.e., $ \dsft $ 
is the set of functions for which the limit in \eqref{:14m} exists for all $
 \sss \in \sSS $ and $ \5 $. 
We define $ \E ^{\XY ,\perp } $ through \eqref{:61y} by replacing $ \DDD $ with $ \DDD ^{\perp }$. 
We define $ \dz ^{\XY } $ by replacing $ \db $ with $ \dz $ in $ \dbXY $. 
We also set $ \db ^{\XY , \perp } $ by replacing 
$ (\db , \EXY ) $ with $ (\dbP , \E ^{\XY ,\perp }) $ in $ \dbXY $. 
\begin{lemma} \label{l:61}
Make the same assumptions as \lref{l:51}. Then the following hold. 
\\
\thetag{1} 
$ ( \E ^{\XY ,\perp } , \db ^{\XY , \perp } )$, 
$ (\EXY , \dbXY )$, and $ (\EXY , \dz ^{\XY } )$ are closable on $ L ^2 (dx \ts \muz ) $. 
\\\thetag{2} 
Let $ (\E ^{\XY ,\perp }, \dom ^{\XY ,\perp } ) $, $ (\EXY , \dXYL )$ and $ (\EXY , \dXYU )$ 
 be the closures of the closable forms in \thetag{1} on $ L ^2 (dx \ts \muz ) $, respectively. Then, 
\begin{align} \label{:61a}& 
 (\E ^{\XY ,\perp }, \dom ^{\XY ,\perp }) = (\EXY , \dXYL ) = (\EXY , \dXYU ) 
.\end{align}
\thetag{3} The diffusion $ ( X , \YY )$ is associated with 
$ (\E ^{\XY ,\perp }, \dom ^{\XY ,\perp }) $ on $ L^2 ( dx \ts \muz ) $, where 
 $ X $ and $ \YY $ are as in \eqref{:61u} and \eqref{:61v}, respectively. 
\end{lemma}
\PF 
Let $ \{ \vartheta_x \}_{x\in\Rd }$ be as in \eqref{:12t}. 
Let $ \iota $ be the transformation on $ \Rd \ts \sSS $ defined by 
$ \iota ( x , \sss ) = ( x , \vartheta_x (\sss ))$. 
%
Using this and \As{A1}, we deduce that 
\begin{align}\label{:61f}&
\muone \circ \iota ^{-1} = dx \ts \muz 
.\end{align}
Hence, $ \iota $ is the unitary transformation between $ \Lmone $ and 
$ L ^2 (dx \ts \muz ) $ such that 
\begin{align}\label{:61h}
(f \circ \iota , g \circ \iota )_{\Lmone } =& (f , g )_{L ^2 (dx \ts \muz ) } 
.\end{align}

Recall that 
$ \vartheta_x (\sss ) = \sum_i \delta_{s_i - x} $ by \eqref{:12t}. Hence 
\begin{align}\label{:61l}&
\PD{}{x} \Big( f ( x , \vartheta_x (\sss )) \Big) = 
\Big(\PD{}{x} f \Big) ( x , \vartheta_x (\sss )) - \Big(\Dsft f \Big) ( x , \vartheta_x (\sss )) 
.\end{align}
Hence from \eqref{:61x}, \eqref{:61l}, and $\iota ( x , \sss ) = ( x , \vartheta_x (\sss )) $, we deduce that 
\begin{align}\notag &
\big( \nabla [ f \circ \iota , g \circ \iota ] + 
 \DDD ^{\perp } [f \circ \iota , g \circ \iota ] \big) \circ \iota ^{-1}= 
 (\nabla - \Dsft ) [f,g] + \DDD ^{\perp } [f,g] 
,\\ &\label{:61m}
\big( \nabla [ f \circ \iota , g \circ \iota ] + 
 \DDD [f \circ \iota , g \circ \iota ] \big) \circ \iota ^{-1}= 
 (\nabla - \Dsft ) [f,g] + \DDD [f,g] 
.\end{align}
From \eqref{:61y}, \eqref{:61f}, and \eqref{:61m}, we see the isometry 
of the bilinear forms such that 
\begin{align}\label{:61o}&
\E ^{[1], \perp } ( f \circ \iota , g \circ \iota ) = \E ^{\XY ,\perp } ( f , g ) 
, \quad 
\Eone ( f \circ \iota , g \circ \iota ) = \EXY ( f , g )
.\end{align}
Using \eqref{:61h} and \eqref{:61o} together with \lref{l:5!}, we obtain \thetag{1}. 

Combining \eqref{:51a} and \eqref{:61o}, we obtain \thetag{2} immediately. 

We regard $ \iota $ as the transformation of $ C([0,\infty ); \Rd \ts \sSS )$, 
denoted by the same symbol $ \iota $, such that 
$ \iota ( X^1 , \XX ^{\diaone } ) = \{ \iota ( X _t^1 , \XX _t^{\diaone } ) \}_{t\in[0,\infty)}$. Then 
\begin{align}\label{:61p}&
 \iota ( X^1 , \XX ^{\diaone } ) = ( X , \YY ) 
. \end{align}
Using \lref{l:51}, \eqref{:61h}, \eqref{:61o}, and \eqref{:61p}, we obtain \thetag{3}. 
\PFEND 

\subsection{Environment processes} \label{s:7}
Let $ \YY $ be the environment process given by \eqref{:61v}. 
Note that $ \YY $ itself is a diffusion. 
Hence, we specify the Dirichlet form associated with $ \YY $. 

Let $ \Dsft $, $ \DDD $, and $ \Dp $ be as in \eqref{:14M}, \eqref{:21r} and \eqref{:33x}, respectively. 

Let $ ( \E ^{\YY ,\perp } , \ \db ^{\YY , \perp } ) $ be the bilinear form defined by 
\begin{align} & \notag 
 \E ^{\YY ,\perp } (f,g) = \int_{\sSS } 
 \Dsft [f,g] + \DDD ^{\perp } [f,g]
d\muz 
,\\ \label{:71z}&
 \db ^{\YY , \perp } = \{ f \in \dsft \cap \db ; \E ^{\YY , \perp } ( f , f ) < \infty ,\, f \in \LMz \} 
.\end{align}
Let $ (\EY , \dbY )$ and $ (\EY , \dzY )$ be the bilinear forms defined by 
\begin{align} & \notag 
\EY (f,g) = \int_{\sSS }\Dsft [f,g] + \DDD [f,g] d\muz 
,\\ &\notag 
 \dbY = \{ f \in \dsft \cap \db ; \, \EY (f,f) < \infty ,\, f \in \LMz \} 
,\\ & \label{:71y}
 \dzY = \{ f \in \dsft \cap \dz ; \, \EY (f,f) < \infty ,\, f \in \LMz \} 
.\end{align}
For a random variable $ Z $ and a probability measure $ \nu $, we write $ Z \elaw \nu $ 
if the law of $ Z $ coincides with $ \nu $. 
\begin{lemma} \label{l:71} 
Make the same assumptions as \lref{l:51}. Then the following hold. 
\\
\thetag{1} 
$ (\E ^{\YY , \perp } ,\db ^{\YY ,\perp } ) $, $ (\EY ,\dbY )$, and $ (\EY ,\dzY )$ are closable on $ \LMz $. 
\\\thetag{2} 
Let $ (\E ^{\YY , \perp } , \dom ^{\YY , \perp } ) $, $ (\EY , \dYL ) $, and $ (\EY , \dYU ) $ be the closures of the closable forms in \thetag{1} on $ \LMz $, respectively. 
Then, 
\begin{align}\label{:71a}&
 (\E ^{\YY , \perp } , \dom ^{\YY , \perp } ) = (\EY , \dYL ) = (\EY , \dYU ) 
.\end{align}
\thetag{3} 
Let $ ( X , \YY )$ be as in \eqref{:61u} and \eqref{:61v}. 
 Suppose that $ ( X _0 , \YY _0) \elaw \zeta \times f d\muz $, 
where $ \zeta $ is a probability measure on $ \Rd $ and 
$ 0 \le f \in \db $ such that $ \int_{\sSS } f d\muz =1$. 
 The distribution of $ \YY $ in $ ( X , \YY )$ is then the same as that of the diffusion $ \YY $ 
 associated with $ (\E ^{\YY , \perp } , \dom ^{\YY , \perp } ) $ such that $ \YY _0 \elaw f d\muz $. 
\end{lemma}
\PF 
For $ \varphi \in C_0^{\infty}(\mathbb{R}^d)$ and $ f \in \db ^{\YY ,\perp } $, we see 
\begin{align}\label{:71f}
\| \varphi \ot f \|_{L ^2 (dx \ts \muz ) } &= \| \varphi \|_{L ^2 (dx ) }\| f \|_{L ^2 (\muz ) }, 
\\ \label{:71g}
\E ^{\XY ,\perp } (\varphi \ot f , \varphi \ot f ) &= 
\| \varphi \|_{L ^2 (dx ) }^2 \E ^{\YY , \perp } ( f , f ) + 
\| \nabla \varphi \|_{L ^2 (dx ) }^2 \| f \|_{L ^2 (\muz ) }^2
.\end{align}
Indeed, \eqref{:71f} is a straightforward calculation and \eqref{:71g} follows from the following. 
\begin{align} \notag &
 (\nabla - \Dsft ) [ \varphi \ot f , \varphi \ot f ] + \Dp [ \varphi \ot f , \varphi \ot f ] 
\\ \label{:71h}& 
= 
\varphi^2 \ot 
 \{\Dsft [f,f] + \DDD ^{\perp } [f,f] \} 
+ \lvert \nabla \varphi \rvert^2 \ot f^2 
- 2 (\varphi \nabla \varphi , f \Dsft f )_{\mathbb{R}^2}
.\end{align}
Integrating \eqref{:71h} over $ \mathbb{R}^d \ts \sSS $ with respect to $ dx \ts \muz $ and using 
\begin{align} &\notag 
\int _{ \mathbb{R}^d \ts \sSS } (\varphi \nabla \varphi , f \Dsft f )_{\mathbb{R}^d }
 dx \ts \muz 
= (
\int _{ \mathbb{R}^d } \varphi \nabla \varphi dx , 
\int _{\sSS } f \Dsft f d\muz 
)_{\mathbb{R}^d} 
= 0 
,\end{align}
we obtain \eqref{:71g}. 

Let $ \| \varphi \|_{L^2(\muz )} \ne 0 $. 
Let $ \{ f_n \} $ be an $ \E ^{\YY ,\perp }$-Cauchy sequence in $ \db ^{\YY ,\perp }$ such that 
$ \limi{n} \| f_n \|_{L^2 (\muz )} = 0 $. 
Then from \eqref{:71f} and \eqref{:71g}, we see that 
$ \{\varphi \ot f_n \}$ is an $ \E ^{\XY ,\perp }$-Cauchy sequence such that 
\begin{align*}&
 \limi{n} \| \varphi \ot f_n \|_{L ^2 (dx \ts \muz ) } = 0
.\end{align*}
By \lref{l:61}\thetag{1}, $ ( \E ^{\XY ,\perp } , \db ^{\YY ,\perp } )$ is closable on 
$ L ^2 (dx \ts \muz ) $. 
Hence we deduce 
\begin{align}& \label{:71k}
\limi{n} \E ^{\XY ,\perp }( \varphi \ot f_n , \varphi \ot f_n) = 0 
.\end{align}
From \eqref{:71g}, \eqref{:71k}, $ \| \varphi \|_{L^2(dx )} \ne 0 $, 
and $ \limi{n} \| f_n \|_{L^2 (\muz )} = 0 $, we deduce 
\begin{align}& \notag 
\E ^{\YY ,\perp } (f_n , f_n ) = 
\frac{1}{\| \varphi \|_{L ^2 (dx ) }^2}
\Big\{ 
\E ^{\XY ,\perp } (\varphi \ot f_n , \varphi \ot f_n )
-
\| \nabla \varphi \|_{L ^2 (dx ) }^2 \| f_n \|_{L ^2 (\muz ) }^2
\Big\}
\to
0 
\end{align}
as $ n \to \infty $. This implies the first claim in \thetag{1}. The proof of the second and third claim in \thetag{1} 
is same. Hence, we omit it. 
Claim \thetag{2} follows from \lref{l:61}\thetag{2} and \eqref{:71g}. 
Claim \thetag{3} follows from \lref{l:61}\thetag{3}, \eqref{:71f}, and \eqref{:71g}. 
\PFEND

\section{Vanishing the self-diffusion matrix}\label{s:IV}
In \sref{s:IV}, we present a sufficient condition of vanishing the self-diffusion matrix. 
The goal of \sref{s:IV}, is to prove \tref{l:93}, which is an \lq\lq in $ \mu_0$-measure'' version of \tref{l:15}\thetag{2}.

\subsection{A sufficient condition for the vanishing of the self-diffusion matrix} \label{s:8}

In \ssref{s:8}, $ \mu $ is a random point field on $ \Rd $ satisfying \As{A1}--\As{A4}. 

Let $ (\EY ,\dbY )$ be the bilinear form defined by \eqref{:71y}. We set 
\begin{align}\label{:81w}
\EYone (f,g) = &\int_{\sSS } \Dsft [f,g] d\muz 
,\quad 
\EYtwo (f,g) = \int_{\sSS } \DDD [f,g] d\muz 
.\end{align}
From this, we have a decomposition of the bilinear form such that 
\begin{align}\label{:81x}&
\EY (f,g) = \EYone (f,g) + \EYtwo (f,g) \quad \text{ for $ f , g \in \dbY $}
.\end{align}
Using \eqref{:81w} and the obvious inequalities $ \EYi (f,f) \le \EY (f,f) $, 
we extend the domain of $ \EYi $ from $ \dbY $ to $ \dYL $, where $ i=1,2$. 
Hence, \eqref{:81x} yields 
\begin{align}\label{:81y}& 
\EY (f,g) = \EYone (f,g) + \EYtwo (f,g) \quad \text{ for $ f , g \in \dYL $}
.\end{align}
Because of \lref{l:71}, $ ( \EY , \dbY )$ is closable on $ \LMz $. 
Meanwhile, each of $ (\EYone ,\dbY )$ and $ (\EYtwo ,\dbY )$ is not necessarily closable 
on $ \LMz $. Still, \eqref{:81y} makes sense for $ f , g \in \dYL $. 

$ \dYL $ is a Hilbert space with inner product $ \EY + (\cdot , * )_{\LMz }$ and 
$ \dYL $ is a pre-Hilbert space with non-negative bilinear form $ \EY $. 
Let $ \sim $ be the equivalence relation on $ \dYL $ such that $ f \sim g $ if and only if 
$ \EY (f-g,f-g) = 0 $. 
The quotient space $ \dYL / \sim $ is a pre-Hilbert space with inner product $ \EYLt $ such that 
\begin{align}\label{:81a}&
\EYLt (\tilde{f}, \tilde{g}) = \EY (f,g) 
 \quad \text{ for } f,g \in \dYL 
,\end{align}
where $ \tilde{f}= f/\sim $ and $ \tilde{g} = g /\sim $. 
The completion $ \dYLt $ of $ \dYL / \sim $ is then a Hilbert space with inner product $ \EYLt $. 

Let $ \Dsft _{ p } $ be as defined in \eqref{:14m}. 
For $ \5 $ and $ g \in \dbY $, we set 
\begin{align}\label{:81c}&
F_{ p } ( g ) = \int_{\sSS } \half \Dsft _{ p } g \, d\muz 
.\end{align}
By the Schwartz inequality, \eqref{:81x}, and \eqref{:81a}, we obtain for any $ g \in \dbY $
\begin{align}\label{:81d}&
\lvert F_{ p } ( g ) \rvert ^2 \le \half \EYone ( g , g ) 
\le \half \EY ( g , g ) 
= \half \EYLt (\tilde{g},\tilde{g}) 
.\end{align}
From \eqref{:81d}, we regard $ F_{ p }$ as a bounded linear functional on 
$ \dYL $ and $ \dYLt $, and we denote it by the same symbol $ F_{ p }$. 
\begin{lemma} \label{l:81} 
For $ \5 $, there exists a unique solution $ \psiP \in \dYLt $ of the equation 
\begin{align} &\notag 
\EYLt (\psiP , g ) = F_{ p } ( g )
\quad \text{ for all } g \in \dYLt 
.\end{align}
\end{lemma}

\PF 
Because we regard $ F_{ p } $ as a bounded linear functional of the Hilbert space 
$ \dYLt $ with inner products $ \EYLt $, \lref{l:81} is obvious from the Riesz theorem. 
\PFEND 

We consider a resolvent equation. For each $ \la > 0 $ and $ \5 $, 
let $ \psilaP \in \dYL $ be the unique solution of the equation such that for any $ g \in \dYL $ 
\begin{align}\label{:82z}&
\EY (\psilaP , g ) + \la (\psilaP , g )_{\LMz } = F_{ p } ( g ) 
.\end{align}

\begin{lemma} \label{l:82} 
For each $ \5 $, the following hold. \\
\thetag{1} 
$ \{ \psilaP \}_{\la > 0 }$ is an $ \EY $-Cauchy sequence in $ \dYL $ satisfying 
\begin{align}\label{:82a}&
\limz{\la , \la ' } \EY ( \psilaP - \psilaaP , \psilaP - \psilaaP ) = 0 
,\\\label{:82b}&
\limz{\la } \la \| \psilaP \|_{\LMz }^2 = 0 
.\end{align}
\thetag{2} 
$ \{ \psilaPt \}_{\la > 0 }$ is an $ \EYLt $-Cauchy sequence in $ \dYLt $ satisfying 
\begin{align}\label{:82c}&
\limz{\la } \EYLt ( \psilaPt - \psiP , \psilaPt - \psiP ) = 0 
.\end{align}
Here, $ \psiP \in \dYLt $ is the limit of $ \{ \psilaPt \}_{\la > 0 }$ in $ \EYLt $. 
\end{lemma}
\PF 
\lref{l:82} follows from the standard argument; see \cite{KV,o.inv1}. 
\PFEND

\begin{lemma} \label{l:83} 
Assume \As{A7}. Then, for $ \5 $, 
\begin{align}\label{:83c}&
\limz{\la } \Big\{ 
\int_{\sSS } 
\0 
\muz (d\yy ) 
+ \EYtwo ( \psilaP ,\psilaP ) \Big\}
= 0 
.\end{align}
In particular, $ \alpha _{p,p} = 0 $ for $ \5 $. 
\end{lemma}

\PF Let $ g \in \dbY $. From \eqref{:81w} and \eqref{:81y}, we have 
\begin{align} \notag &
 \EY ( \ChiLp , g ) - \int_{\sSS } \half \Dsft _{ p } g \, d\muz 
\\ \notag & = 
 \int_{\sSS } 
\half \sumQ 
( \DsftQ \ChiLp )( \Dsft _{ q } g ) 
\, d\muz + 
\EYtwo ( \ChiLp , g ) 
- \int_{\sSS } \half \Dsft _{ p } g \, d\muz 
\\ \notag &
= 
 \int_{\sSS } 
\half \sumQ 
( \DsftQ \ChiLp - \deltaPQ )( \Dsft _{ q } g ) 
\, d\muz + 
\EYtwo ( \ChiLp , g ) 
\to 0 \quad \text{as $ L \to \infty $}
.\end{align}
Here we used \eqref{:14n} and the Cauchy--Schwartz inequality to the last line. 
Hence 
\begin{align} \label{:83h}&
\limi{L} \EY ( \ChiLp , g ) = \int_{\sSS } \half \Dsft _{ p } g \, d\muz 
\quad \text{ for all $ g \in \dbY $}
.\end{align}
Because $ \{ \ChiLp \} $ is an $ \EY $-Cauchy sequence, 
 $ \{ \ChiLpt \} $ is an $ \EYLt $-Cauchy sequence in $ \dYLt $ by \eqref{:81a}. 
By \eqref{:81a} and \eqref{:83h}, 
$ \{ \ChiLpt \} $ is a weak convergent sequence in the Hilbert space $ \dYLt $. 
By \eqref{:81c} and \lref{l:81}, $ \psiP $ is the limit of $ \{ \ChiLpt \} $. Hence, 
\begin{align}\label{:83i}&
\limi{L} \EYLt ( \psiP - \ChiLpt , \psiP - \ChiLpt ) = 0 
.\end{align}

We write $ \EYLtone ( f ) = \EYLtone (f,f)$. Then, 
\begin{align} & \notag
\int_{\sSS } 
\0 
\muz (d\yy ) 
\\ \notag &
\le 2\int_{\sSS } \Big\{
\half \sumQ 
\Big\lvert \DsftQ \psilaP - \DsftQ \ChiLp \Big\rvert ^2 
+ 
\half \sumQ 
\Big\lvert \DsftQ \ChiLp - \deltaPQ \Big\rvert ^2 
\Big\} 
\muz (d\yy ) 
\\ \notag &
= 
2 \EYLtone ( \psilaPt - \psiP +\psiP - \ChiLpt ) 
+ 2\Big\{
\int_{\sSS }
\half \sumQ 
\Big\lvert \DsftQ \ChiLp - \deltaPQ \Big\rvert ^2 
\muz (d\yy ) 
\Big\}\\ \notag &
\le 4
\EYLtone ( \psilaPt - \psiP ) 
 + 4
\EYLtone ( \psiP - \ChiLpt ) 
+ 
2 \Big\{\int_{\sSS } 
\half \sumQ 
\Big\lvert \DsftQ \ChiLp - \deltaPQ \Big\rvert ^2 
\muz (d\yy ) \Big\}
.\end{align}
Taking $ L \to \infty $ and then $ \la \to 0 $ in the last line, 
each term vanishes by \eqref{:82c}, \eqref{:83i}, and \eqref{:14n}. 
Hence, we find that the first term in \eqref{:83c} converges to zero. 

We calculate the second term in \eqref{:83c}. 
We write $ \EYLttwo ( f ) = \EYLttwo (f,f)$. Then, 
\begin{align}\notag 
\EYtwo ( \psilaP ) 
= & 
\EYLttwo (
(\psilaPt - \psiP ) + (\psiP - \ChiLpt ) + \ChiLpt ) 
\\ \notag 
\le &
3 \{
\EYLttwo ( \psilaPt - \psiP )
+
\EYLttwo ( \psiP - \ChiLpt ) 
+
 \EYLttwo ( \ChiLpt ) 
\}
.\end{align}
From \eqref{:82c}, \eqref{:83i},  \eqref{:14n}, and $  \EYLttwo ( \ChiLpt ) =  \EYtwo (\ChiLp )$, 
 the last line converges to zero as $ L \to \infty $ and then $ \la \to 0 $. 
Hence,  the second term in \eqref{:83c} converges to zero. Collecting these, we have \eqref{:83c}. 
From $ \psilaP \in \dbY $, we obtain the second claim. 
\PFEND

\subsection{Invariance principle with the vanishing of the self-diffusion matrix} \label{s:9}

Let $ ( X , \YY )$ be as in \lref{l:61}. 
From \lref{l:61}, we deduce that $ ( X , \YY ) $ is the diffusion associated with 
$ (\EXY , \dXYL ) $ on $ L ^2 (dx \ts \muz ) $. 
We take $ ( X _0, \YY _0) \elaw \2 $, where $ \zeta $ is a probability measure on $ \Rd $. 
The distribution of the second component $ \YY $ of $ ( X , \YY )$ is independent of $ \zeta $ 
from \lref{l:71}. 

Let $ \YY $ be as in \lref{l:71}. 
Then $ \YY $ is the diffusion associated with $ (\EY , \dYL )$ on $ \LMz $. 
If $ \YY _0 \elaw \muz $, then $ \YY $ is a $ \muz $-stationary Markov process and is equivalent in law to 
the second component $ \YY $ in $ ( X , \YY )$. See \lref{l:71}.

The Dirichlet forms $ (\EXY , \dXYL ) $ on $ L ^2 (dx \ts \muz ) $ and 
 $ (\EY ,\dYL )$ on $ \LMz $ are quasi-regular and strongly local (see \ssref{s:72} and \cite{c-f}) . 
 Indeed, the quasi regularity of the upper Dirichlet forms 
 $ (\EXY , \dXYU )$ on $ L ^2 (dx \ts \muz ) $ and $ (\EY ,\dYU )$ on $ \LMz $ were proved in \cite{o.tp}. 
Furthermore, the strong locality of these Dirichlet forms are clear from the structure of the carr\'{e} du champs giving the Dirichlet forms. 
 Hence, using the identities of the upper and lower Dirichlet forms in 
 \lref{l:61} and \lref{l:71}, we obtain the quasi regularity and the strong locality of the lower Dirichlet forms 
 $ (\EXY , \dXYL ) $ on $ L ^2 (dx \ts \muz ) $ and $ (\EY ,\dYL )$ on $ \LMz $. 
Once we have verified the quasi regularity and the strong locality, we can apply the Dirichlet form theory developed in \cite{c-f,fot.2} to the associated diffusion processes $ ( X , \YY )$ and $ \YY $. 


Let $ \psila = (\psilaP )_{ p = 1}^{\dd } $ be the function given by \eqref{:82z}. 
Recall that $ 1 \otimes \psila \in \dXYL $ locally and that $ \psila \in \dYL $. 
Let $ \widetilde{1 \otimes \psila} $ and $ \widetilde{\psi}_{\la } $ be quasi continuous modifications of $ 1 \otimes \psila $ and $ \psila $, respectively. 
We can and do take $ \widetilde{1 \otimes \psila} = 1 \ot \widetilde{\psi}_{\la } $. 

Because $ (\EXY , \dXYL ) $ and $ (\EY ,\dYL )$ are strongly local, 
$ 1 \ot \widetilde{\psi}_{\la } ( X _t , \YY _t ) $ and $ \widetilde{\psi}_{\la } (\YY _t )$ are continuous processes (see \cite[Theorem 4.3.4]{c-f}). 
Let $ \Mla = ( \Mlap )_{ p =1}^{\dd } $, $ \la > 0 $, be a continuous process such that 
\begin{align}\label{:91y}&
 \Mla ( t ) = X _t - X _0 + \widetilde{\psi}_{\la } ( \YY _t ) - \widetilde{\psi}_{\la } ( \YY _0 ) 
- \int_0^t \la \psila (\YY _u ) du 
.\end{align}

Recall that $ ( X _0, \YY _0) \elaw \2 $ and $ \YY _0 \elaw \muz $. 
It is clear that 
the distribution of $ (1 \ot \widetilde{\psi}_{\la } ) ( X _t , \YY _t ) =
 \widetilde{\psi}_{\la } (\YY _t ) $ is independent of $ \zeta $. 
\begin{lemma} \label{l:91} 
Assume \As{A1}--\As{A4} and \As{A7}. 
Let $ ( X , \YY )$ be as in \lref{l:61}. 
Let $ ( X _0, \YY _0) \elaw \2 $. 
Then, $ \{\Mla \}_{ \la >0} $ is a sequence of continuous $ L^2 $-martingales with stationary increments such that 
\begin{align}\label{:91b}&
\limz{\la } E [\lvert \Mla (t)\rvert ^2] = 0 \quad \text{ for each $ t $}
.\end{align}
\end{lemma}
\PF 
Note that $ u_p := x_p \otimes 1 + 1 \otimes \psilaP $ is locally in $ \dXYL $. 
Let $ \widetilde{u} $ be a quasi continuous modification of $ u := (u_p)_{p=1}^d $. 
Applying Fukushima's decomposition to $ u $, we see that the additive functional 
$ A^{[u]} = \widetilde{u} ( X _{\cdot} , \YY _{\cdot} ) - \widetilde{u} ( X _0 , \YY _0 ) $ satisfies the decomposition 
\begin{align}\label{:91d}&
 A^{[u]} = M^{[u]} + N^{[u]}
.\end{align}
Here $ M^{[u]} $ is the martingale part and $ N^{[u]} $ is the continuous additive functional of zero energy. 
We refer to \cite[Theorem 4.2.6]{c-f} for Fukushima's decomposition. 

Because $ (\EXY , \dXYL ) $ is strongly local, $ M^{[u]} $ is a continuous local martingale. 
This follows from the Beurling-Deny formula (see \cite[(4.3.6)]{c-f} and \cite[Theorem 5.5.1]{fot.2}). 

We next calculate $ N^{[u]}$. 
Let $ 0 \le \varphi \in \dbXY $ and $ \int \varphi dx \muz (d\sss ) = 1$. 
Then 
\begin{align} \notag 
 \EXY ( u_p , \varphi ) &=
\int_{\Rd \ts \sSS } (\nabla - \Dsft ) [ u_p , \varphi ] + \DDD [ u_p , \varphi ] \, dx \muz (d\sss ) 
\\\notag & = \int_{\Rd \ts \sSS }
\half \sum_{ q = 1}^{d } \delta_{ p q } \PD{\varphi }{x_{ q }} - 
\half \sum_{ q = 1}^{d } \delta_{ p q } \DsftQ \varphi - 
\half \sum_{ q = 1}^{d } \DsftQ \psilaP \PD{\varphi }{x_{ q }}
\\\notag &\quad \quad \quad 
 + 
\half \sum_{ q = 1}^{d } \DsftQ \psilaP \DsftQ \varphi 
+ 
\DDD [ \psilaP , \varphi ] \, dx \muz (d\sss ) 
\\\notag & = \int_{\Rd \ts \sSS }
\half \PD{\varphi }{x_{ p }} - \half \DsftP \varphi 
- 
\half \sum_{ q = 1}^{d } \DsftQ \psilaP \PD{\varphi }{x_{ q }} 
\\\notag &\quad \quad \quad 
+ 
\half \sum_{ q = 1}^{d } \DsftQ \psilaP \DsftQ \varphi 
+ 
\DDD [ \psilaP , \varphi ] \, dx \muz (d\sss ) 
\\ \notag & = \int_{\Rd \ts \sSS }
 - \half \DsftP \varphi 
+ \half \sum_{ q = 1}^{d } \DsftQ \psilaP \DsftQ \varphi 
+ \DDD [ \psilaP , \varphi ] \, dx \muz (d\sss ) 
\\ & \notag 
= - \int_{\Rd \ts \sSS } \la \psilaP \varphi \, dx \muz (d\sss ) 
\quad \text{ by \eqref{:81c} and \eqref{:82z}}
\\ & \label{:91e} = -
\lim_{t\downarrow 0} \frac{1}{t}
E [ \int_0^t \la \psilaP (\YY _u ) du ] 
\quad \text{ for $ ( X _0, \YY _0) \elaw \varphi dx \times \muz $}
.\end{align}
Using \eqref{:91d} and \eqref{:91e}, we apply Theorem 5.2.4 in \cite{fot.2} to obtain 
\begin{align}\label{:91f}&
N_t^{[u]} = \int_0^t \la \psila (\YY _u ) du 
.\end{align}
Hence from \eqref{:91y}, \eqref{:91d}, and \eqref{:91f}, we deduce 
\begin{align}\label{:91g}&
 \mathit{M}_t^{[u]} = A_t^{[u]} - N_t^{[u]} = \Mla ( t ) 
.\end{align}

Recall that $ u_p = x_p \otimes 1 + 1 \otimes \psilaP $. 
From \eqref{:61y} and a direct calculation similar to \eqref{:91e} with $ \varphi = u_p $, we have 
\begin{align} \notag 
\langle M ^{[u_p]} \rangle _t &= 2 \int_0^t \{
 (\nabla - \Dsft ) [ u_p , u_p ] + \DDD [u_p , u_p ]\} ( X _u , \YY _u ) du 
\\ \notag &= 2 \int_0^t \{
 (\nabla - \Dsft ) [ x_p \otimes 1 + 1 \otimes \psilaP , x_p \otimes 1 + 1 \otimes \psilaP ] 
\\\notag &\quad \quad \quad \quad + 
\DDD[ x_p \otimes 1 + 1 \otimes \psilaP , x_p \otimes 1 + 1 \otimes \psilaP ] \} ( X _u , \YY _u ) du 
\\\notag & = 
2 \int_0^t \{
\sumQ \half 
\Big\lvert \DsftQ \psilaP - \deltaPQ \Big\rvert ^2 + \DDD [\psilaP , \psilaP ] 
\} 
 ( X _u , \YY _u ) du 
\\\label{:91i} & = 
2 \int_0^t \{
\sumQ \half 
\Big\lvert \DsftQ \psilaP - \deltaPQ \Big\rvert ^2 + \DDD [\psilaP , \psilaP ] 
\} 
 ( \YY _u ) du 
.\end{align}

Recall that $ ( X _0 , \YY _0 ) \elaw \2 $ and $ \YY $ is a $ \muz $-stationary diffusion process 
by \lref{l:71}. Using these, \eqref{:91g}, and \eqref{:91i}, we have 
\begin{align}\label{:91n}
E[\lvert \Mlap (t) \rvert ^2 ] &= 2 t \Big\{ \int_{\sSS } 
\half \sumQ 
\Big\lvert \DsftQ \psilaP - \deltaPQ \Big\rvert ^2 \, d\muz + 
 \EYtwo ( \psilaP , \psilaP ) \Big\} 
.\end{align}
Hence, $ E[\lvert \Mlap (t) \rvert ^2 ] < \infty $ by $ \psilaP \in \dYL $. 
Thus, $ \Mlap $ is a continuous, $ L^2$-martingale with stationary increments. 
Applying \lref{l:83} to \eqref{:91n} yields \eqref{:91b}. 
\PFEND

\begin{theorem} \label{l:93}
Assume \As{A1}--\As{A4} and \As{A7}. 
Let $ ( X , \YY )$ be as in \lref{l:61}. 
Then, 
\begin{align}\label{:93a}&
\limz{\epsilon } \epsilon X _{t/ \epsilon ^2 } = 0 
\quad \text{weakly in $ C([0,\infty); \Rd )$ in $ \muz $-measure}
.\end{align}
\end{theorem}
\PF 
Using \eqref{:91y}, we see that for each $ t $ 
\begin{align}\label{:93g}
\epsilon X _{t/\epsilon ^2} - \epsilon X _{0 }& = \epsilon \M _{\epsilon ^2} (t/\epsilon ^2) 
 + \psieW (\YY _{t/ \epsilon ^2}) - \psieW (\YY _{0 } ) - 
 \epsilon \int_{0 }^{t / \epsilon ^2} \epsilon ^2 \psi _{\epsilon ^2} (\YY _u ) du 
.\end{align}
From \lref{l:91}, $ \M _{\epsilon ^2} $ is a continuous, $ L^2$-martingale with stationarity increments. 
Furthermore, $ \M _{\epsilon ^2} $ satisfies \eqref{:91b}. Hence, we find that 
\begin{align}\label{:93h}&
\limz{\epsilon }E [ \lvert \epsilon \M _{\epsilon ^2} (t/\epsilon ^2) \rvert ^2] = \limz{\epsilon }E [ \lvert \M _{\epsilon ^2} (t) \rvert ^2] = 0 
.\end{align}
From the $ \muz $-stationarity of $ \YY $ and \eqref{:82b}, we deduce that 
\begin{align} \notag 
\limz{\epsilon } E[ \lvert \psie (\YY _{t/ \epsilon ^2}) \rvert ^2 ] =& \, \limz{\epsilon } E[ \lvert \psie (\YY _{0}) \rvert ^2] 
\quad \text{ by stationarity of $ \YY $}
\\ \label{:93i} =& \, 
\limz{\epsilon } \epsilon ^2 \int_{\sSS } \lvert \psi _{\epsilon ^2} \rvert ^2 d\muz = 0 
\quad \text{ by \eqref{:82b}}
.\end{align}
Similarly, using the Schwartz inequality, the $ \muz $-stationarity of $ \YY $, and \eqref{:82b}, we obtain 
\begin{align}\notag 
\limz{\epsilon } E [\Big\lvert \epsilon \int_{0 }^{t / \epsilon ^2} \epsilon ^2 \psi _{\epsilon ^2} (\YY _u ) du \Big\rvert ^2 ] \le & \, 
\limz{\epsilon } E [\epsilon ^6 \frac{t}{\epsilon ^2} \int_{0 }^{t / \epsilon ^2} \lvert \psi _{\epsilon ^2} (\YY _u ) \rvert ^2 du ] 
\\ = & \, t^2 \limz{\epsilon }
\epsilon ^2 \int_{\sSS } \lvert \psi _{\epsilon ^2} \rvert ^2 d\muz = 0 
\quad \text{ by \eqref{:82b}}
 \label{:93f} 
.\end{align}
Clearly, $ \limz{\epsilon } E [ \lvert \epsilon X _{0}\rvert ^2] = 0 $. 
Hence, putting \eqref{:93h}--\eqref{:93f} into \eqref{:93g}, we obtain 
\begin{align}\label{:93j}&
\limz{\epsilon } E [\lvert \epsilon X _{t/ \epsilon ^2 }\rvert ^2 ] = 0 \quad \text{ for all }t 
.\end{align}
From \eqref{:93j}, we obtain 
\begin{align}\label{:93m}&
\limz{\epsilon } \epsilon X _{t/ \epsilon ^2 } = 0 
\quad \text{ in f.d.d.\,in $ \muz $-measure}
.\end{align}

Set $ X_{t} = (X_{ p , t})_{ p =1}^d$ and consider 
$ A _t^{[x_p \otimes 1] } ( (X , \YY ) ) = X_{ p , t}-X_{ p , 0}$. 
Applying Fukushima's decomposition to $ x_p \otimes 1 $, we have 
$ A _t^{[x_p \otimes 1] } = \mathit{M}_t^{[x_p \otimes 1]} + N _t ^{[x_p \otimes 1]}$. 
Then 
\begin{align} \notag 
\langle M ^{[x_p \otimes 1]} \rangle _t &= 
2 \int_0^t \{
 (\nabla - \Dsft ) [ x_p \otimes 1, x_p \otimes 1 ] + 
\DDD [x_p \otimes 1 ,x_p \otimes 1 ]\} ( X _u , \YY _u ) 
\} du 
\\\ &= \label{:93n}
t 
.\end{align}
Thus, $ M ^{[x_p \otimes 1]} $ is the standard linear Brownian motion. 

Applying Lyons--Zheng's decomposition \cite{c-f} to $ x_p \otimes 1 $, we see that 
$ A _t^{[x_p \otimes 1] } $ satisfies 
\begin{align}\label{:93o}&
A _t^{[x_p \otimes 1] } = \half \Big\{ \mathit{M}_t^{[x_p \otimes 1]} - 
 \mathit{M}_{t}^{[x_p \otimes 1]}\circ r_t 
 \Big\} 
,\end{align}
where $ r_t $ is the time-reversal operator on the path space on $ [0,\infty)$ such that 
$ r_t (\omega)(s) = \omega (t-s)$ if $ 0\le s \le t $ and  $ r_t (\omega)(s) = \omega (0)$ if $ t \le s  $.  
From \eqref{:93n} and \eqref{:93o}, 
$ \{ \epsilon X _{t/ \epsilon ^2 } \} $ is tight in $ C([0,\infty);\Rd )$. 
Combining this with \eqref{:93m} yields \eqref{:93a}. 
\PFEND

\section{Dual reduced Palm measures and mean-rigid conditioning}\label{s:D}

\subsection{Dual reduced Palm measures} \label{s:X}

In \ssref{s:X}, we introduce the concept of the dual reduced Palm measures. We construct the translation invariant dual reduced Palm measures in \lref{l:X5} and \lref{l:X6}. 


Let $ \SSm = \{ \sss \in \sSS ; \sss (\Rd ) = m \} $ for $ m \in \Nz $ as before. 
Let $ \mutm = \check{\mu }^m \circ \ulab ^{-1}$ for $ m \in \mathbb{N}$, 
where $ \check{\mu }^m $ are defined by \eqref{:13s} and $ \ulab (\mathbf{x}) = \sum_i\delta_{x_i}$ for $ \mathbf{x} = (x_i)_i$. 
Note that $ \mutm $ are measures on $ \SSm $ but usually not probability measures. 
Let $ \mut = \sum_{m=0}^{\infty} \mutm $ be the measure on $ \{ \sss\in \sSS ; \sss (\Rd ) < \infty \} $, where 
$ \mut _0 = \delta_{\mathsf{0}}$ is degenerated to the zero measure $ \mathsf{0}$. 
Let $ \muxx $, $ \xx \in \SSm $, be the reduced Palm measure defined by \eqref{:13x}. 

Recall that $ \SSRm = \{ \sss \in \sSS ; \sss (\SR )= m \} $. 
Note that $ \SSmB \cap \SSRm \subset \piR (\sSS )$ and that $ \sss (\SRc ) = 0 $ for $ \sss \in \SSmB \cap \SSRm $. We repeatedly use this fact in \ssref{s:X}. 
\begin{lemma} \label{l:X1} Assume \As{A2}, \As{A3}, and \eqref{:13f}. 
Then the following hold. 
\begin{align}\label{:X1a} 
&\muR \approx \mutR 
,\\\label{:X1b}
& \muRm \approx \muxx \circ (\piRc )^{-1} 
\end{align}
for each $ \xx \in \SSmB \cap \SSRm $. Here, $ \cdot \approx * $ means $ \cdot $ and $ * $ are mutually  absolutely continuous. 
\end{lemma}
\PF 
From \As{A2}, we have the $ m $-point correlation function $ \rho ^m $ and the density function $ \sigma _R^m $ on $ \SR $ of $ \mu $ such that 
\begin{align}\label{:X1f}
\rho ^m (\mathbf{x}_m) = & 
\sigma _R^m (\mathbf{x}_{m}) +\sum_{n =m+1 }^{\infty}\frac{1}{(n-m)!}
 \int_{(\SR )^{n-m}}
\sigma _R^{n} ((\mathbf{x}_{m},\mathbf{y}_{n-m}) ) \prod_{l=1}^{n-m}dy_{l}
,\end{align}
where $ \mathbf{y}_{n-m}= (y_1,\ldots,y_{n-m}) \in (\Rd )^{n-m} $. 
From \As{A3}, \eqref{:12y} and \eqref{:12z} hold. 
Let $ \mu_{ \Rks }^n $ be the regular conditional measures given by \eqref{:12x}. 
Hence, for $ \mu $-a.s.\,$ \sss $, 
the Radon-Nikodym density $ \sigma_{ \Rks }^n = d\mu_{ \Rks }^n /d\Lambda _{\rR }^n $ 
and the $ n $-point correlation function $ \rho _{ \Rks }^n $ of $ \mu_{ \Rks }^n $ exists. 
Furthermore, there exists a constant $ \Ct (\Rks , n )\label{;X1} > 0 $ such that for any $ \mathbf{x}_n=(x_1,\ldots,x_n) \in \SR ^n $ 
\begin{align*}&
\cref{;X1}^{-1} \exp ( - \beta \sum_{ x_p, x_q \in \SR \atop p < q }^n \Psi (x_p-x_q) ) \le 
 \sigma_{ \Rks }^n (\mathbf{x}_{n}) \le 
\cref{;X1} \exp ( - \beta \sum_{ x_p, x_q \in \SR \atop p < q }^n \Psi (x_p-x_q) )
.\end{align*}
Hence, \eqref{:X1f} also holds for $ \sigma _{ \Rks }^n $ and $ \rho _{ \Rks }^n $.

From \As{A3}, $ \Psi $ is locally bounded from below, and $\Psi ( x ) < \infty $ for $ x \ne 0 $. Using these and that $ \sigma _{\Rks }^n (\mathbf{x}_{n}) $ is symmetric in $ \mathbf{x}_n$, 
we see for $ \mu $-a.s.\,$ \sss $ 
\begin{align*}&
\sigma _{\Rks }^m (\mathbf{x}_{m}) = 0 
\Longleftrightarrow 
\text{$ \sigma _{\Rks }^{n} ((\mathbf{x}_{m},\mathbf{y}_{n-m}) ) = 0 $ for all 
$ \mathbf{y}_{n-m} \in \SR ^{n-m} $, $ n >m $}
.\end{align*}
Hence using this and \eqref{:X1f} for $ \sigma _{ \Rks }^n $ and $ \rho _{ \Rks }^n $, we deduce 
for $ \mu $-a.s.\,$ \sss $ 
\begin{align}\label{:X1g}&
\sigma _{\Rks }^m (\mathbf{x}_{m}) = 0 
\Longleftrightarrow 
 \rho _{\Rks }^m (\mathbf{x}_m) = 0 
.\end{align}
Integrating \eqref{:X1g} for $ \sss $ with respect to $ \mu $ and taking $ k \to \infty $, we obtain from \eqref{:12x} 
\begin{align}\label{:X1h}&
\sigma _{R }^m (\mathbf{x}_{m}) = 0 \Longleftrightarrow 
 \rho ^m (\mathbf{x}_m) = 0 
\end{align}
We deduce from \eqref{:X1f} and \eqref{:X1h} that $ \muR \approx \mutR $, which yields \eqref{:X1a}. 

For $ \mathsf{A} \in \mathcal{B}(\sSS ) $, we set $ \ARc = (\piRc )^{-1}( \piRc (\mathsf{A})) $. Then for $ \xx \in \SSm\cap \SSRm $, 
\begin{align}
&&
\muRmA 
& \notag =
\sum_{n=m}^{\infty}
\int_{\SSRn } \mu ( \ARc \vert \mathsf{t} = \piR (\sss ) ) \muR (d\mathsf{t} )
\\&&&\notag \approx 
\int_{\SSRm } \mu ( \ARc \vert \mathsf{t} \prec \piR (\sss ) ) \mutR (d\mathsf{t} )
&& \text{by \eqref{:X1a}}
\\&&&\notag = 
\int_{\SSRm } \mu ( \ARc + \mathsf{t} \vert \mathsf{t} \prec \piR (\sss ) ) \mutR (d\mathsf{t} ) 
\\&&&\notag = 
\int_{\SSRm } \mu ( \ARc \Vert \mathsf{t} ) \mutR (d\mathsf{t} ) && \text{by }\eqref{:13x} 
\\ &&& \notag 
\approx \int_{\SSRm } \mu ( \ARc \Vert \xx ) \muR (d\mathsf{t} ) 
\approx \mu ( \ARc \Vert \xx ) 
&& \text{by \eqref{:13f}}
.\end{align}
This yields \eqref{:X1b}. 
\PFEND

We note that $ \sSS $ is homeomorphic to a complete separable metric space. 
Hence, for $ \mu $-a.s.\,$ \yy $, we define a random point field $ \muRyy $ on $ \SR $ by 
the regular conditional probability such that 
\begin{align*}&
 \muRyy ( \cdot )= \mu ( \piR (\sss ) \in \cdot \vert \piRc (\sss ) = \piRcy ) 
. \end{align*}
Note that $ \mathsf{A} = \piR (\mathsf{A}) + \piRc (\mathsf{A})$ and that 
$ \muRyy = \mu _{\rR , \piRc (\yy ) }$. Then we have 
\begin{align}\notag 
\mu (\mathsf{A}) &=
\int_{(\piRc )^{-1} (\piRc (\mathsf{A} )) } \muRyy (\piR (\mathsf{A})) \mu (d\yy ) = 
\int_{\piRc (\mathsf{A} ) } \, \muRyy (\mathsf{A}) \muRc (d\yy ) 
.\end{align}
In \lref{l:X2}, $ \mu $ is irreducibly $ k $-decomposable with $ \SSmDDk $ in the sense of \dref{d:D1}. 
%
From \rref{r:14}\thetag{2}, $ \mu (\SSmD) = 0 $ for $ m \ge 1 $. 
Thus, it is not obvious that $ \muRyy $ exists for $ \yy \in \SSmD $, $ m \ge 1 $. 
We resolve this in \lref{l:X2}. 

Let $ \mutmD $ be a random point field such that 
\begin{align}&\label{:X2w}
\text{$ \mutmD \approx \muxx $ for some $ \xx \in \SSmB $}
.\end{align} 
From \eqref{:13a}--\eqref{:13f}, we have 
\begin{align}&\label{:X2x}
 \mutmD (\SSmD ) = 1
,\\&\label{:X2y}
\mutmD \approx \muxx \quad \text{ for all } \xx \in \SSmB 
.\end{align}
Note that 
$ \SSRm + \piR (\yy ) = \{ \sss + \piR (\yy ) ; \sss \in \SSRm \} \subsetneqq \SSRm $ if $ \yy (\SR ) \ge 1 $. 
We set 
\begin{align}&\label{:13e}
\SSzDfyy = \{ \sss \in \SSzD ; \yy \prec \sss \} 
.\end{align}
If $ \mu $ is $ k $-decomposable, then $ \SSzDfyy - \yy \subset \SSm $ for $ \yy \in \SSmD $ from \eqref{:13b}. 

\begin{lemma} \label{l:X2}
Assume \As{A2}, \As{A3}, and \eqref{:13f}. 
Let $ \mu $ be irreducibly $ k $-decomposable with $ \SSmDDk $ in the sense of \dref{d:D1}. 
Then for $ 1 \le m \le k $, 
\begin{align}\label{:X2a}&
\text{$ \muRyy $ exists for \yas}
,\\\label{:X2b}&
\text{$ \muRyy (\SSRm + \piR (\yy )) = 1 $ for \yas}
.\end{align}
\end{lemma}
\PF 
Using \eqref{:X2y} and \eqref{:X1b}, we have for all $ \xx \in \SSmB $ such that $ \xx (\SR ) = m $
\begin{align} \notag 
 \mutmD \circ (\piRc )^{-1}& \approx \muxx \circ (\piRc )^{-1} 
\\ \label{:X2h}& \approx \muRm \ll \muRc 
.\end{align}
From $ \muRyy = \mu _{\rR , \piRcy }$, $ \muRyy $ exists for $ \muRc $-a.s.\,$ \yy $. 
Combining this with \eqref{:X2h}, we obtain \eqref{:X2a}. 
From \eqref{:13b}, $ \SSzDfyy - \yy \subset \SSm $ for $ \yy \in \SSmD $. 
Note that we condition $ \mu $ on $ \SRc $ as $ \muRyy = \mu (\cdot \vert \piRc (\sss ) = \piRcy ) $. 
Hence, 
$ \SSzDfyy - \yy \subset \SSRm $ for $ \mutmD $-a.s.\,$ \yy \in \SSmD $. 
We have thus obtained \eqref{:X2b}. 
\PFEND

By \eqref{:X2a}, $ \muRyy $ exists for \yas\ for $ 1 \le m \le k $. 
Let $ \rho_{\rR , \yy}^m $ be the $ m $-point correlation function of $ \muRyy $. 
Then $ \rho_{\rR , \yy}^m $ exists on $ \SR $ and satisfies 
\begin{align}\label{:X3a}&
\int_{\sSS }\rho_{\rR , \yy}^m (x_1,\ldots,x_m) \muRc (d\yy ) = \rho ^m (x_1,\ldots,x_m) 
\quad \text{ on } \SRm 
.\end{align}
From \eqref{:X3a}, the $ m $th factorial moment measure $ \mut _{\rR , \yy}^m $ of $ \muRyy $ exists on $ \SR $ for \yas. Hence, the reduced Palm measure $ \muRyy (\cdot \Vert \xx ) $ of $ \muRyy $ exists for $ \mut _{\rR , \yy}^m $-a.e.\,$ \xx $ for \yas. 
Let $ \SSR $ be the configuration space over $ \SR $. We regard $ \SSRm $ as a subset of $ \SSR $. 
Then by definition, for \yas, 
\begin{align}& \notag 
\int_{\mathsf{A}} \muRyy ( \mathsf{B} \Vert \xx ) \mut _{\rR , \yy}^m (d\xx ) = 
\int_{\mathsf{A}}
	 \muRyy (\mathsf{B}+\xx \vert \xx \prec \sss ) \mut _{\rR , \yy}^m (d\xx ) 
\end{align}
for all $\mathsf{A} \in \mathcal{B}(\SSRm )$, and $ \mathsf{B} \in \mathcal{B}(\SSR ) $, 
and, equivalently, for $ \mut _{\rR , \yy}^m $-a.e.\,$ \xx $, 
\begin{align}\label{:X3y}&
 \muRyy ( \cdot \Vert \xx ) = \muRyy (\cdot +\xx \vert \xx \prec \sss ) 
.\end{align}
Originally, $ \muRyy (\cdot \Vert \xx ) $ is the random point field on $ \SR $. 
We can regard $ \muRyy ( \cdot \Vert \xx )$ as a random point field on $ \Rd $ by taking 
\begin{align}\label{:X3z}&
 \muRyy ( \cdot \Vert \xx ) \circ (\piRc )^{-1}= \delta_{\piRc(\yy )}
.\end{align}
We denote this extension by the same symbol $ \muRyy (\cdot \Vert \xx ) $. 
By definition, $ \muRyy (\cdot \Vert \xx ) $ degenerates into $ \piRc (\yy )$ outside $ \SR $. 
In the following, $ \muRyy (\cdot \Vert \xx ) $ is thus a random point field on $ \Rd $. 
Let $ \SSyy = \{ \sss \in \sSS ; \yy \prec \sss \} $ and $ \SSzDfyy = \SSyy \cap \SSzD $. 
\begin{lemma} \label{l:X3}
Make the same assumptions as \lref{l:X2}. Then for $ 1 \le m \le k $ 
\begin{align}\label{:X3f}&
 \muRyy ( \SSRm + \piRcy \Vert \piRy ) = 1 \quad \text{ for \yas}
.\end{align}
\end{lemma}
\PF 
The reduced Palm measure $ \muRyy (\cdot \Vert \xx ) $ satisfies for $ 1 \le m \le k $
\begin{align}\label{:X3b}&
\text{$ \muRyy ( \sSS (\xx + \piRc (\yy ) ) - \xx \Vert \xx ) =1 $ for $ \mut _{\rR , \yy}^m $-a.e.\,$ \xx $}
.\end{align}
Similarly as \eqref{:X1h}, we deduce 
$ \sigma_{\rR , \yy}^m (\mathbf{x}_{m}) = 0 \Longleftrightarrow \rho_{\rR , \yy}^m (\mathbf{x}_m) = 0 $. 
Hence, 
\begin{align}\label{:X3c}&
 \muRyy \approx \mut _{\rR , \yy} \quad \text{ for \yas} 
.\end{align}
Here $ \mut _{\rR , \yy} = \sum_{ n =1}^{\infty} \mut _{\rR , \yy}^{ n }$. 
Taking $ \xx = \piRy $ in \eqref{:X3b} and using \eqref{:X3b} and \eqref{:X3c}, we see that the reduced Palm measure $ \muRyy (\cdot \Vert \piRy ) $ satisfies 
\begin{align}&\label{:X3e}
 \muRyy ( \SSyy - \piRy \Vert \piRy ) =1 \quad \text{ for \yas}
.\end{align}
Using \eqref{:X3y}, we can replace $ \SSyy $ in \eqref{:X3e} by $ \SSzDyy =\SSyy \cap \SSzD $. 
Thus we obtain 
\begin{align}\label{:X3g}&
\text{ $ \muRyy (\SSzDfyy - \piRy \Vert \piRy ) = 1 $ for \yas}
.\end{align}
Note that $ \muRyy = \mu _{R , \piRc (\yy )}$. 
By \eqref{:X3z}, all removed particles of $ \muRyy ( \cdot \Vert \xx ) \circ (\piRc )^{-1} $ 
are in $ \SR $, that is, $ \xx (\SR ) = m $. Hence, from \eqref{:13b} and \eqref{:X3z}, we have for \yas 
\begin{align}\label{:X3i}& 
 \SSzDfyy - \yy \subset \SSRm \quad \text{ under } \muRyy 
.\end{align}
Hence from \eqref{:X3g} and \eqref{:X3i}, we obtain 
\begin{align}\notag & 
1 = \muRyy ( \SSzDfyy - \yy \Vert \piRy ) \le \muRyy ( \SSRm \Vert \piRy ) 
\quad \text{for \yas}
.\end{align}
This implies \eqref{:X3f}. 
\PFEND

For $ \rR \in \mathbb{N}$ and $ \yy \in \sSS $, let $ \IR $ and $ \Iinfty $ be the $ \sigma $-fields 
such that 
\begin{align} \label{:X9q}
\IR & = \sigma [\SSzD (\piR (\yy ))] \vee \sigma [\piRc ] ,\quad 
\Iinfty = \bigcap_{\rR = 1}^{\infty} \IR 
.\end{align}

\begin{lemma} \label{l:X9} 
Make the same assumptions as \lref{l:X2}. 
Let $ \nuR{\cdot } $ be the regular conditional probability measure with respect to $ \IR $ 
for $ \rR \in \mathbb{N}\cup \{ \infty \} $ and $ \yy \in \SSmD $. 
Then the following hold for $ 1 \le m \le k $ and \yas. 
\\
\thetag{1} For each $ \mathsf{B} \in \mathcal{B}(\sSS ) $, 
\begin{align}
\label{:X9z}&
 \nui{\mathsf{B} }(\sss ) = \limi{\rR } \nuR{\mathsf{B} } (\sss ) 
\text{ for $ \mu $-a.s.\,$ \sss $ and in $ \Lone $}
.\end{align}
\thetag{2}
For $ \mu $-a.s.\,$ \sss $, 
\begin{align}\label{:X9y}&
 \nui{\mathsf{B} } (\sss ) = \limi{\rR } \nuR{\mathsf{B} } (\sss ) \quad \text{ for all }
 \mathsf{B} \in \mathcal{B}(\sSS )
.\end{align}
\thetag{3} 
Let $ \mathsf{B} \in \mathcal{B}(\sSS ) $ be such that $ \mathsf{B} \subset \SSmB \cap \SSRm $ for some $ \rR \in \mathbb{N}$. Then $ \nui{\mathsf{B} } (\sss ) $ is constant and 
$  \nui{ \SSm + \yy } (\sss ) = 1 $  for $ \mu $-a.s.\,$ \sss $. 
\end{lemma}
\PF 
We easily see that $ \IR \supset \IRR $ for $ \rR \le \rR ' \le \infty $ and that 
$ \{ \mu (\mathsf{B} \vert \IR ) \} $ is bounded in $ \Lm $. 
Hence using the martingale convergence theorem (cf. \cite[I (2.4) Corollary]{RY}), we have \eqref{:X9z}, which implies \thetag{1}. 

From \thetag{1}, it is easy to see that $\{ \nuR{\cdot } (\sss ) \} _{ \rR \in \mathbb{N}} $ is tight for $ \mu $-a.s.\,$ \sss $. 
Hence, we denote an arbitrary convergent subsequence by the same symbol $ \nuR{\cdot } (\sss ) $ and its limit by $ \mu ' (\cdot ) (\sss ) $. 
Note that the measurable space $ (\sSS , \mathcal{B}(\sSS ) ) $ is countably determined, that is, any probability measures on $ (\sSS , \mathcal{B}(\sSS ) )$ are determined by a countable system of elements of 
$ \mathcal{B}(\sSS ) $ (cf. \cite[p.14]{IW}). 
Let $ \{ \mathsf{B}_n \}_{n\in \mathbb{N}} $ be a countable system of subsets determines the probabilities on $ (\sSS , \mathcal{B}(\sSS ) )$. Then from \eqref{:X9z}, for $ \mu $-a.s.\,$ \sss $, 
\begin{align}\label{:X9i}&
\mu ' (\mathsf{B}_n) (\sss ) = \limi{\rR } \nuR{\mathsf{B}_n} (\sss ) = \nui{\mathsf{B}_n} (\sss ) 
\quad \text{ for all } n \in \mathbb{N}
.\end{align}
Hence from \eqref{:X9i}, we deduce $ \mu ' (\cdot ) (\sss ) = \nui{\cdot } (\sss ) $ for $ \mu $-a.s.\,$ \sss $. This implies \thetag{2}.

Let $ T _1 = \sS _{1} $ and $ \TQ = \SQ \backslash \sS _{{\qQ -1}} $ for $ \qQ \ge 2 $. 
For $ m $ and $ \rR $, let 
\begin{align*}&
 \NNNR = \{ (\nQ )_{\qQ = 1}^{\rR }; 0\le \nQ \le m , 
\sum_{\qQ =1}^{\rR } \nQ = m \} 
,\\&
 \TTRyn = \bigcap_{\qQ = 1}^{\rR } 
 \{ \sss \in \SSzDfyy ; \sss (\TQ ) = \yy (\TQ ) +\nQ \}
, \ \nnn =(\nQ )_{\qQ = 1}^{\rR } \in \NNNR 
,\\&
\UURy = \bigcup_{ \nnn \in \NNNR } \TTRyn 
.\end{align*}
Then $ \TTRyn \cap \TTRynn = \emptyset $ for $ \nnn \ne \nnn '$. Thus, 
$ \{ \TTRyn \}_{\nnn \in \NNNR } $ is a partition of $ \UURy $. 
This together with $ \mathsf{B} \subset \SSmB \cap \SSRm $ yields $ \mathsf{B} + \yy \subset \UURy $ and 
\begin{align}\label{:X9b}&
 \nuR{\mathsf{B} + \yy } (\sss ) = 
\sum_{\nnn \in \NNNR } \nuR{( \mathsf{B} + \yy ) \cap \TTRyn } (\sss ) 
.\end{align}
Using \eqref{:X9q} and $ \mathsf{B} \subset \SSmB \cap \SSRm $, we have 
\begin{align} &\label{:X9!}&
 \nuR{ ( \mathsf{B} + \yy ) \cap \TTRyn } (\sss ) 
= 
\begin{cases}
1 & \text{ for $ \mu $-a.e.\,} \sss \in ( \mathsf{B} + \yy ) \cap \TTRyn 
\\
0 & \text{ for $ \mu $-a.e.\,} \sss \notin ( \mathsf{B} + \yy ) \cap \TTRyn 
\end{cases}
.\end{align}
From \eqref{:X9b}, \eqref{:X9!}, $ \mathsf{B} + \yy \subset \UURy $, and that 
$ \{ \TTRyn \}_{\nnn \in \NNNR } $ is a partition of $ \UURy $, we see 
\begin{align}\notag 
 \nuR{ \mathsf{B} + \yy } (\sss ) &= 1 \quad \text{ for $ \mu $-a.e.\,} \sss \in \UURy 
.\end{align}
Using this, \eqref{:X9y}, and that $ \{\UURy \}_{\rR =1}^{\infty} $ is increasing, we see 
\begin{align} \label{:X9d} 
 \nui{\mathsf{B} + \yy } (\sss ) &= \limi{\rR} \nuR{ \mathsf{B} + \yy } (\sss ) 
= 1 \quad \text{ for $ \mu $-a.e.\,} 
\sss \in \bigcup_{\rR =1}^{\infty}\UURy 
.\end{align}
From \eqref{:X9d} and $ \mu (\cup_{\rR =1}^{\infty}\UURy ) = 1 $, we obtain the first claim in \thetag{3}. 
Using \eqref{:X3y} and \lref{l:X3}, we see for \yas\, and any 
$ \rR \in \mathbb{N}$, 
\begin{align} &\notag 
 \nuR{ \SSm + \yy } (\sss ) = 1 \quad \text{ for $ \mu $-a.e.\,} \sss \in \UURy 
.\end{align}
From this, \eqref{:X9d}, and $ \mu (\cup_{\rR =1}^{\infty}\UURy ) = 1 $, we get the second claim in \thetag{3}. 
\PFEND

We now construct the dual reduced Palm measure conditioned on infinite-many points $ \yy \in \SSmD $. For \yas\, and $ \mathsf{B} \in \mathcal{B}(\SSm ) $, let 
\begin{align}\label{:X9h}&
 \mu ( \mathsf{B} \Vert \yy ) := \nui{\mathsf{B} + \yy } 
.\end{align}
From \lref{l:X9}\thetag{3}, $ \nui{\mathsf{B} + \yy } (\sss )$ is constant and 
$ \nui{\SSm + \yy } = 1 $ for $ \mu $-a.s.\,$ \sss $. 
Thus $ \mu ( \mathsf{B} \Vert \yy ) $ is independent of $ \sss $ and satisfies $ \mu ( \SSm \Vert \yy ) = 1$. 
Hence, we extend the domain of $ \mu ( \cdot \Vert \yy ) $ to $ \mathcal{B}(\sSS ) $ 
by $ \mu ( \cdot \Vert \yy ) = \mu ( \cdot \cap \SSm \Vert \yy )$. 

We regard \eqref{:X4z} and \eqref{:X4x} in \lref{l:X4} as the dual relation to $ \mu ( \SSmD \Vert \yy ) = 1 $ for $ \yy \in \SSm $ in \eqref{:13c} for the original reduced Palm measure. 
This relation is a result of decomposability of $ \mu $. 
We call $ \mu ( \cdot \Vert \yy ) $ the dual reduced Palm measure because $ \mu ( \cdot \Vert \yy ) $ is conditioned at an element $ \yy $ of the dual set $ \SSmD $ to $ \SSm $ and $ \mu ( \cdot \Vert \yy ) $ is supported on $ \SSm $ such that $ \mu ( \SSm \Vert \yy ) =1 $ for $ \yy \in \SSmD $. 

We convert the convergence of the Palm measures in \lref{l:X9} to that of the reduced Palm measures. 
\begin{lemma} \label{l:X4} 
Make the same assumptions as \lref{l:X2}. Let $ 1 \le m \le k $. 
\\
\thetag{1} 
For \yas, 
\begin{align}
\label{:X4z}&
\mu ( \mathsf{B} \Vert \yy ) = \limi{\rR } \muRyy ( \mathsf{B} \Vert \piRy ) 
\quad \text{ for all } \mathsf{B} \in \mathcal{B}(\sSS ) 
.\end{align}
\thetag{2} 
$ \mu ( \mathsf{B} \Vert \cdot ) $ is a $ \mathcal{B}(\sSS )$-measurable function for each $ \mathsf{B} \in \mathcal{B}(\sSS ) $. 
\\\thetag{3} For \yas, $ \mu ( \cdot \Vert \yy ) $ is a probability measure on $ (\sSS ,\mathcal{B}(\sSS ) ) $ such that 
\begin{align}\label{:X4x}& 
\mu ( \SSm \Vert \yy ) = 1 
.\end{align}
\end{lemma}
\PF
Because $ \mu ( \cdot \Vert \yy ) = \mu ( \cdot \cap \SSm \Vert \yy )$ by definition, we assume 
$ \mathsf{B} \subset \SSm $. We see 
\begin{align}\label{:X4a}
 \mu ( \mathsf{B} \Vert \yy ) &= \nui{\mathsf{B} + \yy } \quad \text{ by \eqref{:X9h}}
\\\notag & 
= \limi{\rR } \nuR{ \mathsf{B} + \yy}) (\sss ) 
\text{ by \lref{l:X9}}
\\\notag & 
= \limi{\rR } \nuR{ \piR(\mathsf{B} ) + \yy } (\sss ) 
\text{ by $ \mathsf{B} \subset \SSm $}
.\end{align}
Here the convergence takes place for $ \mu $-a.s.\,$ \sss $ and in $ \Lone $. 
From \eqref{:X9q} 
\begin{align}\notag 
& \mu ( \piR(\mathsf{B} ) + \yy \vert \IR ) (\sss ), \quad \piRc (\sss ) = \piRc (\yy )
, \\&\notag 
= \muRyy ( \piR(\mathsf{B} ) + \piR (\yy ) \vert \piR (\yy ) \prec \sss ) 
\\\label{:X4b}&
= \muRyy ( \piR(\mathsf{B} ) \Vert \piR (\yy ) \prec \sss ) 
.\end{align}
From \eqref{:X4a}, \eqref{:X4b}, and $ \mathsf{B} \subset \SSm $, we see \eqref{:X4z}. 
From \eqref{:X3y} and \eqref{:X4z}, we obtain \thetag{2}. 
From \lref{l:X9}\thetag{3} and \eqref{:X9h}, we have \thetag{3}. 
 \PFEND

\begin{lemma}	\label{l:X5} 
Make the same assumptions as \lref{l:X2}. 
In addition, assume that $ \mu $ and $ \SSzD $ are translation invariant. 
Then, for \yas\, and $ 1 \le m \le k $, 
\begin{align} & \notag 
\muyy = \mu ( \cdot \Vert \vartheta _{ \aaa } ( \yy ) ) \circ \vartheta_{ \aaa }^{-1} 
\quad \text{ for all } \aaa \in \Rd 
.\end{align}
\end{lemma}
\PF 
Because $ \SSzD $ is translation invariant, we have 
\begin{align}\label{:X5f}&
 \SSzDP \vartheta_{ \aaa } ( \yy )) = \vartheta_{ \aaa } (\SSzDyy )
.\end{align}
Then using \eqref{:X4z}, \eqref{:X5f}, and \As{A1}, 
we have, for \yas\ and all $ \aaa \in \Rd $, 
\begin{align}
\notag 
\mu ( \cdot \Vert \yy ) &
=
 \limi{\rR } \muRyy (\cdot \Vert \piRy ) \quad \text{ by \eqref{:X4z}}
\\\notag & =
 \limi{\rR } \muRyy \big( 
\vartheta_{ \aaa } ^{-1} ( \vartheta_{ \aaa } (\cdot ) ) 
 \Vert 
\piR ( 
\vartheta_{ \aaa } ^{-1} ( \vartheta_{ \aaa } ( \yy ) ) )
 \big)
\\\notag & =
 \limi{\rR } \mu _{\rR , \vartheta_{ \aaa } (\yy ) } 
\big( 
 \vartheta_{ \aaa } (\cdot ) 
 \Vert 
\piR ( \vartheta_{ \aaa } (\yy ) ) 
\big) 
\quad \text{ by \eqref{:X5f} and \As{A1}}
\\&\notag 
= 
\mu ( \cdot \Vert \vartheta_{ \aaa }(\yy ) ) \circ \vartheta_{ \aaa }^{-1} \quad \text{ by \eqref{:X4z}}
.\end{align}
We have thus completed the proof of \lref{l:X5}. 
\PFEND

\begin{lemma} \label{l:X6} 
Assume \As{A1}--\As{A3} and \As{A6}. Then \As{A5} holds with $ \{ \SSzD , \SSoneD \} $ such that both $ \SSzD $ and $\SSoneD $ are translation invariant. 
\end{lemma}
\PF 
Let $ \Ke \subset \sSS $ be an increasing sequence of compact sets such that $ \mu (\Ke ) > 1-\epsilon $ as $ \epsilon \downarrow 0 $. We set 
\begin{align}& \notag 
\KKKRe = \bigcup_{ \lvert x \rvert \le \rR } \vartheta_x (\Ke ), \ 0 \le \rR \le \infty ,
\quad \KKKI = \bigcup_{0 < \epsilon < 1 } \KKKIe 
.\end{align}
Note that the translation operator $ \vartheta_x $ on $ \sSS $ is a homeomorphism for each $ x \in \Rd $ and that $ \sss \mapsto \vartheta_x (\sss )$ is continuous in $ x \in \Rd $ for each $ \sss $. 
Hence, $ (x,\sss ) \mapsto \vartheta_x (\sss ) $ is continuous. 

Let $ \rR < \infty $. 
Because $ \{ \lvert x\rvert \le \rR \} \times \Ke $ is compact and the map 
$ (x,\sss ) \mapsto \vartheta_x (\sss ) $ is continuous, we deduce that $ \KKKRe $ is a compact set in $ \sSS $. 
Hence, $ \KKKI $ is a Borel set because 
$ \KKKI $ is the increasing limit of $ \KKKIe $ as $ \epsilon \downarrow 0 $ and 
$ \KKKIe $ is the increasing limit of compact sets $ \KKKRe $ as $ \rR \uparrow \infty $. 

From $ \mu (\cup_{ 0 < \epsilon < 1 } \Ke ) = 1 $ and $ \cup_{ 0 < \epsilon < 1 } \Ke \subset \KKKI $, 
we see $ \mu (\KKKI ) = 1 $. By construction, $ \KKKIe $ is translation invariant for each $ \epsilon $. 
Hence, $ \KKKI $ is translation invariant. 
Let $ \xx = \delta_x$ and $\KKKIxx = \{ \sss \in \KKKI ; \xx \prec \sss \} $. 
Then we have for $ \mutone $-a.e.\,$ \xx $ 
\begin{align}\label{:X6o}&
 \mu (\KKKIxx - \xx \Vert \xx ) = 1 
.\end{align}
Recall that $ \SSone = \{ \xx = \delta _x ; x \in \mathbb{R}^d \} $. 
Because $ \mu $, $ \KKKI $, and $ \SSone $ are translation invariant, 
we can refine the property above so that \eqref{:X6o} holds for \xasone. Let 
\begin{align}\label{:X6p}&
 \KKKIone = \bigcup_{\xx \in \SSone } \{\KKKIxx - \xx \}
.\end{align}
Then $ \KKKIone $ is translation invariant. 
Because $ \KKKIxx - \xx \subset \KKKIone $ and $ \mu (\KKKIxx - \xx \Vert \xx ) = 1$ for all $ \xx \in \SSone $, 
we obtain 
\begin{align}\label{:X6P}&
 \KKKIone \in \overline{\mathcal{B} (\sSS )}^{\muxx } \text{ and } \mu ( \KKKIone \Vert \xx ) = 1 
 \text{ for all } \xx \in \SSoneB 
\end{align}

Let $ \SSoneD = \KKKIone $ and $ \SSzD = (\KKKIone + \SSone ) \backslash \KKKIone $. 
Then $ \SSzD \cap \SSoneD= \emptyset $, and thus we obtain \eqref{:13a}. \eqref{:13b} follows from 
$ \SSzD = (\KKKIone + \SSone ) \backslash \KKKIone \subset \KKKIone + \SSone = \SSoneD + \SSone $. 

From \As{A6}\thetag{1} and \eqref{:X6P}, $ \mu ( \KKKIone ) = 0 $. 
From this and $ \mu (\KKKI ) = 1 $, we deduce 
\begin{align}\notag &&
\mu (\SSzD ) &\ge \mu (\KKKIone + \SSone ) - \mu ( \KKKIone ) 
&&\text{by } \SSzD = (\KKKIone + \SSone ) \backslash \KKKIone 
 \\\notag && &
= \mu (\KKKIone + \SSone ) && \text{by $ \mu ( \KKKIone ) = 0 $}
\\&&&\notag 
= \mu ( \bigcup_{\xx \in \SSone } \{\KKKIxx - \xx \} + \SSone ) &&\text{by \eqref{:X6p}}
\\&&&\notag 
\ge \mu ( \bigcup_{\xx \in \SSone } \KKKIxx ) &&\text{by } \{\KKKIxx - \xx \} + \SSone \supset \KKKIxx 
\\&&&= \mu (\KKKI ) = 1 &&\text{by } \bigcup_{\xx \in \SSone } \KKKIxx = \KKKI 
\label{:X6r}
.\end{align}
Hence, \eqref{:13c} follows from \eqref{:X6P} and \eqref{:X6r}. 
By construction, $ \SSzD $ and $ \SSoneD $ are translation invariant. 
Irreducibility follows from \As{A6}\thetag{2}, which completes the proof. 
\PFEND

\subsection{Mean-rigid conditioning}\label{s:J}
In \ssref{s:J}, we introduce the concept of the mean-rigid $ \sigma $-field $ \Gi $. 
We define the functions $ \NR $ and $ \MMR $ on $ \sSS $ such that 
\begin{align}\label{:J1m}&
\NR (\sss ) = \sss (\SR ) , \quad \MMR (\sss ) = \sum_{s_i \in \SR } s_i 
.\end{align}
Let $ T _1 = \sS _{1} $ and $ \TR = \SR \backslash \sS _{{\rrr -1}} $ for $ \rrr \ge 2 $. 
Replacing $ \SR $ by $ \TR $ in \eqref{:J1m}, we define $ \NTr $ and $ \MTr $. 
For a function $ f $, a $ \sigma $-field $ \mathcal{F} $, and a random point field $ \mu $, 
we say $ f $ is $ \mathcal{F} $-measurable for $ \mu $-a.s.\,if a $ \mu $-version of $ f $ is $ \mathcal{F} $-measurable. 
\begin{definition}\label{d:m2} 
\thetag{1} 
A random point field $ \mu $ on $ \Rd $ is said to be number-rigid if, for each $ \rR \in \mathbb{N} $, the function $ \NR (\sss ) $ (resp.\,$ \NTr $) is $ \sigma [\piRc ] $-measurable 
(resp.\,$ \sigma[\pi _{\TR }^c]$-measurable) 
for $ \mu $-a.s. 
\\\thetag{2}
A random point field $ \mu $ on $ \Rd $ is said to be mean-rigid if $ \mu $ is number-rigid and, for each $ \rR \in \mathbb{N} $, the function $ \MMR (\sss ) $ (resp.\,$ \mathrm{M}_{\TR } (\sss ) $) is $ \sigma [\piRc ] $-measurable (resp.\,$ \sigma[\pi _{\TR }^c]$-measurable) for $ \mu $-a.s.
\end{definition}
If $ \mu $ is mean-rigid, then, by definition, 
for each $ \rR \in \mathbb{N} $ and for $ \mu $-a.s.\,$ \sss $, there exist 
$ a = a (\piRc (\sss )) $, $ b = b (\piRc (\sss ))$, 
$ a' = a' (\pi _{\TR }^c (\sss )) $, and $ b' = b' (\pi _{\TR }^c (\sss ))$ such that 
\begin{align}\notag &
\mu ( \{ \sss \in \sSS ; \MMR (\sss ) = a (\piRc (\sss )) ,\, \NR (\sss ) = b (\piRc (\sss )) \} ) = 1 
,\\ \notag & 
\mu ( \{ \sss \in \sSS ; \MTr (\sss ) = a' (\pi _{\TR }^c (\sss )) ,\, 
\NTr (\sss ) = b' (\pi _{\TR }^c (\sss )) \} ) = 1 
.\end{align}
In \cite{gp,bu}, rigidity is posed for all Borel sets with Lebesgue-negligible boundary. 
Our concept of number and mean rigidity is slightly weaker than that in \cite{gp,bu}. 
We pose rigidity only for $ \SR $ and $ \TR $. This reduction is enough for our purpose. 

Mean rigidity is a critical property for sub-diffusivity. 
The Ginibre random point field is number-rigid but not mean-rigid. 
Hence, we introduce the algorithm to transform a function on $ \sSS $ with 
a number-rigid random point field to a function that is measurable concerning a mean-rigid $ \sigma $-field. 

For $ \rR \in \mathbb{N} $, let $ \GR $ and $ \mathcal{H}_{\rr } $ 
be the sub $ \sigma $-fields of $ \mathcal{B}(\SSoneD ) $ given by 
\begin{align}\label{:J1u}&
 \GR = \sigma [\NR , \MMR , \piRc ] 
,\quad \mathcal{H}_{\rr } = \sigma [\NTr , \MTr , \pi _{\TR }^c ] 
.\end{align}
We set the mean-rigid $ \sigma $-field $ \Gi $ as follows. 
\begin{align}\label{:J1w}&
\Gi = \bigcap _{\rR = 1}^{\infty} \GR 
.\end{align}
\begin{lemma} \label{l:J1}\thetag{1} 
$ \GR \supset \mathcal{G}_{\rR +1} $ and $ \mathcal{H}_{\rr } \supset \GR $ for each $ \rR \in \mathbb{N} $. 
\\
\thetag{2} 
$ \vartheta_{\aaa } (\Gi ) = \Gi $ for all $ \aaa \in \Rd $, where 
$ \vartheta_{\aaa } $ is the translation on $ \sSS $ in \eqref{:12t}. 
\end{lemma}
\PF
Recall that $ \TRR = \SRR \backslash \SR $. 
Then we see $ \pi _{\rR }^c = \pi _{\TRR } + \pi _{\rr + 1 }^c $. 
Hence, 
\begin{align} \notag 
 \sigma [\NR , \MMR , \piRc ] & = \sigma [\NR , \MMR , \pi _{\TRR } , \pi _{\rR +1}^c ] 
\supset 
 \sigma [ \mathrm{N}_{\rR + 1} , \mathrm{M}_{\rR + 1} , \pi _{\rR +1 }^c ] 
.\end{align}
This together with \eqref{:J1u} yields $ \GR \supset \mathcal{G}_{\rR +1} $. 

From $ \TR = \SR \backslash \sS _{\rR -1 }$, we find 
$ \pi _{\TR }^c = \pi _{\rr -1 } + \piRc $. 
Hence, we see 
\begin{align}& \notag 
 \sigma [\NTr , \MTr , \pi _{\TR }^c ] =
\sigma [\NTr , \MTr , \pi _{\rr - 1 } , \piRc ] 
\supset 
 \sigma [ \NR , \MMR , \piRc ] 
.\end{align}
This together with \eqref{:J1u} yields $ \mathcal{H}_{\rr } \supset \GR $. 
We thus obtain \thetag{1}. 

Without loss of generality, we set $ \lvert \aaa \rvert < 1 $. Let $ \SR (\aaa ) =\vartheta _{\aaa } (\SR ) $. Then, 
\begin{align}\label{:J1e}&
 \sS _{\rR - \lvert \aaa \rvert } \subset 
 \SR (\aaa ) \subset 
 \sS _{\rR + \lvert \aaa \rvert } 
.\end{align}
Replacing $ \SR $ with $ \SR (\aaa ) $ in \eqref{:J1m}, we define $ \NRa $ and $ \MRa $. 
Let $ \pi_{ \rR , \aaa }^c = \pi _{ \SR (\aaa ) ^c}$. Then from \eqref{:J1e}, we find that 
\begin{align}\label{:J1f}&
\mathcal{G}_{\rR - \lvert \aaa \rvert } 
\supset 
\sigma [\NRa , \MRa , \pi_{ \rR , \aaa }^c ] \supset 
\mathcal{G}_{\rR + \lvert \aaa \rvert } 
.\end{align}
From \eqref{:J1u}, we find that 
$ \vartheta_{\aaa } (\GR ) = \sigma [\NRa , \MRa , \pi_{ \rR , \aaa }^c ] $. 
From this and \eqref{:J1f}, 
\begin{align}\label{:J1g}&
\bigcap_{\rR =1}^{\infty}\mathcal{G}_{\rR - \lvert \aaa \rvert } 
\supset 
\bigcap_{\rR =1}^{\infty} 
 \vartheta_{\aaa } (\GR ) 
 \supset 
\bigcap_{\rR =1}^{\infty}
\mathcal{G}_{\rR + \lvert \aaa \rvert } 
.\end{align}
Using \eqref{:J1u} and \eqref{:J1g} and noting $ \bigcap_{\rR =1}^{\infty} \vartheta_{\aaa } (\GR ) = 
 \vartheta_{\aaa } ( \bigcap_{\rR =1}^{\infty} \GR ) $, 
we obtain 
$ \Gi \supset \vartheta_{\aaa } (\Gi ) \supset \Gi $, which yields \thetag{2}. 
\PFEND

Let $ \mu _{\aaa } $ be the reduced Palm measure of $ \mu $ conditioned at $ \aaa \in \Rd $ as before. 
Let $ \mu _{\aaa } (\cdot \vert \GR )$, $ \rR \in \Ni $, be the regular conditional probabilities. 
\begin{lemma} \label{l:J2}
$ \mucGi = \mu _{\aaa } (\, \cdot\, \vert \Gi ) \circ \vartheta_{\aaa }^{-1}
$ for each $ \aaa \in \Rd $. 
\end{lemma}
\PF
From the martingale convergence theorem, \eqref{:J1w}, and \lref{l:J1}\thetag{1}, 
\begin{align} &\notag 
 \muAGi (\sss ) = \limi{\rR } \muAGR (\sss ) 
\quad \text{ in $ L^1(\muz )$ and $ \muz $-a.s.\,$ \sss $ for any $ A \in \mathcal{B}(\sSS ) $} 
.\end{align}
Hence, the translation invariance of $ \mu $ and \lref{l:J1}\thetag{2} yield the claim. 
\PFEND

\begin{lemma} \label{l:J3} 
Assume that $ f $ is $ \Gi $-measurable. Then $ f \in \dbP $ and $ \Dp [f,f] = 0 $. 
\end{lemma}
\PF We write $ \sss = \sum_i \delta_{s_i}$ as before. 
Let $ \Upsone , \Upstwo , \Upsthree $ be as in \ssref{s:4}. 

Suppose $ \upsilon \in \Upsone $. 
Then for $ \sss \in \TT _{\rR }^m $, 
\begin{align} &\notag 
\upsilon = \PD{}{s_i} - \frac{1}{ m }\PD{}{\Gamma (\TR )} ,\quad s_i \in \TR 
.\end{align}
We note 
$ \sss (\TR ) = \mathrm{N}_{\TR } (\sss ) $ and 
$ \sum_{s_i \in \TR }s_i= \mathrm{M}_{\TR } (\sss )$. 
Note that $ f $ is $ \Gi $-measurable by assumption. 
Then, $ f $ is measurable with respect to 
$ \mathcal{H}_{\rR } = \sigma [\NTr , \MTr , \pi _{\TR }^c ] $ from \eqref{:J1u} and \lref{l:J1}\thetag{1}. 
Hence, we see that $ f $ is a function of $ \NTr (\sss ) $, $ \MTr (\sss ) $, and $ \pi _{\TR }^c (\sss )$. 
Note that $ \NTr (\sss ) = m $ on $ \TT _{\rR }^m $. 
Let $ f _{\TR , \sss }^m $ be the function defined after \eqref{:14e} for $ f $ and $ A = \TR $. 
Then $ f _{\TR , \sss }^m = f_{\TR , \pi _{\TR }^c (\sss )}^m $ from \eqref{:14c}. 
Let $ \tilde{f} _{\TR , \sss }^m $ be the function defined on $   \{ x \in \Rd ; x = x_1+\cdots + x_m , (x_i)_{i=1}^m \in (\TR )^m \} $ such that 
\begin{align}&\label{:J3b} 
 f _{\TR , \sss }^m (s_1,\ldots,s_m) = \tilde{f} _{\TR , \sss }^m (s_1+\cdots + s_m)
.\end{align}
The vector $ \upsilon $ is perpendicular to $ {\partial}/{\partial \Gamma (\TR )} $ 
as we see in \eqref{:31u}. Hence from \eqref{:J3b}, %
\begin{align}\label{:J3c}&
 \upsilon f _{\TR , \sss }^m (s_1,\ldots,s_m) = 
 \upsilon \tilde{f} _{\TR , \sss }^m (s_1+\cdots + s_m)= 0 
.\end{align}

Suppose $ \upsilon \in \Upstwo $. Then 
$ \sss (\sSS _{\rrr - 1 }) = m $ and $ \sss (\SR ) = m + n $. We find 
\begin{align}& \label{:J3h}
\upsilon = 
 \frac{1}{\sqrt{ m }} 
 \PD{}{\Gamma ( \Upp )} - \frac{\sqrt{m}}{ m + n } \PD{}{\Gamma ( \Up ) }
.\end{align}
From \eqref{:J1w}, $ f $ is $ \mathcal{G}_{\rR } $-measurable. 
Hence, $ f $ is a function of 
$ \mathrm{N}_{\rR } (\sss ) $, $ \mathrm{M}_{\rR } (\sss ) $, and $ \piRc (\sss )$. 
Note that $ \NR (\sss ) = m+n$ and $ \mathrm{M}_{\rR } (\sss ) = s_1 + \cdots + s_{m+n} $ on 
$ \sSS _{\rrr - 1 }^m \cap \TT _{\rrr }^n $. 
Then we have a function $ \tilde{f} _{\rR , \sss }^{m+n} $ on 
$\{ x \in \Rd ; x = x_1+\cdots + x_{m+n} , (x_i)_{i=1}^{m+n} \in \SR ^{m+n} \}$ such that 
\begin{align}&\label{:J3k}
 f _{\rR , \sss }^{m+n} (s_1,\ldots,s_{m+n}) = 
 \tilde{f} _{\rR , \sss }^{m+n} (s_1+\cdots + s_{m+n})
.\end{align}
From \eqref{:31v} and \eqref{:J3h}, $ \upsilon $ is perpendicular to $ \partial / \partial \Gamma (\SR )$. Hence from  this and \eqref{:J3k}, 
\begin{align}\label{:J3n}&
 \upsilon f _{\rR , \sss }^{m+n} (s_1,\ldots,s_{m+n}) = 0 
.\end{align}


 Suppose $ \upsilon \in \Upsthree $. Then  $ \sss (\TR ) = m $ and $ \sss (\SR ) = n $. Hence, we see 
 \begin{align}& \label{:J3o}
 \upsilon = 
 \frac{1}{\sqrt{ m }} 
 \PD{}{\Gamma ( \TR )} - \frac{\sqrt{m}}{ n } \PD{}{\Gamma ( \Up ) }
 .\end{align}
 From \eqref{:J1w}, $ f $ is $ \mathcal{G}_{\rR } $-measurable. 
 Hence, $ f $ is a function of 
 $ \mathrm{N}_{\rR } (\sss ) $, $ \mathrm{M}_{\rR } (\sss ) $, and $ \piRc (\sss )$. 
 Note that $ \NR (\sss ) = n$ and $ \mathrm{M}_{\rR } (\sss ) = s_1 + \cdots + s_{n} $ on $ \SSR ^n $. 
 Thus, we find a function $ \tilde{f} _{\rR , \sss }^{n} $ on 
 $\{ x \in \Rd ; x = x_1+\cdots + x_{n} , (x_i)_{i=1}^{n} \in \SR ^n \}$ such that 
 \begin{align}&\label{:J3l}
 f _{\rR , \sss }^{n} (s_1,\ldots,s_{n}) = 
 \tilde{f} _{\rR , \sss }^{n} (s_1+\cdots + s_{n})
 .\end{align}
 From \eqref{:31v} and \eqref{:J3o}, $ \upsilon $ is perpendicular to 
 $ \partial / \partial \Gamma (\SR )$. Hence from this and \eqref{:J3l}, 
 \begin{align}\label{:J3p}&
 \upsilon f _{\rR , \sss }^{n} (s_1,\ldots,s_{n}) = 0 
 .\end{align}

Putting \eqref{:J3c}, \eqref{:J3n}, and \eqref{:J3p} together and recalling \eqref{:4p}, we obtain 
\begin{align*}&
\DDD _{\upsilon} [f,f] = 0 \quad \text{ for all } \upsilon \in \Upsilon 
.\end{align*}
Hence using \eqref{:4r}, we conclude $ f \in \dbP $ and 
$ \Dp [f,f] = \sum_{ \upsilon \in \Upsilon } \DDD _{\upsilon} [f,f] = 0 $. 
\PFEND

\section{Proof of the main theorems (Theorems \ref{l:11}--\ref{l:15})}\label{s:E}
In \sref{s:E}, we complete the proof of the main theorems. 
In \ssref{s:Q}, we prove Theorems \ref{l:13}--\ref{l:15}. 
Subsections \ref{s:K}--\ref{s:6V} are devoted to the preparation of the proof of \tref{l:11}. 
In \ssref{s:K}, we present the solution $ \XRp $ to the stochastic differential equation associated with the Dirichlet form $ ( \EP _{\rrrr }, \dQbP ) $ on $ \Lm $. 
In Subsections \ref{s:K3} and \ref{s:K4}, we prepare a sequence of lemmas concerning the convergence of $ \{ \XRp \}_{\rR } $. Using these lemmas, we prove the limit points of the sequence $\{ \XRp \}_{\rR } $ are solutions of \eqref{:10a} in \pref{l:6V} in \ssref{s:6V}. 
With the uniqueness of weak solutions of \eqref{:10a}, we prove that the Ginibre random point field satisfies \As{A4} in \pref{l:K8} and complete the proof of \tref{l:11} in \ssref{s:65}.

\subsection{Proof of Theorems \ref{l:13}--\ref{l:15}} \label{s:Q}
In \ssref{s:Q}, we prove Theorems \ref{l:13}--\ref{l:15}. 

Assume \As{A1}--\As{A3} and \As{A5}. 
Thus, $ \mu $ is translation invariant and irreducibly one-decomposable with $ \SSDDone $. 
Let $ \muz $ be the reduced Palm measure conditioned at the origin. 
Let $ \mutoneD $ be the random point field introduced in \eqref{:X2w} for $ m = 1 $. 
Because of \eqref{:X2y}, $ \mutoneD \approx \muz $. 

For $ \mutoneD $-a.s.\,$ \yy \in \SSoneD $, let $ \muyy $ be the dual reduced Palm measure defined in \lref{l:X4} for $ m = 1 $. 
Using \eqref{:X4x} and $ \mutoneD \approx \muz $, we have 
\begin{align}\label{:Q1o}&
\mu ( \SSone \Vert \yy ) = 1 \quad \text{ for $ \muz $-a.s.\,$ \yy $}
.\end{align}
From \eqref{:Q1o} and $ \SSone = \{\delta_x; x \in \mathbb{R}^d\} $, 
 set the probability measure $ \muyC $ on $ \Rd $ by 
\begin{align}& \label{:Q1p}& \quad 
\muiyC ( \cdot ) = \mu ( \{ \delta_x ; x \in \cdot \} \Vert \yy ) \quad 
\quad \text{ for $ \muz $-a.s.\,$ \yy $}
.\end{align}

For $ L \in \mathbb{N} $, let $ \map{\xiL }{\mathbb{R} }{\mathbb{R} }$ be a smooth non-decreasing function such that 
\begin{align}\label{:Q1q}&
 0 \le \xiL ' (t) \le 2 
,\quad 
 \xiL (t) = 
 \begin{cases}
 L+1 & L + 2 \le t\\
t & \lvert t\rvert \le L \\
 -L - 1 & t \le -L-2
. \end{cases}
\end{align}
Let $ \Gi $ be as in \eqref{:J1w}. 
Let $ x = (\xP ) \Pdd \in \Rd $. 
We set $ \chiL = (\chiLp ) \Pdd $ by 
\begin{align}\label{:Q1r}&
\chiLp ( \sss ) =
\intSS 
\Big\{
 \int_{\Rd }\xiL ( \xP ) \muiyC (dx ) 
\Big\} 
\muz (d\yy \vert \Gi )(\sss )
.\end{align}
Here 
$ \muz ( \cdot \, \vert \Gi )(\sss ) $ is the regular conditional probability of $ \muz $ concerning $ \Gi $. 
We note that $ \chiLp $ is $ \Gi $-measurable by construction and that 
$ \chiLp $ is neither a continuous nor local function on $ \sSS $. 
Let $ \DsftQ $ be as in \eqref{:14m}. 
Let 
\begin{align} &\notag 
 \domtrn = \bigcap _{ q = 1 }^d \Big\{ f \in \Lmz ; \DsftQ f (\sss ) 
 \text{ exists for $ \muz $-a.s.\,$ \sss \in \SSoneD $} %
\Big\} 
.\end{align}

\begin{lemma}\label{l:Q1} \thetag{1} $ \chiLp \in \domtrn $ and 
$ \lvert \DsftQ \chiLp (\sss ) \rvert \le 2 $ for $ p , q = 1,\ldots , d $. 
\\\thetag{2} 
$ \{ \chiLp \} $ is an $ \EYone $-Cauchy sequence as $ L \to \infty $ satisfying 
for $ p , q = 1,\ldots , d $ 
\begin{align} &\notag 
\limi{ L } \DsftQ \chiLp (\sss ) = \delta_{ p , q } 
\quad \text{ for $ \muz $-a.s.\,$ \sss $ and in $ \Lmz $}
.\end{align}
\end{lemma}
\PF 
From \eqref{:Q1p} and \lref{l:X5}, we have for $ \muz $-a.s.\,$ \yy $ 
\begin{align}\label{:Q1e}&
\muiyC (dx )= \muiyCC (d(x- \aaa ))
.\end{align}
 From \lref{l:J2}, we deduce
\begin{align}\label{:Q1f}&
 \muz (d\yy \vert \Gi )(\sss ) = 
\mu _{0 - \aaa } (\vartheta _{\aaa } (\cdot ) \in d\yy \vert \Gi )(\vartheta _{\aaa }(\sss ))
.\end{align}
Hence using \eqref{:Q1r}, \eqref{:Q1e}, and \eqref{:Q1f}, we obtain for $ \muz $-a.s.\,$ \sss $ 
\begin{align} \notag &
\frac{1}{h} \{
\chiLp ( \vartheta _{ h \eQ }(\sss )) -\chiLp ( \sss ) 
\}
\\\notag =& 
\frac{1}{h} \intSS 
\Big\{
 \int_{\Rd } \xiLP \muiyCh (dx) - \int_{\Rd } \xiLP \muiyC (dx) 
\Big\}
\muz (d\yy \vert \Gi )(\sss )
\\ \notag = &
\frac{1}{h} \intSS 
\Big\{
 \int_{\Rd } \xiLh \muiyC (dx) - \int_{\Rd } \xiLP \muiyC (dx) 
\Big\}
\muz (d\yy \vert \Gi )(\sss )
\\ \label{:Q1g} = & \intSS 
 \int_{\Rd }
\frac{1}{h}
\Big\{ \xiLh - \xiLP 
\Big\}\muiyC (dx) 
\muz (d\yy \vert \Gi )(\sss )
.\end{align}
Using \eqref{:Q1q}, \eqref{:Q1r}, and \eqref{:Q1g}, we find $ \chiLp \in \domtrn $ 
and that for $ \muz $-a.s.\,$ \sss $ 
\begin{align} \label{:Q1h}&
 \DsftQ \chiLp (\sss ) = 
\deltaPQ \intSS \int_{\Rd } \xiL ' ( \xP ) \muiyC (dx) \muz (d\yy \vert \Gi )(\sss )
.\end{align}
Hence, $ \lvert \DsftQ \chiLp (\sss ) \rvert \le 2 $ from \eqref{:Q1q} and \eqref{:Q1h}. 
We have thus obtained \thetag{1}.

From \eqref{:Q1q}, \eqref{:Q1h}, and the Lebesgue convergence theorem, 
we obtain \thetag{2}. 
\PFEND

\noindent {\em Proof of \tref{l:15}}. 
Let $ \dbP $ be as in \eqref{:42q}. By \eqref{:Q1r}, $ \chiLp $ is $ \Gi $-measurable. 
Hence from \lref{l:J3}, we have 
\begin{align}\label{:Q0a}& 
 \chiLp \in \dbP , \quad \Dp [\chiLp , \chiLp ] = 0 
.\end{align}
From \lref{l:Q1}\thetag{1}, we see 
\begin{align}\label{:Q0A}&
\chiLp \in \domtrn , \quad \vert \DsftQ \chiLp \vert \le 2 
.\end{align}
From \eqref{:Q0a} and \eqref{:Q0A}, we obtain $ \chiLp \in \dbYp $, where $ \dbYp $ is as in \eqref{:71z}. 

By \lref{l:71}\thetag{2}, $ (\E ^{\YY , \perp } , \dom ^{\YY , \perp } ) = (\EY , \dYL )$. 
Thus, $ \dYL $ is the closure of $ \dbYp $ with respect to $ \E ^{\YY , \perp } $. 
This implies $ \chiLp \in \dbYp \subset \dYL$. 
%
Hence for each $ \epsilon > 0 $ and $ L \in \mathbb{N} $, we find $ \chiLpE\in \dbY $ such that 
\begin{align} \notag & 
\int_{\sSS } \half \sumQ \Big\lvert \DsftQ \chiLp - \DsftQ \chiLpE \Big\rvert ^2 
\, d\muz < \epsilon 
,\\\label{:Q0d} &
\EYtwo ( \chiLp - \chiLpE , \chiLp - \chiLpE ) < \epsilon 
.\end{align}
By \lref{l:Q1}\thetag{2} and \eqref{:Q0a}, 
 $ \{ \chiLp \} $ is an $ \EY $-Cauchy sequence such that 
\begin{align} \notag &
\limi{L}\int_{\sSS } \half \sumQ \Big\lvert \DsftQ \chiLp - \deltaPQ \Big\rvert ^2 \, d\muz = 0 
,\\ &\label{:Q0c}
\EYtwo ( \chiLp , \chiLp ) = 0 \quad \text{ for each }L 
.\end{align}
Using $ \chiLpE \in \dbY $, \eqref{:Q0d}, and \eqref{:Q0c}, we see \As{A7}. 
Thus, we obtain \thetag{1}. 

From \tref{l:93}, we have 
$ \limz{\epsilon } \epsilon X _{t/ \epsilon ^2 } = 0 $
weakly in $ C([0,\infty); \Rd )$ in $ \muz $-measure. 
Combining this with Theorem 3.7 in \cite{o.nc}, we obtain \thetag{2}. 
\qed

\smallskip
\noindent {\em Proof of Theorems \ref{l:13} and \ref{l:14}}. 
We obtain \tref{l:13} from \tref{l:15}. 
 \tref{l:14} follows from \tref{l:13} and \lref{l:X6}. 
\qed 

\subsection{Preparation for proof of \tref{l:11}} \label{s:K} 
In \ssref{s:K}, we present the stochastic differential equation describing the dynamics given by $ ( \EP _{\rrrr }, \dQbP ) $ on $ \Lmg $.

Let $ ( \EP _{\rrrr }, \dQbP ) $ be the Dirichlet form defined in \lref{l:41} for $ \mu = \mug $. 
The carr\'{e} du champ of the Dirichlet form is then given by  $ \DDDR ^{\perp }$. From \eqref{:32a},  
\begin{align} \label{:K2a}&
\Dp _{\rrrr } [f,f] (\sss ) = \DDD _{\rrrr } [ f , f ] (\sss ) -
\mathbb{U}_{\rrrr }^{\gamma }[f,f] (\sss ) 
.\end{align}
Let $ \lab = (\labi )_{i \in \mathbb{N} }$ be the label on $ \SSsi $ such that 
$ \lvert \labi (\sss ) \rvert \le \lvert \labii (\sss ) \rvert $ for all $ i \in \mathbb{N} $ as defined before \tref{l:11}. 
We consider the stochastic differential equation describing the labeled dynamics $ \XRp $ given by 
 $ ( \EP _{\rrrr }, \dQbP ) $ on $ \Lmg $. We take $ \XRp _0 = \lab (\sss )$. 

Let $ \bbb = (1/2) \dgin $ be the drift coefficient of ISDE \eqref{:10a} given by 
\begin{align}\label{:K1x}&
\bbb ( x , \sss ) 
= \limi{ \rR } \sum_{\lvert x -s_i \rvert < \rR } \frac{ x -s_i }{\lvert x -s_i \rvert ^2 }
\quad \text{ in } \Llocmugonep , \ 1 \le p < 2 
.\end{align}
Formula \eqref{:K1x} follows from Theorem 61 in \cite{o.isde}. Recall that $ \mug $ is tail trivial. 
Hence, $\dgin $ satisfies \eqref{:13i} not only for $ \varphi \in \dibone $ but also for $ \varphi \in \di \ot \dbb $. 
Then for any $\varphi \in C_{0}^{\infty} (\Rtwo )\ot \dbb $ satisfying the Neumann boundary condition on $ \partial \SR \ts \sSS $, it holds that 
\begin{align}&\notag 
\int _{\SR \times \mathsf{S} } \dmu (x,\sss )\varphi (x,\sss ) \muone _{\mathrm{Gin}, \rR , \yy } (dx d\sss ) =
 - \int _{\SR \times \mathsf{S} } \nabla_x \varphi (x,\sss ) \muone _{\mathrm{Gin}, \rR , \yy } (dx d\sss ) 
,\end{align}
where $ \muone _{\mathrm{Gin}, \rR , \yy } $ is the one-Campbell measure of 
$ \mug (\cdot \vert \piRc (\sss ) =\piRc (\yy ))$. 
From this and \eqref{:K2a}, the stochastic differential equation of $ \XRp $  is given by 
\begin{align}\notag 
d X_t^{\Rpi } =& dB_t^i + 
 \bbb 
( \XRpDti ) dt + 
 \nR (X_t^{\Rpi } ) d L_t^{\Rpi }
\\ \notag &
- \frac{1}{\m }\sum_{k=1}^\m d B_t^k 
- \frac{1}{\m } \sum_{k=1}^\m \bbb 
( \XRpDt ) dt 
\\\notag &
 - 
\frac{1}{\m } \sum_{k=1}^\m 
 \nR (X_t^{\rR ,\perp , k } ) d L_t^{\rR ,\perp , k } ,\ \quad 
1 \le i \le \m 
,\\ \notag 
d L_t^{\Rpi } = &1_{\partial \SR } (X_t^{\Rpi } ) d L_t^{\Rpi } 
,\ \quad 1 \le i \le \m 
,\\ \label{:K2c} 
X_t^{\Rpi } = & X_0^{\Rpi } , \ \quad \m < i < \infty 
.\end{align}
Here, $ \mathsf{X}_t^{ \RpiD } = \sum_{j \ne i }^{\infty} \delta _{X_t^{ R , \perp ,j}}$, $ L_t^{\Rpi }$ are non-negative increasing processes, and $ \nR $ is the inner normal vector on $ \partial \SR $. 
Then, $ \XRp$ is frozen outside $ \SRover $. 

Let $ \vp \mug (\cdot) = \int_{\cdot } \vp (\xx ) \mug (d\xx )$, where $ \vp \in \Lmug $, 
 $ \vp \ge 0 $, $ \int_{\sSS } \vp d\mug = 1$ and 
$ \Ct = \sup_{\mathsf{x \in \sSS }}\vert \vp (\mathsf{x}) \vert < \infty \label{;63} $. 
We assume 
\begin{align}\label{:K2d}&
 \XRp _0 \elaw (\vp \mug ) \circ \lab ^{-1}
.\end{align}
We prove that $ \{\XRp \}$ converges to the solution of \eqref{:10a} in \pref{l:6V}.

\subsection{The main terms of perpendicular SDEs} \label{s:K3}

The goal of \ssref{s:K3} is to prove the convergence of the main terms of coefficients of the stochastic differential equations \eqref{:K2c} (see \pref{l:4X}). 

Let $ \XRp = (X ^{\Rpi })_{i=1}^{\infty} $, $ \X _0^{\Rp } = \lab (\sss ) $, be the solution of \eqref{:K2c} with initial distribution \eqref{:K2d}. From Lyons--Zheng's decomposition \cite[p.284]{c-f}, $ X ^{\Rpi } $ can be written as the sum of the martingale additive functionals of the diffusion associated with the Dirichlet form $ ( \EP _{\rrrr }, \dQbP ) $ on $ \Lmg $.  For $ 1\le i \le \sss (\SRover ) $, 
\begin{align}\label{:K3x}& 
 X_t^{\Rpi } - X_0^{\Rpi }
= \half \Big\{ \mathit{M}_t^{[x_i\ot 1 ]} - \mathit{M}_t^{[x_i\ot 1 ]}\circ r_t \Big\} 
,\end{align}
where $ r_t $ is the time-reversal operator on the path space on $ [0,\infty)$ such that 
$ r_t (\omega)(s) = \omega (t-s)$ if $ 0\le s \le t $ and  $ r_t (\omega)(s) = \omega (0)$ if $ t \le s  $.  We see 
\begin{align}\label{:K3y} 
 \mathit{M}_t^{[x_i\ot 1 ]} & =B_t^i -\frac{1}{\sss (\SRover ) }\sum_{k=1}^{\sss (\SRover ) } B_t^k 
.\end{align}
We write $ x_i = (x_{i,p})_{p=1,2}$ and $ \mathit{M}^{[x_{i}\ot 1 ]}  = (\mathit{M}^{[x_{i,p}\ot 1 ]})_{p=1,2}$. 
Because $ (B^i)_{i\in\mathbb{N} }$ is the standard Brownian motion in $ (\mathbb{R} ^2)^{\mathbb{N} }$ 
we have 
\begin{align}\label{:K3Y}&
\langle \mathit{M}^{[x_{i,p}\ot 1 ]} ,\mathit{M}^{[x_{i,q}\ot 1 ]} \rangle _t = 
\Big( 1 - \frac{1}{\sss (\SRover ) }\Big) \delta_{p,q} t 
.\end{align}
We note that $ X_t^{\Rpi } = X_0^{\Rpi } = \lab ^i(\sss )$, $ 0 \le t < \infty $, for each $ i > \m $.

\begin{lemma} \label{l:KA} 
There exists a positive constant $ \Ct \label{;76} > 0 $ satisfying the following. 

\begin{align}\label{:K5s} & 
\sup_{ \rR , i \in \mathbb{N} } E [ \lvert X_t^{\Rpi } -X_u^{\Rpi } \rvert ^4 ] 
 \le \cref{;76} \lvert t - u \rvert ^ 2 , \text{ } 0 \le t,u < \infty 
,\\ \label{:K5c}&
\limi{a} \liminfi{ R }
P ( \max_{1\le i \le n } \maxT \lvert X_t^{\Rpi } \rvert \le a ) = 1 
\quad \text{ for each } n , T \in \mathbb{N} 
,\\\label{:K5d} &
\limi{l} \inf_{ R \in \mathbb{N} }
P ( \IrT ( \XRp ) \le l ) = 1 \quad \text{ for each } r , T \in \mathbb{N} 
.\end{align}
Here $ \IrT $ is as in \eqref{:61w}. 
\end{lemma}
\begin{proof} 
From \eqref{:K3x}--\eqref{:K3Y} and the martingale inequality, we obtain \lref{l:KA}. 
(see Lemma 4.2 in \cite{k-o-t.udf} for detail). 
\end{proof} 

In general, a family of probability measures $ \{ m_a \} $ in a Polish space is compact under the topology of weak convergence if and only if for any $ \epsilon > 0 $ there exists a compact set $ K $ such that $ \inf_a m_a (K) \ge 1- \epsilon $. We call $ \{ m_a \} $ tight in this case. 
The tightness of probability measures on a countable product of Polish spaces under the product topology follows from the tightness of the distributions of each component. 
We say a family of random variables is tight if the family of their distributions is tight. 

We equip $\Cinfty $ with the metric given by \eqref{:Z3x}. 
Hence, the tightness in $ \CT $ for all $ T \in \mathbb{N} $ is equivalent to the tightness in $\Cinfty $. 
Let $ \upath $ be the map from $\Cinfty $ to the space of measure valued paths such that 
\begin{align}& \notag 
 \upath (\mathbf{w})_t :=\ulab (\mathbf{w}(t))= \sumii{i} \delta_{w^i(t)} 
\quad \text{ for $ \mathbf{w}=(w^i)_{i=1}^{\infty}$}
.\end{align}
Note that $ \upath (\Cinfty ) \not\subset \CiSS $. 
We endow $ \CiSS $ with the metric $ \rho_{\mathrm{upath}}$ defined by \eqref{:Z3z}. 
Then $ \CiSS $ becomes a complete separable metric space. For $ \XRp = (X ^{\Rpi })_{i=1}^{\infty}$, 
we set 
\begin{align} &\notag 
\XXRp = \upath (\XRp ) , \quad \XXRpDi = \sum_{j\ne i }^{\infty} \delta_{X ^{\Rpj }}
.\end{align}

\begin{lemma}\label{l:KB}
\thetag{1} 
$ \{ \XRp \}_{\rR \in \mathbb{N}} $ is tight in $ \Cinfty $. 
\\\thetag{2} 
$ \{ \XXRp \}_{\rR \in \mathbb{N}} $ and $ \{ \XXRpDi \}_{\rR \in \mathbb{N}} $ are tight in $ \CiSS $ for each $ i \in \mathbb{N}$. 
\end{lemma}
\PF
Using \eqref{:K2d} and \eqref{:K5s}, we see that $ \{ \XRp \}_{\rR \in \mathbb{N}} $ is tight in $ \CRtwoN $. 
We have thus obtained \thetag{1}. Let 
\begin{align*}&
\mathcal{NBJ}=\{ \mathbf{w}\in \CiRdN ; 
\text{$ \IRT (\mathbf{w}) < \infty $ for each $ R , T \in \mathbb{N}$} \}
.\end{align*}
From \lref{l:Z3}, $ \map{\upath }{\mathcal{NBJ}}{\CiSS }$ is continuous. 
From \eqref{:K5d}, $ \XRp $ takes value in $ \mathcal{NBJ}$. 
Hence using \thetag{1} and the continuity of $ \upath $ on $ \mathcal{NBJ}$, we see 
$ \{ \XXRp \}_{\rR \in \mathbb{N}} $ is tight in $ \CiSS $. 
Note that 
\begin{align*}&
 \XXRpDi _t = \XXRp _t - \delta_{X ^{\Rpi }_t }
.\end{align*}
Hence, the tightness of 
$ \{ \XXRpDi \}_{\rR \in \mathbb{N}} $ follows from the tightness of $ \{ \XXRp \}_{\rR \in \mathbb{N}} $ and 
$ \{ \XRp \}_{\rR \in \mathbb{N}} $. This completes the proof of \thetag{2}. 
\PFEND

Using \lref{l:41} and \lref{l:KB}, we characterize the limit points of $ \XRp $ in terms of the Dirichlet form $ (\EP , \dP ) $ on $ \Lmg $. 
\begin{lemma} \label{l:6V1} 
Let $ \Xp $ be an arbitrary limit point of $ \{ \XRp \}_{R\in\mathbb{N} }$ under the convergence in law in $\Cinfty $. Let $ \XXp = \upath (\Xp )$. Then the following hold. 
\\
\thetag{1} 
$ \XXp $ is associated with the Dirichlet form $ (\EP , \dP ) $ on $ \Lmg $ and 
\begin{align}&\label{:K4y}
\limi{\rR } \XXRp = \XXp \quad \text{ weakly in $ \CiSS $} 
.\end{align}
\thetag{2} $ \Xp $ is unique in law in $\Cinfty $ and 
\begin{align} \label{:K4x}&
\limi{\rR } \XRp = \X ^{\perp } \quad \text{ weakly in $ \CRtwoN $} 
.\end{align}
\end{lemma}
\begin{proof}
From \lref{l:41}, $ (\EP , \dP ) $ is the increasing limit of $ ( \EP _{\rrrr }, \dQbP ) $ on $ \Lmg $. 
Hence, the distribution of $ \XXRp $ converges in finite-dimensional distributions to that of $ \XXp $ associated with $ (\EP , \dP ) $ on $ \Lmg $. From this and \lref{l:KB}\thetag{2}, we obtain \eqref{:K4y}, which yields \thetag{1}. 

Note that $ \Xp = \lpath (\XXp )$. From this and \thetag{1}, $ \Xp $ is unique in law. 
From \lref{l:KB}\thetag{1}, $ \{ \XRp \}_{R\in\mathbb{N} }$ is tight in $\Cinfty $. 
Hence, \eqref{:K4x} holds. 
\end{proof}

We use the cut-off functions $ \chi_{\rsp }= \chi_{(\mathsf{p},\mathsf{q},\mathsf{r})} $ defined by \thetag{11.14} in \cite[p.1228]{o-t.tail} with $ m = 1 $. 
To simplify the notation, we set $ \llL = \limi{\mathsf{r}} \limi{\mathsf{s}}\limi{\mathsf{p}} $. 
The parametrization $ \rsp =(\mathsf{r},\mathsf{s},\mathsf{p}) $ comes from the construction of $ \chi_{\rsp }$. 
\begin{lemma}[{\cite{o-t.tail}}] \label{l:K4}
Let $ \bbb $ be as in \eqref{:K1x}. There exist $ \bbb _{\rsp } \in C_b (\RtwoSS ) $ such that 
\begin{align} \label{:K4c}&
\llL \| \Qbbrsp \|_{\Lmugonep } = 0 ,\ 1 \le p < 2 , \qQ \in \mathbb{N}
.\end{align}
\end{lemma}
\PF 
By construction, $ 0 \le \chi_{\rsp } \le 1 $, $ \llL \chi_{\rsp } = 1 $ on the support of $ \chi_{\rsp }$, and 
$ \chi_{\rsp } \in C_b (\RtwoSS ) $. Let $ \bbb _{\rsp }= (1/2) \chi _{\rsp } \dgin $. 
Then we obtain \eqref{:K4c} using the Lebesgue convergence theorem. 
 $ \bbb _{\rsp } \in C_b (\RtwoSS ) $ follows from the construction of $ \chi_{\rsp }$. 
\PFEND

\begin{lemma} \label{l:44}
Let $ \bRp $ and $ \bp\in \CRtwoN $ be such that 
\begin{align}&\label{:44x} 
\bRp = \pL \int_0^{\cdot} \bbb \XXXRiu du \pR 
, 
\bp = \pL \int_0^{\cdot} \bbb \XXXsiu du \pR 
.\end{align}
Then we have, weakly in $ \CRtwoN \ts \CRtwoN $, 
\begin{align}& \label{:44d} 
\limi{\rR }\pL \XRp , \mathbf{b}^{\Rp } \Big) 
= \pL \Xp , \bp \Big) 
.\end{align}
\end{lemma}
\PF 
Using \lref{l:50} for $ \mug $, we see $ ( \ER ^{[1],\perp } , \db ^{[1],\perp } )$ is closable on $ \Lmugone $. 
Let $ ( \ER ^{[1],\perp } , \dR ^{[1],\perp } )$ be its closure on $ \Lmugone $. 
Let $ \XXXRone $ be the diffusion associated with $ ( \ER ^{[1],\perp } , \dR ^{[1],\perp } )$ on $ \Lmugone $. 
Recall that $ \XRp = (X ^{\Rpi })_{i=1}^{\infty} $ is a solution to \eqref{:K2c}. 
Then each one-labeled process $ \XXXRi $ is equivalent in law to the diffusion $ \XXXRone $ 
with the same initial distribution as $ \XXXRi $. 

From \eqref{:K2d}, $ \XRp _0 \elaw (\vp \mug ) \circ \lab ^{-1}$. Hence, 
\begin{align}& \label{:44e}
\XXXRonez \elaw \vp \mug \circ (\lab ^i (\sss ), \sum_{j\ne i}^{\infty} \delta_{\lab ^j (\sss )})^{-1}
.\end{align}
From these, we have for all $ 0 \le u < \infty $ 
\begin{align} \notag
\EE \big[ 
\SQXu \big\vert  &( \bbrsp ) \XXXRiu \big\vert \big]
= 
\EE \big[ 
\SQXoneu \big\vert ( \bbrsp ) \XXXRoneu \big\vert \big]
\\ \notag& \le 
 \cref{;63}\int_{\SQ \ts \mathsf{S} } 
 \big\vert ( \bbrsp ) (x , \yy ) \big\vert 
 \mugone (dx d\mathsf{y}) \quad \text{ by \eqref{:44e}}
\\ \label{:44p} &\le  \cref{;63} \cref{;44b} \| \Qbbrsp \|_{\Lmugonep } 
.\end{align}
Here $ \cref{;63} = \sup_{\mathsf{x \in \sSS }}\vert \vp (\mathsf{x}) \vert $ 
and $ \Ct = \mugone (\SQ \ts \sSS )^{\frac{p-1}{p}} \label{;44b}$. 
From \eqref{:44p} and \eqref{:K4c}, 
\begin{align}& \notag 
\llL \limsupi{\rR } 
\EE \Big[\int_0^t \SQXu \big\vert ( \bbrsp ) \XXXRiu \big\vert du \Big] 
%
 \\ \label{:44b} 
\le & t \cref{;63} \cref{;44b} \llL \| \Qbbrsp \|_{\Lmugonep }
= 0 
.\end{align}
In the same fashion as \eqref{:44b}, we have from \eqref{:K4c} 
\begin{align}\label{:44c} & 
\llL \EE \Big[\int_0^t 1_{\Sr }(\Xsiu ) \big\vert ( \bbrsp ) \XXXsiu \big\vert du \Big]=0 
.\end{align}

We set 
\begin{align*}&
\brspRp  = \pL \int_0^{\cdot} \bbb _{\rsp } \XXXRiu du \pR , \quad 
\brsp  = \pL \int_0^{\cdot} \bbb _{\rsp }(X_u^{\si }, \mathsf{X}_u^{\iD }) du \pR 
.\end{align*}
Using \lref{l:6V1} and $\bbb _{\rsp } \in C_b (\RtwoSS ) $, we deduce 
\begin{align} \label{:44B}
\limi{\rR }&\pL \XRp , \brspRp \Big) = \pL \Xp , \brsp \Big) 
\end{align}
weakly in $ \Cinfty \ts \Cinfty $. From \eqref{:44b}--\eqref{:44B}, we have 
\begin{align}\notag 
 \limi{\rR } \pL \XRp , \bRp \Big) 
 & 
= \limi{\rR } \pL \XRp , \brspRp + ( \bRp - \brspRp ) \Big) 
\\ \notag &=
\llL \limi{\rR } \pL \XRp , \brspRp \Big) \quad \text{ by \eqref{:44b}}
\\ \notag &=
\llL \pL \Xp , \brsp \Big) \quad \text{ by \eqref{:44B} }
\\ & \notag 
= \pL \Xp , \bp \Big) 
\quad \text{ by \eqref{:44c}}
\end{align}
weakly in $ \Cinfty \ts \Cinfty $. This implies \lref{l:44}. 
\PFEND

\begin{lemma} \label{l:K6} 
For each $ i , T \in \mathbb{N} $ and $ \epsilon > 0 $, 
\begin{align}\label{:K6a}&
\limi{\rR } P \Big( \maxT \Big\lvert 
 \int_0^t
 \nR (X_u^{\rR ,\perp , i} ) d L_u^{\rR ,\perp , i}\Big\rvert \ge \epsilon \Big) = 0
.\end{align}
\end{lemma}

\PF 
Using \eqref{:K3x}, \eqref{:K3Y}, and the martingale inequality, we deduce for all $ h \in \mathbb{N} $ 
\begin{align} \notag &
 P \Big(
\max_{0\le t \le T} \Big\lvert X_t^{\rR ,\perp , i}- X_0^{\rR ,\perp , i} \Big\rvert \ge h 
\Big)
 \\ \notag 
\le\, &
 P \Big(
\max_{0\le t \le T} \Big\lvert \mathit{M}_t^{[x_i\ot 1 ]}  \Big\rvert \ge \frac{h}{2} 
\Big)+
 P \Big(
\max_{0\le t \le T} \Big\lvert \mathit{M}_t^{[x_i\ot 1 ]}\circ r_t  \Big\rvert \ge \frac{h}{2} 
\Big)
\\ 
=\,& 2 
 P \Big( 
\max_{0\le t \le T} \Big\lvert \mathit{M}_t^{[x_i\ot 1 ]}  \Big\rvert \ge \frac{h}{2} 
\Big)
\label{:K6f} 
\le \cref{;C5}\int_{\frac{h}{2} }^{\infty} e^{- t^2/\cref{;C5a}}dt 
.\end{align}
Here $ \Ct \label{;C5}$, $ \Ct \label{;C5a} > 0 $ are constants independent of $ \rR , i , h \in \mathbb{N} $. 
 $ L_u^{\rR ,\perp , i }$ increases only when $ X_u^{\rR ,\perp , i}$ is on the boundary $ \partial \SR $. 
Hence, \eqref{:K6a} follows from \eqref{:K6f}. 
\PFEND

\begin{proposition}\label{l:4X}
Weakly in $ (C([0,\infty); \Rtwo ) ^{\mathbb{N}} )^3 $, 
\begin{align} 
\limi{\rR } 
\Big(\XRp , \bRp 
 , (\int_0^{\cdot } \nR (X_u^{\rR ,\perp , i} ) d L_u^{\rR ,\perp , i} )_{i=1}^{\infty} 
 \Big) 
& \label{:46a}
= \Big(\Xp , \bp ,\mathbf{0} \Big) 
.\end{align}
\end{proposition}

\begin{proof}
\eqref{:46a} follows from \eqref{:44d} and \eqref{:K6a} immediately. 
\end{proof}

\subsection{Vanishing of the $ \gamma $-coefficient terms}\label{s:K4}

In \ssref{s:K3}, we proved the convergence of the main terms of stochastic differential equations \eqref{:K2c}. 
The goal of \ssref{s:K4} is to prove the remainder terms certainly vanish. 
To simplify \eqref{:K2c}, we set 
\begin{align}\label{:K2g}&
\BRt := - \frac{1}{\sss (\SRover ) } \sum_{k=1}^{\sss (\SRover ) } B_t^k
,\\ \label{:K2h}&
\KRt := -
\int_0^t\frac{1}{\sss (\SRover ) } \sum_{k=1}^{\sss (\SRover ) } \bbb 
( \XRpDu ) du 
,\\\label{:K2i} & 
\LRt := - \int_0^t
 \frac{1}{\sss (\SRover ) } \sum_{k=1}^{\sss (\SRover ) }
 \nR (X_u^{\rR ,\perp , k } ) d L_u^{\rR ,\perp , k }
.\end{align}
Using these notations, we rewrite the first equation in \eqref{:K2c} as 
\begin{align}\notag &
 X_t^{\Rpi } - X_0^{\Rpi } = B_t^i + 
\int_0^t
 \bbb 
( \XRpDui ) du 
\\ & \label{:K2k}\quad \quad 
+ 
\int_0^t
 \nR (X_u^{\Rpi } ) d L_u^{\Rpi }
+ \BRt + \KRt + \LRt ,
\quad 1 \le i \le \sss (\SRover ) 
.\end{align}
The sum of the last three terms comes from $ - \mathbb{U}_{\rrrr }^{\gamma } $ in \eqref{:K2a}, and is called the $ \gamma $-coefficient term. We shall prove that the last three terms vanish as $ \rR \to \infty $.

\begin{lemma} \label{l:K1} For $ \sss = \sum_{i}\delta_{s_i}$, we set 
$ \sssiD =\sum_{j\ne i } \delta_{s_j}$. Let $ 1 \le p < \infty $. 
\begin{align}\label{:K1i}&
\limi{\rR } \frac{\lvert \SRover \rvert }{\sss(\SRover )} = \frac{1}{\rho _{\mathrm{gin}}^1 } =\pi 
\quad \text{ in }L^p (\mug ) \text{ and for } \mu\text{-a.s.\,}\sss 
,\\ \label{:K1a}&
\limi{\rR } 
 \int _{\sSS }
\Big\lvert \frac{1}{\sss (\SRover ) } \sum_{s_i \in \SRover } \bbb ( s_i , \sssiD ) \Big\rvert 
\mug (d \sss )
= 0 
.\end{align}
\end{lemma}
\PF
$ \mug $ is tail trivial \cite{ly.18,o-o.tt}. Hence, $ \mug $ is ergodic under the translation 
$ \{ \vartheta_x \}_{x\in\mathbb{R} ^2 }$, where $ \vartheta_x$ is as in \eqref{:12t}. 
 From the ergodic theorem, we obtain \eqref{:K1i}. 

Let $ \mathbb{I} = [0,1)^2$ and $ F (\sss ) = \sum_{s_i \in \mathbb{I} } \bbb ( s_i , \sssiD ) $. 
Recall that $ \bbb \in \Llocmugonep $, $ 1\le p < 2 $, by \eqref{:K1x} and that 
$ \mug $ is reflection and translation invariant in $ \Rtwo $. 
Hence, we see 
\begin{align}\label{:K1f}&
\int_{\sSS } \lvert F (\sss ) \rvert ^p \mug (d\sss ) < \infty 
,\quad 
\int_{\sSS } F (\sss ) \mug (d\sss ) = 0 
.\end{align}

Note that $ \vartheta_x (\mathbb{I}) \cap \vartheta_y (\mathbb{I}) = \emptyset $ for $ x , y \in \mathbb{Z} ^2 $ such that $ x \ne y $. Let 
\begin{align*}&
I_R = \{ z \in \mathbb{Z} ^2 ; \vartheta _z ( \mathbb{I} ) \subset \SRover \} 
, \ 
J_R = \{ z \in \mathbb{Z} ^2 ; \vartheta _z ( \mathbb{I} ) \not\subset \SRover , 
 \vartheta _z ( \mathbb{I} ) \cap \partial \SRover \ne \emptyset \} 
.\end{align*}
We then have 
\begin{align}\label{:K1g} & 
 \sum_{s_i \in \SRover } 
\bbb ( s_i , \sssiD ) 
= \sum_{ z \in I_R } F ( \vartheta _z ( \sss )) + 
 \sum_{ z \in J_R }
\sum_{s_i \in \vartheta _z ( \mathbb{I} ) \cap \SRover }
 \bbb ( s_i , \sssiD ) 
.\end{align}
For a set $ A $, we denote the volume by $ \lvert A \rvert $. 
Note that 
\begin{align}\label{:K1h}&
\limi{\rR } \frac{\sum_{ z \in I_R } \lvert \vartheta _z ( \mathbb{I} ) \rvert }{\lvert \SRover \rvert } = 1 
,\quad 
\limi{\rR } \frac
{
 \sum_{ z \in J_R }
\lvert \vartheta _z ( \mathbb{I} ) \rvert 
}
{\lvert \SRover \rvert } = 0
.\end{align}
Using the ergodic theorem, we obtain from \eqref{:K1i}, \eqref{:K1g}, and \eqref{:K1h} 
\begin{align} &\notag 
 \limi{\rR } 
 \frac{1}{\sss (\SRover ) } \sum_{s_i \in \SRover } 
\bbb ( s_i , \sssiD ) 
= \pi \int_{\sSS } F (\sss ) \mug (d\sss )
\quad \text{ in } L^1 (\mug )
.\end{align}
Combining this with the equality in \eqref{:K1f}, we obtain \eqref{:K1a}. 
\PFEND

\begin{lemma} \label{l:73} For each $ T \in \mathbb{N} $ and $ \epsilon > 0 $, we have the following. 
\begin{align}
\label{:K2x}&
\limi{\rR } P \Big( \maxT \big\lvert \BRt \big\rvert \ge \epsilon \Big) = 0 
,\\ \label{:K2y}&
\limi{\rR } P \Big( \maxT \lvert \KRt \rvert \ge \epsilon \Big) = 0 
,\\ \label{:K5b}& 
\limi{\qQ }
P \Big( \maxT \Big\lvert 
\frac{1}{\sss (\SQover )} \kSQ 
\int_0^t \bbb ( \XpDuk ) du \Big\rvert \ge \epsilon \Big) = 0 
,\\ \label{:K3b}&
\limi{ Q } 
P \Big( 
\maxT \Big\lvert 
\frac{1}{\sss (\SQover )} 
 \kSQ \Big( X_t^{\perp , k} - X_0^{\perp , k} \Big) \Big\rvert \ge \epsilon \Big) = 0
.\end{align}
\end{lemma}
\PF 
From the ergodic theorem, 
$ \limi{\rR } \sss (\SRover )/ \lvert \SRover \rvert = \rho_{\mathrm{Gin}}^1 = 1/\pi $ 
for $ \mug $-a.s.\,and in $ L ^2 (\mug ) $. 
Recall that $ B ^k $, $ k \in \mathbb{N} $, in \eqref{:K2g} are independent Brownian motions. 
Then, we see that for each $ T >0 $ and $ \epsilon > 0 $ 
\begin{align} &\notag 
\limi{\rR } P \Big( \maxT \big\lvert \BRt \big\rvert \ge \epsilon \Big) = 
\limi{\rR } P \Big( \maxT 
\Big\lvert \frac{1}{\sss (\SRover ) } \sum_{k=1}^{\sss (\SRover ) } B_t^k \Big\rvert \ge \epsilon \Big) = 0 
.\end{align}
This implies \eqref{:K2x}. 

We note that $ \XXRp $ is the diffusion process associated with $ ( \EP _{\rrrr }, \dQbP ) $ on $ \Lmg $ with initial distribution $ \vp \mug $. 
Let $ E_{\mug }^{\Rp }$ be the expectation of $ \XXRp $ with initial distribution $ \mug $ (instead of $ \vp \mug $). 
For each $ T \in \mathbb{N} $ and $ \epsilon > 0 $, 
\begin{align} 
 P \Big( \maxT \lvert \KRt \rvert \ge \epsilon \Big) = 
&\notag 
 P \Big( \maxT \Big\lvert 
\int_0^t\frac{1}{\sss (\SRover ) } \sum_{k=1}^{\sss (\SRover ) } 
\bbb ( \XRpDu ) du 
 \Big\rvert \ge \epsilon \Big) 
\\ \notag \le & 
 P \Big( 
\int_0^T \Big\lvert \frac{1}{\sss (\SRover ) } \sum_{k=1}^{\sss (\SRover ) } 
\bbb 
( \XRpDu ) \Big\rvert du 
 \ge \epsilon \Big) 
\\ \notag \le &
\frac{1}{\epsilon } 
 E \Big[ 
\int_0^T \Big\lvert \frac{1}{\sss (\SRover ) } \sum_{k=1}^{\sss (\SRover ) } 
\bbb 
( \XRpDu ) \Big\rvert du 
 \Big] 
\\ \notag \le &
\frac{\cref{;63}T}{\epsilon } 
 E _{\mug }^{\Rp } \Big[ 
\Big\lvert \frac{1}{\sss (\SRover ) } \sum_{k=1}^{\sss (\SRover ) } 
\bbb 
( \XpDz ) \Big\rvert \Big] 
\\ \notag = &
\frac{\cref{;63}T}{\epsilon } 
 \int _{\sSS }
\Big\lvert \frac{1}{\sss (\SRover ) }
 \sum_{k=1}^{\sss (\SRover ) }
\bbb ( s_k , \ssskD ) \Big\rvert 
\mug (d \sss )
.\end{align}
Here we used the fact that $ \XXRp $ is $ \mug $-symmetric and $ 0 \le \vp \le \cref{;63}$ for the forth line. 
Hence applying \eqref{:K1a} to the last line, we obtain \eqref{:K2y}.

From \lref{l:6V1}\thetag{1}, $ \XX ^{ \perp }$ is a $ \mug $-reversible Markov process 
associated with $ (\EP , \dP ) $ on $ \Lmg $ such that $ \XX _0^{ \perp }\elaw \vp \mug $. 
Hence, we deduce from \eqref{:K1a} 
\begin{align}& \notag 
\limi{\qQ}
E \Big[
\maxT \Big\lvert 
 \int_0^t \frac{1}{\sss (\SQover ) } \kSQ 
 \bbb 
( \XpDu ) 
 du
 \Big\rvert \Big]
\\ \label{:K5i}
 \le &\limi{\qQ} \cref{;63}T 
 \int _{\sSS }
\Big\lvert \frac{1}{\sss (\SQover ) } \kSQ 
 \bbb 
( s_k, \ssskD ) \Big\rvert 
\mug (d \sss ) = 0 
.\end{align}
Using Chebyshev's inequality and \eqref{:K5i}, we obtain \eqref{:K5b}.

From \eqref{:K3x}--\eqref{:K3Y} and \eqref{:K2x}, 
for each $ k \in \mathbb{N} $ such that $ k \le \sss (\SQover ) $, 
\begin{align}\notag 
 X_t^{\perp , k } - X_0^{\perp , k }& = 
\limi{\rR } \Big(  X_t^{\rR ,\perp , k} - X_0^{\rR ,\perp , k} \Big)
=\limi{\rR }  \half \Big( \mathit{M}_t^{[x_i\ot 1 ]} - \mathit{M}_t^{[x_i\ot 1 ]}\circ r_t \Big)
\\ & \label{:K3g} 
= \half \Big( B_t^{ k } - B_{t}^{ k }\circ r_t \Big)
.\end{align}
From \eqref{:K2x} and \eqref{:K3g}, we obtain \eqref{:K3b}. 
\PFEND

\begin{lemma} \label{l:K7}
For each $ T \in \mathbb{N} $ and $ \epsilon > 0 $, the following hold. 
\begin{align}\label{:K7a}&
\limi{\rR } P \Big( \maxT \lvert \LRt \rvert \ge 6 \epsilon \Big) 
= 0 
.\end{align}
\end{lemma}
\PF
Using \eqref{:K2k}, we deduce for each $ 1 \le i \le \sss (\SRover )$ 
\begin{align} \notag 
\LRt = \, &
 X_t^{\Rpi } - X_0^{\Rpi } 
- B_t^{ i } 
- \int_0^t 
 \bbb ( \XRpDui ) du 
\\\label{:K7f} &
 - \int_0^t \nR (X_u^{\Rpi } ) d L_u^{\Rpi }
 - \BRt - \KRt 
.\end{align}
Note that the left-hand side of \eqref{:K7f} is independent of $ 1\le i \le \sss (\SRover ) $. 
We fix $ Q $ and sum over $ 1\le i \le \SQover $. 
Note that $ \BQt = (1/\sss (\SQover ))\sum_{i=1}^{\sss (\SQover )} B_t^i $. 
Then, for any $ Q \le R $, 
\begin{align} \notag 
\LRt &= \frac{1}{\sss (\SQover ) } \Big\{ \iSQ X_t^{\Rpi } - X_0^{\Rpi }
- B_t^{ i } 
- \int_0^t \bbb 
( \XRpDui ) du 
\\ &\notag \quad \quad \quad \quad 
 - \int_0^t \nR (X_u^{\Rpi } ) d L_u^{\Rpi } - \BRt - \KRt \Big\}
\\&\notag 
= \frac{1}{\sss (\SQover ) } \iSQ \Big( X_t^{\Rpi } - X_0^{\Rpi } \Big) -
 \frac{1}{\sss (\SQover ) } \iSQ \int_0^t \bbb ( \XRpDui ) du
\\ \label{:K7h} & \quad \quad \quad \quad 
 - \BQt - \frac{1}{\sss (\SQover ) } \int_0^t \nR (X_u^{\Rpi } ) d L_u^{\Rpi } - \BRt - \KRt 
.\end{align}
From \pref{l:4X}, \lref{l:K6}, and \eqref{:K2y}, the last three terms in \eqref{:K7h} vanish as $ \rR $ goes to infinity. 
Hence, we deduce for each $ Q \in \mathbb{N} $
\begin{align}\notag &%
\limsupi{\rR } P \Big( \maxT \lvert \LRt \rvert \ge 6 \epsilon \Big) 
\le 
P \Big( 
\maxT \Big\lvert 
\frac{1}{\sss (\SQover )}\Big\{
 \iSQ 
X_t^{\perp , i} - X_0^{\perp , i} \Big\} \Big\rvert \ge \epsilon \Big) 
\\& \notag 
+ 
P \Big( 
\maxT \Big\lvert 
\frac{1}{\sss (\SQover )} \iSQ 
\Big\{ 
\int_0^t \bbb ( \XpDui ) du 
\Big\} \Big\rvert 
 \ge \epsilon \Big) 
+
 P \Big( \maxT \big\lvert \BQt \big\rvert \ge \epsilon \Big) 
.\end{align}
Taking $ \qQ \to \infty $ and applying \eqref{:K3b}, \eqref{:K2x}, and \eqref{:K5b}, we obtain \eqref{:K7a}. 
\PFEND

\subsection{Identification of the ISDE for $ \Xp $}\label{s:6V}


Let $ \XRp $ be the solution of \eqref{:K2c}. We take the initial distribution of $ \XRp $ as \eqref{:K2d}. 
In \pref{l:4X}, we proved that $ \XRp $ converge to the continuous process $ \Xp $. 
In \ssref{s:6V}, we identify the ISDE that the limit point $ \Xp $ satisfies. 
Indeed, in \pref{l:6V}, we shall prove that $ \Xp $ satisfies \eqref{:10a}. 

\begin{proposition}\label{l:6V}
$ \Xp $ is a weak solution of ISDE \eqref{:10a} with initial distribution $ \vp \mug \circ \lab ^{-1}$ satisfying
\begin{align}\label{:6Va}&
P ( \ulab (\Xp _t ) \in \cdot ) 
\ll \mu \quad \text{ for all } t > 0 
,\\ \label{:6Vb}&
 P (\IRT (\Xp ) < \infty ) = 1 \quad \text{ for each $ \rR , T \in \mathbb{N}$}
,\\ \label{:6Vc}&
 P (\upath ( \Xp ) \in \WSsiNE ) = 1 
.\end{align}
\end{proposition}
\begin{proof} 
Using \eqref{:K2c} and \eqref{:K2g}--\eqref{:K2i}, we deduce for each $ i \le \m $
\begin{align}	\notag &
X_{t }^{\Rpi } - X_0^{\Rpi } - 
\int_0^{{t }} \bbb (X_u^{\Rpi }, \mathsf{X}_u^{\RpiD }) du 
\\ \notag 
= & B_{t }^{i} + \int_0^{t } \nR (X_u^{\Rpi } ) d L_u^{\Rpi }
+ ( \BR + \KR + \LR )(t ) 
\end{align}
and $ X_{t }^{\Rpi } = X_0^{\Rpi } $, $ 0 \le t < \infty $, for $ i > \m $. 
Let 
\begin{align*}&
\mathcal{R} ^{R,i}(t) = 
\begin{cases}
 \int_0^t \nR (X_u^{\rR ,\perp , i} ) d L_u^{\rR ,\perp , i} + ( \BR + \KR + \LR )(t ) & i \le \m \\
0 & i >\m 
\end{cases}
.\end{align*}
We set $ \mathcal{R} ^R (t) = (\mathcal{R} ^{R,i}(t) )_{i=1}^{\infty}$. 
Let $ \mathbf{B}=(B^i)$ be the $ \RtwoN $-Brownian motion. Then, 
\begin{align}&\label{:47o} 
\XRp (t) - \XRp (0) - \bRp (t) = \mathbf{B}(t) + \mathcal{R} ^R (t) 
,\end{align} 
Using \pref{l:4X}, we have, in law in $ \Cinfty $, 
\begin{align}\label{:47k}&
\limi{\rR} (\XRp (t) - \XRp (0) - \bRp (t) ) = \Xp (t) - \Xp (0) - \bp (t) 
.\end{align}

From \lref{l:K6}, \eqref{:K2x}, \eqref{:K2y}, and \lref{l:K7}, 
for each $ i , T \in \mathbb{N}$ and $ \epsilon > 0 $, 
\begin{align}&\notag 
\limi{\rR } P \Big( \maxT \Big\lvert 
 \int_0^t
 \nR (X_u^{\rR ,\perp , i} ) d L_u^{\rR ,\perp , i} \Big\rvert \ge \epsilon \Big) = 0 
,\\ &\notag 
\limi{\rR } P \Big( \maxT \lvert ( \BR + \KR + \LR )(t ) \rvert \ge \epsilon \Big) = 0 
.\end{align}
Hence, for each $ i , T \in \mathbb{N}$ and $ \epsilon > 0 $, 
\begin{align*}&
\limi{\rR } P \Big( \maxT \Big\lvert \mathcal{R} ^{R,i} (t) \rvert \ge \epsilon \Big) = 0 
.\end{align*}
From this, we obtain 
\begin{align}& \label{:47s}
\limi{\rR } ( \mathbf{B}(t) + \mathcal{R} ^R (t) ) = \mathbf{B}(t) \quad \text{ in law in $ \Cinfty $}
.\end{align}
From \eqref{:47o}--\eqref{:47s}, we obtain in $ \Cinfty $ 
\begin{align}\label{:47t}&
\Xp (t) - \Xp (0) - \bp (t) \elaw \mathbf{B}(t) 
.\end{align}
From \eqref{:K1x}, \eqref{:44x}, and \eqref{:47t}, we see that $ \Xp $ is a weak solution of ISDE \eqref{:10a}. 

By \lref{l:6V1}, the unlabeled dynamics $ \XXp $ is associated with the Dirichlet form $ (\EP , \dP ) $ on $ \Lmg $ and $ \ulab (\Xp _0) \elaw \vp d\mug $. Hence, \eqref{:6Va} is obvious. 

Let $ i \le \m $. From \eqref{:K3x} 
\begin{align} & \notag 
 X_t^{\Rpi } - X_0^{\Rpi }
= \half \Big\{ \mathit{M}_t^{[x_i\ot 1 ]}  - \mathit{M}_t^{[x_i\ot 1 ]}\circ r_t \Big\}
\quad \text{ for } 0 \le t < \infty 
.\end{align} 
Taking $ \rR \to \infty $ and using \eqref{:K2x}, we have for each $ T \in \mathbb{N}$
\begin{align}\label{:47g}&
X_t^{\perp , i } - X_0^{\perp , i } = \half \Big\{ B_t^{i} - B_{t}^i \circ r_t \Big\}
\quad \text{ for } 0 \le t < \infty 
.\end{align}
Here $ r_T $ is the time reversal on $ [0,\infty) $. Using \eqref{:47g}, we obtain \eqref{:6Vb}. 

Let $ \SSs = \{ \sss \in \sSS ; \sss (\{ x \} ) \in \{ 0,1 \} \text{ for all }x \in \Rtwo \} $. 
From \cite[Proposition 5.1]{k-o-t.ifc}, %
\begin{align} \label{:47v}&
P ( \XXp _t \in \SSs \text{ for all } t \in [0,\infty) ) =1 
.\end{align}
From \eqref{:K4y}, we see $ \XXp \in \CiSS $. From this and \eqref{:47v}, we obtain \eqref{:6Vc}. 
\end{proof}

\subsection{Proof of \tref{l:11}}\label{s:65}

In \ssref{s:65}, we complete the proof of \tref{l:11}. 
We begin by quoting a uniqueness result of weak solutions of ISDEs. 

We need the IFC condition \As{IFC} to state the uniqueness result. Because this condition requires further preparation to state, we present \As{IFC} in \ssref{s:73}. 
For the reader's convenience, we give a full version of \lref{l:74'} and related concepts in \ssref{s:73}. 
Let $ \XX = \lpath (\mathbf{X})$ and \\
\As{AC} $ P \circ \XX _t^{-1}$ is absolutely continuous with respect to $\mu $ for all $ 0 < t < \infty $. \\
\As{SIN} $ P ( \XX \in \WSsiNE ) = 1$. \\
\As{NBJ} $ P (\IRT (\mathbf{X}) < \infty ) = 1 $ for each $ R , T \in \mathbb{N}$. \\
Here $ \WSsiNE $ and $ \IRT $ are as in \eqref{:12b} and \eqref{:61w}, respectively. 

\begin{lemma}[{\cite[Corollary 3.2, p.1158]{o-t.tail}}] \label{l:74'}
Assume that $ \mu $ is tail trivial. Then the uniqueness in law of weak solutions $ \X $ of \eqref{:14j} with initial distribution $ \nu = \varphi d\mu \circ \lab ^{-1}$ holds under the constraints of \IASN. 
\end{lemma}

\begin{proposition} \label{l:K8} The Ginibre random point field satisfies \As{A4}. 
\end{proposition}

\begin{proof}

From \pref{l:6V}, $ \Xp $ is a weak solution of \eqref{:10a} with initial distribution $ \nu = \varphi d\mu \circ \lab ^{-1}$ satisfying \ASN. 
Using \cite[Theorems 3.2, 3.3]{k-o-t.ifc}, we see that $ \Xp $ satisfies \As{IFC} for \eqref{:10a}. 
 
Let $ \overline{\mathsf{X}}$ be the diffusion associated with $ \EDg $ on $ \Lmg $. 
Let $ \overline{\mathbf{X}} = \lpath ( \overline{\mathsf{X}}) $ be the associated labeled process with initial distribution $ \nu = \varphi d\mu \circ \lab ^{-1}$. 
Then $ \overline{\mathbf{X}} $ is a weak solution of \eqref{:14j} satisfying \IASN\ \cite{o-t.tail}. 

We have thus seen that $ \Xp $ and $ \overline{\mathbf{X}} $ are weak solutions of \eqref{:10a} with initial distribution $ \nu = \varphi d\mu \circ \lab ^{-1}$ satisfying \IASN. 
It is known that $ \mug $ is tail trivial \cite{o-o.tt,ly.18}. Hence using \lref{l:74'}, we deduce
\begin{align}\label{:K!e}&
\text{$ \X ^{\perp } = \overline{\X } $ in law in $ C([0,\infty );\SN )$}
.\end{align}
From \eqref{:K!e}, \As{NBJ}, and \lref{l:Z3}, we obtain 
 \begin{align}\label{:K!f}&
\text{$ \upath (\X ^{\perp }) = \upath (\overline{\X }) $ in law in $ C([0,\infty );\sSS )$}
\end{align}
for any initial distributions $ \varphi d\mu $ such that $ \vp \in \Lmg $, $ \varphi \ge 0 $, and $ \int \varphi d\mu = 1$. 

From \lref{l:6V1} and \cite{o.isde,o-t.tail}, $ \upath (\X ^{\perp }) $ and $ \upath (\overline{\X }) $ are 
 associated with the Dirichlet forms $ (\EP , \dP ) $ and $ \EDg $ on $ \Lmg $, respectively. 
Thus, we deduce \As{A4} from \eqref{:K!f}. 
\end{proof}

\noindent {\em Proof of \tref{l:11}. } 
It is well known that $ \mug $ satisfies \As{A1} and \As{A2}. 
 \As{A3} was proved in \cite[Theorem 2.3]{o.rm}. 
 \As{A4} follows from \pref{l:K8}.  \As{A6} follows from \lref{l:12}. 
Thus, the Ginibre random point field satisfies all the assumptions of \tref{l:14}, which yields \tref{l:11}. 
\qed 

\section{Appendices}\label{s:VII}
\subsection{Proof of \eqref{:12a}}\label{s:Z}
%
\begin{lemma} \label{l:Z1}
Each $ \ww = \{ \ww _t \} \in \WSsi $ can be written as \eqref{:12a}. 
\end{lemma}
\begin{proof}
Let $ \sS _{ x , \epsilon } = \{ s \in \Rd ; \lvert s - x \rvert < \epsilon \} $. 
Set $ \epsilon (t,x)= \sup\{ \epsilon > 0 ; \ww _{ t } (\sS _{ x , \epsilon } ) =1 \} $. 
Let $ (\tau _0 , x_0) $ be such that $ \ww _{\tau _0} (\{ x_0 \} ) = 1 $. 
Then $ \epsilon_0:=\epsilon (\tau _0 , x_0) > 0 $ because $ \ww _{\tau _0} \in \SSsi $. 
Let 
\begin{align} \notag & 
\tau _1 = \sup \{ t \ge \tau _0 ; \ww _t ( \sS _{ x_0 ,\epsilon_0 /2} ) = 1 \} 
.\end{align}
We see 
$ \ww _{\tau _0} (\partial \sS _{ x_0 ,\epsilon_0 /2}) = 0 $ and 
$ \ww _{\tau _0} (\sS _{ x_0 ,\epsilon_0 /2}) = 1 $. 
Hence from the continuity of $ \ww _t $ at $ \tau _0 $ under the vague topology, we obtain 
\begin{align} & \notag 
\lim_{t\to \tau _0} \ww _{t} (\sS _{ x_0 ,\epsilon_0 /2}) 
= \ww _{\tau _0} (\sS _{ x_0 ,\epsilon_0 /2}) = 1
.\end{align}
We thus obtain $ \tau _1 > \tau _0$ because $ \ww _{t}(\cdot)$ is an integer valued measure.

For each $ t \in [\tau _0 , \tau _1) $, there exists a unique $ w (t) \in \sS _{ x_0 ,\epsilon_0 /2} $ such that 
\begin{align}\label{:Z1h}&
\ww _{t} ( \{ w (t) \} ) = 1
.\end{align}
Thus, $ w (t)$ is a function defined on $ [\tau _0,\tau _1)$ such that 
$ \ww_t (\cdot \cap \sS _{ x_0 ,\epsilon_0 /2} ) = \delta_{ w (t)}$. 

Because $ \ww $ is an $ \sSS $-valued continuous process under the vague topology, 
$ w (t)$ is a continuous function in $ t \in [\tau _0 , \tau _1) $. Indeed, if $ w (t)$ is not continuous at 
$ t_* \in [\tau _0 , \tau _1) $, then there exists a sequence $ \{ t_n \} \subset [\tau _0 , \tau _1) $ such that 
$ \limi{n} t_n = t_* $ and that 
\begin{align*}&
\limi{n} w (t_n) = x_{\infty} \ne w (t_*) 
.\end{align*}
Let $ w (t_*) = x_* $ and $ \epsilon_* = \vert x_{\infty} - x_* \vert $. 
Then $ \epsilon_* < \epsilon_0 $ because $ x_{\infty} , x_* \in \overline{\sS }_{x_0,\epsilon_0/2}$

From $ \epsilon_*/2 < \epsilon_0/2 $ and $ t_* < \tau_1$, we see 
$ \ww _{t_*} (\partial \sS _{x_* , \epsilon_*/2 } ) = 0$. Hence, 
\begin{align}\label{:Z1j}&
\limi{n} \ww _{t_n} (\sS _{x_* , \epsilon_* /2}) = 
\ww _{t_*} ( \sS _{x_* , \epsilon_*/2}) 
= 1
.\end{align}
From $\ww _{t_n} (\{ w (t_n) \} ) = 1$, $w(t_n) \in \sS _{ x_0 ,\epsilon_0 /2} $, and 
$ \ww _{t_n} (\sS _{ x_0 ,\epsilon_0 /2}) = 1$, we see 
\begin{align}\label{:Z1k}&
\limi{n} \ww _{t_n} ( \sS _{x_* , \epsilon_*/2}) = 0
.\end{align}
Combining \eqref{:Z1j} and \eqref{:Z1k} yields contradiction. 
Thus, $ w (t)$ is continuous on $ [\tau _0,\tau _1)$.

Suppose $ \tau _1 < \infty $. Then the limit points of $ w (t)$ as $ t \uparrow \tau _1 $ is single. 
Indeed, if there exist limit points $ z(1) \ne z(2) $, then 
for all sufficiently small $ \epsilon > 0 $, 
\begin{align}& \label{:Z1l}
 \overline{\sS }_{ z(1) , \epsilon } \cap \overline{\sS }_{ z(2) ,\epsilon } = \emptyset , \quad 
\ww _{\tau _1} ( \overline{\sS }_{ z(i) , \epsilon } \backslash \{ z(i) \} ) 
= 0
 \quad \text{ for } i=1,2
.\end{align}
Because $ \overline{\sS }_{ z(1) , \epsilon } $ and $ \overline{\sS }_{ z(2) ,\epsilon }$ are closed sets, we have 
\begin{align} \label{:Z1m} 
& 1 \le 
\limsup_{ t \uparrow \tau _1} \ww _t ( \overline{\sS }_{ z(i) , \epsilon } ) \le 
 \ww _{\tau _1} ( \overline{\sS }_{ z(i) , \epsilon }) 
 \quad \text{ for } i=1,2
.\end{align}
Let $ \SSs = \{ \sss \in \sSS ; \sss (\{ x \} ) \in \{ 0,1 \} \text{ for all }x \in \Rd \} $. 
Because $ \ww _{\tau _1}$ is $ \SSs $-valued and \eqref{:Z1m} holds for all sufficiently small $ \epsilon > 0 $, 
we deduce from \eqref{:Z1l} 
\begin{align} & \label{:Z1n}
\ww _{\tau _1} (\{ z(i) \} )=1 , i = 1,2 
.\end{align}
We take $ \epsilon > 0 $ such that $ \sS _{ z(1) ,\epsilon }\cup \sS _{ z(2) ,\epsilon } \subset \sS _{ x_0 ,\epsilon_0 } $. 
For $ \epsilon > 0 $, we set 
\begin{align} \notag &
\tau _{\epsilon }^* = \inf \{ t \le \tau _1 ; \ww _t ( \sS _{ z(i) ,\epsilon } ) = 1 , i = 1,2 \} 
.\end{align}
Similarly as the proof of $ \tau _1 > \tau _0$, we can prove $ \tau _{\epsilon }^* < \tau _1$ for sufficiently small $ \epsilon > 0 $. 
Using $ \tau _{\epsilon }^* < \tau _1$ and the fact that $ {\sS }_{ z(1) , \epsilon } \cup {\sS }_{ z(2) ,\epsilon }$ is an open set, 
we obtain 
\begin{align}
\label{:Z1q}&
1 \ge 
\liminf_{t \uparrow \tau _1 }
 \ww _t ( {\sS }_{ z(1) , \epsilon } \cup {\sS }_{ z(2) ,\epsilon }) 
\ge 
 \ww _{\tau _1} ( {\sS }_{ z(1) , \epsilon } \cup {\sS }_{ z(2) ,\epsilon }) 
.\end{align}
Combining \eqref{:Z1n} and \eqref{:Z1q} yields contradiction. 
Hence, the function $ w (t)$ on $ [\tau _0,\tau _1)$ can be extended to the continuous function 
on $ [\tau _0,\tau _1]$ such that $ w (\tau _1) = \lim_{t \uparrow \tau _1} w (t)$. 

We set $ x_1 := w (\tau _1) $. Then $ \ww _{\tau _1} (\{ x_1 \} ) = 1$. 
Applying the same method to $ ( \tau _1, x_1 )$, we find $ (\tau _2 , x_2)$ such that 
$ w (t)$ can be extended to $ w \in C([\tau _0 , \tau _2); \Rd ) $ while keeping \eqref{:Z1h}. 
In addition, if $ \tau _2 < \infty $, then we can extend to $ w (t)$ to $ w \in C([\tau _0 , \tau _2]; \Rd ) $ 
and $ x_2 := w (\tau _2)$ satisfies $ \ww _{\tau _2} (\{ x_2 \} ) = 1$. 
Repeating this procedure, we have an increasing sequence $ \{\tau _k \}$, $ k=0,1,2,\ldots $, such that 
 $ w (t) $ can be extended to $ w \in C([\tau _0 ,\tau _{\infty}); \Rd ) $ while keeping \eqref{:Z1h}, 
where $ \tau _{\infty} := \limi{k} \tau _k $. 

For $ (\tau , x )$ such that $ \ww _{\tau } (\{ x \} ) = 1$, we denote by $ (I (\tau , x ) , w^{(\tau , x ) } ) $ the pair of the rectangle $ [\tau _0 , \tau _{\infty})$ and the element $ w $ of $ C([\tau _0, \tau _{\infty}); \Rd )$ constructed as above. 

Two such pairs $ (I (\tau , x ) , w^{(\tau , x )} )$ and $ (I (\tau ', x ') , w^{(\tau ', x ' )})$ satisfy 
that one of each is an extension of another one, or 
that the graphs given by the pairs are disjoint, that is, 
\begin{align*}&
\text{$ \{ ( t , w^{(\tau , x )}(t) ) ; t \in [\tau , \tau _{\infty}) \} \cap 
\{ ( t , w^{(\tau ', x ')}(t) ) ; t \in [\tau ', \tau _{\infty}') \} = \emptyset $}
.\end{align*}
Hence, we set the equivalence relation $ \sim $ such that 
$ (\tau , x ) \sim (\tau ', x ') $ if the former holds. 

For $ (\tau , x )/\sim $, we assign the pair $ (I , w)$ such that $ I=[a,b)$, where 
$ a = \inf\{ \tau ' ; (\tau ',x') \sim (\tau , x ) \} $ and $ b = \tau _{\infty}$. 
Here, if $ (\tau , x )\sim (\tau ', x ') $, then $ \tau _{\infty} = \tau _{\infty}'$.

Note that $ w^{(\tau ', x ')}(\tau ') = x' $ by construction. 
If $ \tau = 0 $ or 
\begin{align*}&
 \liminf_{ \tau ' \to 0} \{\lvert x' \rvert ; (\tau , x )\sim (\tau ', x ') \} < \infty 
,\end{align*}
we retake $ I = [0,b)$ instead of $ I = (0,b)$. 
Let $ w \in C( I ; \Rd )$ be the extension of all $ w^{(\tau ', x ')}(t)$ such that 
$ (\tau ', x ') \sim (\tau , x ) $. Using $ ( I , w )$ thus constructed for each equivalence class, 
we obtain the representation \eqref{:12a}. 
\end{proof}

\begin{lemma} \label{l:Z2}
For a label $ \lab $, we have a unique $ \map{\lpath }{ \WSsiNE }{ C([0,\infty);(\mathbb{R}^d)^{\mathbb{N}}) }$ 
such that $ \lab (\ww _0) = \lpath (\ww )_0 $. 
\end{lemma}
\begin{proof}
For $ \ww \in \WSsiNE $, let $ ( w^i , [0,\infty ))$ be as in \eqref{:12b}. 
Then from \lref{l:Z1}, the set $ \{ ( w^i , [0,\infty )) , i \in \mathbb{N} \} $ is unique. 
Hence, we take the unique bijection $ \map{\sigma }{\mathbb{N}}{\mathbb{N}}$ such that 
$ ( w_0 ^{\sigma (i)})_{i\in \mathbb{N}} = \lab (\ww _0)$. 
We set $ \lpath (\ww ) = ( w ^{\sigma (i)})_{i\in \mathbb{N}}$. 
Because $ w ^{\sigma (i)} \in C([0,\infty); \mathbb{R}^d )$ for all $ i \in \mathbb{N}$ 
and $ (\Rd )^{\mathbb{N}}$ is endowed with the product topology, the map 
$ t \mapsto \lpath (\ww )_t = (w ^{\sigma (i)}(t))_{i\in \mathbb{N}} \in \RdN 
$ is continuous. 
\end{proof}

Let $ \Vert v - w \Vert _T = \sup_{0\le t \le T} \vert v (t) -w (t) \vert $ and 
$ \mathbf{v}=(v^i)_{i=1}^{\infty}$. 
We set 
\begin{align} \notag &
\rho _{\mathrm{path}} ( v , w ) = \sum_{T = 1}^{\infty} \frac{1}{2^{T }} (1 \wedge \Vert v - w \Vert _T )
,\\& \label{:Z3x}
\rho _{\mathrm{lpath}} ( \mathbf{v} , \mathbf{w} )= 
\sum_{i = 1}^{\infty} \frac{1}{2^{i }} (1 \wedge \rho _{\mathrm{path}} ( v^i , w^i ) ) 
.\end{align}
We endow $ \CiRdN $ with $\rho _{\mathrm{lpath}} $. $ \CiRdN $ is a complete separable metric space under $\rho _{\mathrm{lpath}} $. Let $ \rho_{\mathrm{Proho}}$ be the Proholov metric on $ \sSS $: 
\begin{align*}&
 \rho_{\mathrm{Proho}} (\sss , \sss ') = \inf \{ \epsilon > 0 ; \sss ( B ) \le \sss ' (B ^{\epsilon }) + \epsilon , 
 \sss ' (B ) \le \sss (B ^{\epsilon }) + \epsilon 
\} 
.\end{align*}
Here $ B^{\epsilon }$ is the open $ \epsilon $-neighborhood of $ B $. 
Let $ f_R \in C_0(\Rd )$ be such that $ 0 \le f_R \le 1 $, $f_R (x)=1$ on $ \SR $, and $f_{R+1} (x)=0 $ 
on $ \SRR^c$. We set the measure $ f\sss $ by $ f \sss = \sum_i f(s_i) \delta_{s_i}$. 
Then the vague topology on $ \sSS $ is given by the separable, complete metric defined by 
\begin{align}& \label{:Z3y}
\rho_{\mathrm{vague}} (\sss , \sss ') = \sum_{\rR = 1}^{\infty}
\frac{1}{2^{\rR }} \Big( 1 \wedge \rho_{\mathrm{Proho}} (f_R \sss , f_R \sss ') \Big) 
.\end{align}
We refer to \cite{Kal} for these metrics. We set 
\begin{align}\label{:Z3z}&
\rho_{\mathrm{upath}} (\ww , \ww ') = 
 \sum_{ T = 1}^{\infty}\frac{1}{2^{ T }} 
 \Big( 1 \wedge \supT \rho_{\mathrm{vague}} ( \ww _t , \ww _t ') \Big) 
.\end{align}
Then $ \CiSS$ is a complete separable metric space under $ \rho_{\mathrm{upath}}$. 

We set $ \mathbf{w} = (w^i)_{i=1}^{\infty}$ and $ w^i = \{ w^i (t) \} $. 
For $ \mathbf{w} \in \CiRdN $, let 
\begin{align} &\notag 
\IRT (\mathbf{w}) = \sup\{ i \in \mathbb{N}; \min_{t\in[0,T]} \lvert w^i(t)\rvert \le \rR \}
,\\&\notag 
\text{$ \mathcal{NBJ}=\{ \mathbf{w}\in \CiRdN ; 
\text{$ \IRT (\mathbf{w}) < \infty $ for each $ R , T \in \mathbb{N}$} \} $}
.\end{align}
For $ \mathbf{w}=(w^i)$, we set $ \upath (\mathbf{w}) = \mathsf{w}$, where 
$ \mathsf{w} $ is such that $ \mathsf{w}_t = \sum_{i\in\mathbb{N}} \delta_{w^i (t)}$. 
\begin{lemma} \label{l:Z3} 
\thetag{1} $ \upath (\mathbf{w}) \in \CiSS $ for $ \mathbf{w}\in \mathcal{NBJ} $. 
\\\thetag{2} 
$ \map{\upath }{\mathcal{NBJ}}{C([0,\infty);\sSS )}$ is continuous. 
\end{lemma}
\begin{proof}
Let $ \varphi \in C_0(\Rd )$ be such that $ \varphi (x) = 0 $ for $ x \notin \SR $. 
From $ \mathbf{w} \in \mathcal{NBJ} $, we have $ \IRT (\mathbf{w}) < \infty $. 
Let $ t \in [0,T)$. Then, 
\begin{align*}
\lim_{u\to t} \int_{\Rd } \varphi (x) \mathsf{w}_u(dx) &= 
\lim_{u\to t} \sum_{i=1}^{ \infty }\varphi (w^i(u)) = 
\lim_{u\to t} \sum_{i=1}^{ \IRT (\mathbf{w}) }\varphi (w^i(u))
\\& =
\sum_{i=1}^{ \IRT (\mathbf{w}) }\varphi (w^i(t)) =
\sum_{i=1}^{\infty}\varphi (w^i(t)) = \int_{\Rd } \varphi (x) \mathsf{w}_t(dx) 
.\end{align*}
Thus, $ \mathsf{w} = \{ \mathsf{w}_t \} $ is continuous at $ t $, this implies \thetag{1}. 

Suppose $ \limi{n} \rho _{\mathrm{lpath}}( \mathbf{w}_n , \mathbf{w} ) = 0 $. Then 
\begin{align}\label{:Z3e}&
\limi{n} \Vert w_n^i - w^i \Vert _T = 0 
\quad \text{ for all $ i , T \in \mathbb{N}$}
.\end{align} 
Note that 
\begin{align} \notag 
\sup_{0 \le t \le T }& \rho_{\mathrm{vague}} ( \mathsf{w}_n (t) , \mathsf{w} (t)) 
= 
\sup_{0 \le t \le T } \sum_{\rR = 1}^{\infty} \frac{1}{2^R}
\Big( 1 \wedge 
\rho_{\mathrm{Proho}} ( f_R \mathsf{w}_n (t) , f_R \mathsf{w} (t)) \Big)
\\ \le & \notag 
 \sum_{\rR = 1}^{\infty} \frac{1}{2^R}
\Big( 1 \wedge \sup_{0 \le t \le T } 
\rho_{\mathrm{Proho}} ( f_R \mathsf{w} (t) , f_R \mathsf{w} (t)) 
\Big) 
\\ \le & \label{:Z3f}
 \sum_{\rR = 1}^{\infty} \frac{1}{2^R}
\Big( 1 \wedge \sum_{i=1}^{\IRRT (\mathbf{w}_n) \vee \IRRT (\mathbf{w}) }
\Vert w_n^i - w^i \Vert _T 
\Big) 
.\end{align}
From \eqref{:Z3e} and 
$ \mathbf{w}_n , \mathbf{w} \in \mathcal{NBJ}$, we deduce for each $ \rR , T \in \mathbb{N}$
\begin{align}\label{:Z3g}&
\sup_{n} \{\IRRT (\mathbf{w}_n) \vee \IRRT (\mathbf{w}) \} < \infty 
.\end{align}
From \eqref{:Z3e}--\eqref{:Z3g}, 
$ \limi{n} \sup_{0 \le t \le T } \rho_{\mathrm{vague}} ( \mathsf{w}_n (t) , \mathsf{w} (t)) 
= 0 $ for each $ T \in \mathbb{N}$, which implies \thetag{2}. 
%
\end{proof}

\subsection{Strongly local, quasi-regular Dirichlet forms}\label{s:72}

\DN\fFF{\mathsf{F}}
\DN\aAA{\mathsf{A}}
\DN\bBB{\mathsf{B}}
\DN\kKK{\mathsf{K}}
\DN\capamu{\mathrm{Cap}}
\DN\oOO{\mathcal{O} }

In \ssref{s:72}, we recall the concept of strong local, quasi-regular Dirichlet forms following \cite{c-f}. 
%
Let $ \mathcal{O} $ be the family of all open subsets of $ \sSS $. 
For $ \aAA \in \mathcal{O} $, let 
$ \mathcal{L}_{\aAA ,1}=\{ f \in \mathscr{D} ; f \ge 1 \ \mu \text{-a.e.\,on }\aAA \} $. 
We set $ \mathcal{O}_0=\{ \aAA \in \mathcal{O} ; \mathcal{L}_{\aAA ,1} \ne \emptyset \} $. 
Let $ \E _1 = \E + (\cdot,*)_{\Lm }$. 
For an open set $ \aAA \in \mathcal{O} $, we set 
\begin{align*}&
\capamu (\aAA ) = 
\begin{cases}
\inf \{ \E _1(f,f); f \in \mathcal{L}_{\aAA ,1}\}, & \aAA \in \mathcal{O}_0 
\\
\infty , & \aAA \notin \mathcal{O}_0 
.\end{cases}
\end{align*}
For any set $ \bBB \subset \sSS $, we set 
$ \capamu (\bBB ) = \inf \{ \capamu (\aAA ) ; \aAA \in \mathcal{O},\, \aAA \supset \bBB \} 
$. 
The quantity $ \capamu (\bBB )$ is called one-capacity. 

An increasing sequence of closed sets $ \{ \fFF_{k} \} $ is called an $ \E $-nest if, 
for any compact set $ \kKK $, $ \limi{k} \capamu (\kKK \backslash \fFF_{k}) = 0 $. 
%
A function $ f $ is called $ \E $-quasi-continuous if for any $ \epsilon >0 $, there exists an open set $ \oOO $ with $ \capamu (\oOO ) < \epsilon $ such that $ f \vert _{\sSS \backslash \oOO }$ is finite and continuous. 
We call a subset $ \mathfrak{N} \subset \sSS $ an $ \E $-polar set if 
there exists an $ \E $-nest $ \{ \fFF_{k} \} $ such that 
$ \mathfrak{N} \subset \cap_k ( \sSS \backslash \fFF_{k} )$.

A Dirichlet form $ (\E , \mathscr{D} ) $ on $ L^2(\sSS ,\mu ) $ is called quasi-regular if: \\ 
\thetag{1} there exists an $ \E $-nest $ \{ \fFF_{k} , k \ge 1 \} $ consisting of compact sets; 
\\\thetag{2} 
there exists $ \E + ( \cdot , *)_{ L^2(\sSS ,\mu ) }$-dense subset of $ \mathscr{D} $ 
whose elements have $ \E $-quasi-continuous $ \mu $-version; 
\\\thetag{3} 
there exists $ \{ f_k , k \ge 1 \} \subset \mathscr{D} $ having $ \E $-quasi-continuous $ \mu $-versions $ \{ \tilde{f}_k , k \ge 1 \} \subset \mathscr{D} $ and 
an $ \E $-polar set $ \mathfrak{N} \subset \sSS $ such that 
$ \{ \tilde{f}_k , k \ge 1 \} \subset \mathscr{D} $ 
separates the points on $ \sSS \backslash \mathfrak{N} $. 

A Dirichlet form $ (\E , \mathscr{D} ) $ on $ L^2(\sSS ,\mu ) $ is called strongly local if $ \E ( f , g ) = 0 $ for any $ f , g \in \mathscr{D} $ \ such that $ f $ is constant on a neighborhood of the support of $ g $. 
A diffusion process is a family of Markov processes with continuous sample path and has the strong Markov property. We say a diffusion process is conservative if it has an infinite lifetime.

\subsection{IFC condition and the uniqueness of weak solutions of ISDEs}\label{s:73} 

In \ssref{s:73}, we quote the IFC condition and the result of the uniqueness of weak solutions of ISDE from \cite{o-t.tail}. 

Let $ \mathbf{X}= (X ^i)_{i\in\mathbb{N}}$ be an $ \SN $-valued continuous process. 
For $ \mathbf{X}$ and $ i\in\mathbb{N}$, let $ \mathsf{X}$ and $ \mathsf{X}^{\diai } $ as 
$ \mathsf{X}_t = \sum_{i\in\mathbb{N}} \delta_{X_t^i} $ and 
$ \mathsf{X}_t^{\diai } = \sum_{j\in\mathbb{N} ,\ j\not=i } \delta_{X_t^j} $. 

Let $ \SSsde $ be a Borel subset of $ \sSS $ such that $ \SSsde \subset \SSsi $. 
Let $\ulab ^{[1]}$ be the map on $\RdSS $ such that $ \ulab ^{[1]} (x,\sss ) = \delta_x + \sss $. 
Let $ \SSSsde \subset \SN $ and $ \SSSsdeone \subset \RdSS $ be such that 
$ \SSSsde = \ulab ^{-1} (\SSsde ) $ and $ \SSSsdeone =( \ulab ^{[1]}) ^{-1} (\SSsde ) $. 

Let 
$ \map{\sigma }{\SSSsdeone }{\mathbb{R}^{d^2}}$ and $ \map{\bbb }{\SSSsdeone }{\Rd }$ be Borel measurable functions. We consider coefficients 
$ \sigma $ and $ \bbb $ defined only on a suitable subset $ \SSSsdeone $ of $ \RdSS $. 
We consider an ISDE to $ \mathbf{X} =(X^i)_{i\in\mathbb{N}}$ with state space $ \SSSsde $ such that 
\begin{align}\label{:73a}&
dX_t^i = \sigma (X_t^i,\mathsf{X}_t^{\diai }) dB_t^i + 
\bbb (X_t^i,\mathsf{X}_t^{\diai }) dt \  (i\in\mathbb{N}),\quad   
\mathbf{X}_t \in \SSSsde , \forall t \in [0,\infty )
.\end{align}

\begin{definition}[weak solution]\label{d:51} 
By a weak solution of ISDE \eqref{:73a}, we mean an $ \RdN \ts \RdN $-valued stochastic process $ \XB $ 
defined on a probability space $ (\Omega , \mathcal{F}, P )$ 
with a reference family $ \{ \mathcal{F}_t \}_{t \ge 0 } $ such that

\noindent 
\thetag{i} $ \mathbf{X}=(X^i)_{i=1}^{\infty} $ is an $ \{ \mathcal{F}_t \}_{t \ge 0 }$-adapted, $ \SSSsde $-valued continuous process. 
\\
\thetag{ii} $ \mathbf{B} = (B^i)_{i=1}^{\infty}$ is an $ \RdN $-valued 
\FtB 
with $ \mathbf{B}_0 = \mathbf{0}$, 
\\
\thetag{iii} 
the family of measurable $ \{ \mathcal{F}_t \}_{t \ge 0 } $-adapted 
processes $ \Phi $ and $ \Psi $ defined by 
\begin{align}\notag &
\Phi ^i(t,\omega ) = \sigma 
(X_t^i(\omega ),\mathsf{X}_t^{\diai }(\omega ))
,\quad 
\Psi ^i(t,\omega ) = \bbb 
(X_t^i(\omega ),\mathsf{X}_t^{\diai }(\omega ))
\end{align}
belong to $ \mathcal{L}^{2} $ and $ \mathcal{L}^1 $, respectively. 
Here $ \mathcal{L}^{p} $ is the set of all measurable $ \{ \mathcal{F}_t \}_{t \ge 0 } $-adapted 
processes $ \alpha $ such that 
$ E[ \int_0^T \vert \alpha (t,\omega) \vert ^p dt ] < \infty $ for all $ T $. 
Here we can and do take a predictable version of $ \Phi ^i $ and 
$ \Psi ^i$ (see pp 45-46 in \cite{IW}). 
\\
\thetag{iv} with probability one, the process $ \XB $ satisfies for all $ t $ 
\begin{align}\notag &
X_t^i - X_0^i = 
\int_0^t \sigma (X_u^i,\mathsf{X}_u^{\diai }) dB_u^i 
 + 
\int_0^t 
\bbb (X_u^i,\mathsf{X}_u^{\diai }) du \quad (i\in\mathbb{N}) 
.\end{align}
\end{definition}

Let $ \mathbf{X} = (X^i)_{i\in\mathbb{N}}$ be a weak solution to \eqref{:73a} starting at $ \mathbf{s} = \lab (\sss )$. Here $ \map{\lab }{\SSsi }{\SN }$ is the label introduced in \sref{s:1}. 
For $ \mathbf{X}$, we introduce the family of finite-dimensional stochastic differential equations as follows. 

Define $ \map{ \sigmaXms } {[0,\infty) \ts \Rdm }{\mathbb{R}^{d^2}}$ and 
 $ \map{ \bbbXms }{[0,\infty) \ts \Rdm }{\Rd }$ such that, for $ (u,\mathbf{v}) \in \Rdm $ and 
$ \mathsf{v} = \sum_{i=1}^{m-1} \delta_{v_i} \in \sSS $, where 
$\mathbf{v}=(v_1,\ldots,v_{m-1}) \in (\Rd )^{m-1} $, 
\begin{align} &\notag %
 \sigmaXms ( t, (u, \mathbf{v})) = 
 {\sigma } (u , \mathsf{v} + \mathsf{X}_t^{m*}) ,\quad 
\bbbXms ( t, (u, \mathbf{v})) = {\bbb } (u , \mathsf{v} + \mathsf{X}_t^{m*})
.\end{align}
Here $\mathsf{X}_t^{m*} = \sum_{i=m+1}^{\infty} \delta_{X_t^i }$. 
We write $ \lab (\mathsf{s}) =(s_i)_{i\in\mathbb{N}}= \mathbf{s} $ and 
$ \mathsf{s}_m^* = \sum_{i=m+1}^{\infty} \delta_{s_i}$. 
Recall that $ \mathbf{X}_0 = \lab (\mathsf{s}) $. 
We have $ \mathsf{X}_0^{m*} = \mathsf{s}_m^* $ by construction. 
The coefficients $ \sigmaXms $ and $ \bbbXms $ depend on 
both unlabeled path $ \mathsf{X}^{m*} $ and the label $ \lab $. 
Let 
\begin{align} \notag &
\SSSsdemtw 
= \{ \mathbf{s}^m = (s_1,\ldots,s_m) \in \Rdm ; 
\ulab (\mathbf{s}^m) + \ww _t^{m*} \in \SSsde 
 \} 
,\end{align}
where $ \ww _t^{m*}=\sum_{i=m+1}^{\infty}\delta_{w_t^i}$ for 
$ \ww _t = \sum_{i=1}^{\infty} \delta_{w_t^i}$. 
By definition, $ \SSSsdemtw $ is a time-dependent domain in $ \Rdm $ given by $ \ww _t^{m*}$. 

We set 
$ \mathbf{Y}^m=(Y^{m,i})_{i=1}^m $, $ \mathbf{Y}^{m,\diai } = (Y^{m,j})_{j\not=i}^m $, and 
$ \mathsf{Y}_t^{m,\diai }=\sum_{j\not=i}^m \delta_{Y_t^{m,j}}$. 

We consider the stochastic differential equation to $ \mathbf{Y}^m $ with random environment $ \mathsf{X}$ defined on $ \OFpsF $ such that 
\begin{align} &\notag 
dY_t^{m,i} = 
\sigmaXms (t, (Y_t^{m,i},\mathbf{Y}_t^{m,\diai })) dB_t^i + 
 \bbbXms (t, (Y_t^{m,i},\mathbf{Y}_t^{m,\diai })) dt
,\\ &\notag 
 \mathbf{Y}_t^m \in \SSSsdemt \quad \text{ for all } t 
,\\\label{:73h}&
 \mathbf{Y}_0^{m} = \mathbf{s}^m ,
\quad \text{ where $ \mathbf{s}^m=(s_1,\ldots,s_m) $ for $ \mathbf{s}=(s_i)\in\SN $}
.\end{align}
For $ \X = (X^i)_{i\in\mathbb{N}}$, let $ \mathbf{X}^{m*} = (0,\ldots,0,X^{m+1}, X^{m+2},...) $. 
The first $ m $ components of $ \mathbf{X}^{m*} $ are constant paths $ 0 $. A triplet $ (\mathbf{Y}^m,\mathbf{B}^m ,\mathbf{X}^{m*} )$ of continuous processes on $ \OFpsF $ satisfying \eqref{:73h} is called a weak solution.

Let $ \PPPm = \Ps \circ (\lBlhatm )^{-1}  $ be the distribution of $ (\lBlhatm ) $. 
Let $ \WSN = C([0,\infty);\RdN )$ and 
$ \WRdzm = \{ \mathbf{w} \in \WSN ; \mathbf{w}(0) =\mathbf{0} \}$. 
We set 
\begin{align}\notag & 
 \Btm := \Btm (\WRdzm ) = \sigma [\mathbf{w}(s) ; 0\le s \le t ]
,\\&\notag 
 \Bt (\WWdm )= \sigma [(\mathbf{v}(s) ,\mathbf{w}(s) ) ; 0\le s \le t ]
,\\&\notag 
\Ehatmt = \overline{\Bt (\WWdm ) }^{\PPPm }
,\quad 
\Ehatm = 
\overline{\mathcal{B} (\WWdm ) } ^{\Pt ^m }
.\end{align}
\begin{definition}
\label{d:41} 
\thetag{1} 
 $\mathbf{Y} ^{m}$ is called a strong solution of \eqref{:73h} for $ \XB $ \uPs\ if 
$ (\mathbf{Y}^{m},\mathbf{B}^m,\mathbf{X}^{m*})$ satisfies \eqref{:73h} and 
there exists a $ \Ehatm $-measurable function $ \map{\Fms }{\WWdm }{ \WRdzm } $ 
such that $ \Fms $ is $ \Ehatmt /\Btm $-measurable for each $ t $, and $ \Fms $ satisfies %
$ \mathbf{Y}^m = \Fms (\lBlhatm )$ $ \Ps $-a.s.
\\\thetag{2} 
The SDE \eqref{:73h} is said to have a unique strong solution for $ \XB $ \uPs\ 
if there exists a function $ \Fms $ satisfying the conditions in \dref{d:41}\thetag{1} and, 
for any weak solution $ (\hat{\mathbf{Y}}^m ,\mathbf{B}^m,\mathbf{X}^{m*})$ of  \eqref{:73h} \uPs, 
$ \hat{\mathbf{Y}}^m = \Fms (\lBlhatm ) $ for $ \Ps $-a.s. 
\end{definition}

The function $ \Fms $ in \dref{d:41}\thetag{1} is called a strong solution starting at $ \mathbf{s}^m $. 
 \eqref{:73h} is said to have a unique strong solution $ \Fms $ 
if $ \Fms $ satisfies the condition in \dref{d:41}\thetag{2}. 
The function $ \Fms $ is unique for $ \PPPm $-a.s.

We introduce the IFC condition of $ \XB $ defined on $ \OFPF $: 
 

\noindent \As{\iFc} \ 
The SDE \eqref{:73h} 
has a unique strong solution $ \Fms (\mathbf{B}^m,\mathbf{X}^{m*})$ 
for $ \XB $ under $ \Ps $ for $ P\circ \mathbf{X}_0 ^{-1}$-a.s. $ \mathbf{s}$ 
for all $ m \in \mathbb{N}$, where $ \Ps = P (\cdot \vert \mathbf{X} _0 = \mathbf{s})$. 

%

We say that the uniqueness in law of weak solutions for \eqref{:73a} with $ \nu $ holds if 
the laws of the processes $ \mathbf{X}$ and $ \mathbf{X}'$ in $ \WSN = C([0,\infty); \SN )$ coincide for any 
 weak solutions $ \mathbf{X}$ and $ \mathbf{X}'$ with initial distributions $ \nu $. 

\begin{lemma}[Corollary 3.2, p.1158] \label{l:74}
Assume that $ \mu $ is tail trivial. 
Then the uniqueness in law of weak solutions of \eqref{:73a} with initial distribution $ \nu $ holds under the constraints of \As{IFC}, \As{AC}, \As{SIN}, and \As{NBJ}. 
\end{lemma}

\section*{Acknowledgments}
The author thanks Shota Osada for his comment on dual reduced Palm measures and Stuart Jenkinson, PhD, from Edanz (https://jp.edanz.com/ac) for editing a draft of this manuscript. 
This work was supported by JSPS KAKENHI Grant Numbers JP16H06338, JP20K20885, JP21H04432, and JP21K13812. 


\end{document}